%% file: ConwaySpheres.tex
\documentclass{article}
\usepackage[utf8]{inputenc}
\usepackage[english]{babel}
\usepackage{amssymb}
\usepackage{amsthm}
\usepackage{amsmath}
\usepackage{tikz-cd}
\usepackage{setspace}
\usepackage[hidelinks]{hyperref}
\usepackage{cleveref}
\usepackage{booktabs}
\usepackage{float}
\usepackage{longtable}
\usepackage[margin=1.2in]{geometry}

\usepackage{graphicx}
\usepackage{subcaption}
\usepackage{longtable}
\graphicspath{ {images/} }

\DeclareMathOperator{\tr}{tr}

\newcommand{\Ad}{\operatorname{Ad}}
\newcommand{\sltwo}{\mathfrak{sl}_2}

\newcommand{\im}[1]{\mathrm{im} #1}

\theoremstyle{plain}
\newtheorem{theorem}{Theorem}[section]

\theoremstyle{definition}
\newtheorem{definition}[theorem]{Definition}

\theoremstyle{plain}
\newtheorem{prop}[theorem]{Proposition}

\theoremstyle{plain}
\newtheorem{corollary}[theorem]{Corollary}

\theoremstyle{plain}
\newtheorem{lemma}[theorem]{Lemma}

\theoremstyle{remark}
\newtheorem{remark}[theorem]{Remark}

\theoremstyle{remark}

\theoremstyle{plain}
\newtheorem{conjecture}[theorem]{Conjecture}

\theoremstyle{plain}

\newcounter{main}

\theoremstyle{plain}
\newtheorem{mainthm}[main]{\textbf{Theorem}}
\theoremstyle{plain}

\title{Real ideal points, Conway spheres, and left-orderable Dehn fillings}
\author{Yi Wang\thanks{Department of Mathematics, University of Illinois Urbana-Champaign \ Email: \texttt{yiwang20@illinois.edu}}}
\date{}

\begin{document}
	
\maketitle

\begin{abstract}
We study real ideal points of $SL_2(\mathbb C)$-character varieties of knot exteriors containing essential Conway spheres. Under explicit deformation-theoretic hypotheses on the two complementary tangles, we show that essential Conway spheres are detected via Culler-Shalen theory with real ideal points; such ideal points are then deformed into arcs of $SL_2(\mathbb R)$ representations on both sides. We then compute the asymptotic behavior of the representations on these arcs and compute the translation numbers of their lifts to $\widetilde{PSL}_2(\mathbb R)$ representations. Combining this with a result of Gao, we conclude that for certain knots with essential Conway spheres, all sufficiently large positive and negative rational fillings have left-orderable fundamental group. This verifies the $L$-space conjecture for large-slope Dehn fillings of knots which were previously unknown in the literature. We verify the criteria for three infinite families of knots assembled from torus-trivial and twist-trivial tangles, then count the resulting real ideal points, determine the longitudinal translation numbers of all constructed branches, and identify exactly the zero-translation arcs. Computations for the $28$ distinct verified census exteriors listed in Table \ref{tab:instances} agree with the certified branches. 
\end{abstract}

\tableofcontents

\section{Introduction}\label{sec:intro}

\input{intro}

\section{Character varieties, limiting characters, and non-realized gluing data}\label{sec:cullershalen}

\input{cullershalen}

\section{Holonomy extension loci}\label{sec:h00}

\input{h00}

\section{Real ideal points detecting Conway spheres}\label{sec:idealpoints}

\input{smoothness}

\section{Translation numbers of real arcs}\label{sec:trans}

\input{trans}

\section{Torus and twist-trivial tangles}\label{sec:twobridge}

\input{twobridge}

\section{Families of knots with large orderable slopes}\label{sec:families}

\input{families}

\section{Future directions}\label{sec:future}

\input{future}

\appendix

\input{computations}

\bibliography{conway}
\bibliographystyle{plain}

\end{document}

%% file: intro.tex
For a knot $K \subset S^3$, we write its complement as $S^3 \setminus K$, understood as the complement of a tubular neighborhood of the knot in the 3-sphere. A \emph{Conway sphere} for a knot $K \subset S^3$ is an embedded 2-sphere meeting $K$ transversely in four points; it is \emph{essential} in $S^3 \setminus K$ if the corresponding four-holed sphere in the knot exterior is incompressible, boundary-incompressible, and not boundary-parallel. A group is \emph{left-orderable} if it is nontrivial and admits a total order invariant under left multiplication. Finally, $M(r)$ denotes the Dehn filling of $M$ along the slope $r \in \mathbb{Q} \cup \{1/0\}$ in the standard meridian-longitude basis; we say a Dehn filling is left-orderable if its fundamental group is.

\medskip

The \emph{$L$-space conjecture} of Boyer-Gordon-Watson \cite{bgw} predicts that the fundamental group of an irreducible rational homology sphere is left-orderable if and only if the manifold is not an $L$-space. Essential Conway spheres in knot complements provide a natural testing ground for studying various aspects of this conjecture. For instance, Lidman-Moore-Zibrowius \cite{lmz} proved that a knot with an essential Conway sphere has no nontrivial $L$-space surgery. Combining this theorem with the $L$-space conjecture results in the following prediction, which serves as the primary motivation for this paper.

\begin{conjecture}\label{conj:main}
If $M = S^3 \setminus K$ is a knot complement with an essential Conway sphere, $\pi_1(M(r))$ is left-orderable for every $r \in \mathbb{Q}$. 
\end{conjecture}

In this paper, we develop a mechanism which starts from $SL_2(\mathbb R)$ representations on fundamental groups of the two complementary tangles and produces a real ideal point on the character variety detecting the Conway sphere. From this ideal point, we produce real deformation arcs of $SL_2(\mathbb R)$ characters that can be used to study Dehn fillings. As the first applications, we completely implement this mechanism for three infinite families of knots.

\medskip

This approach is a continuation of the program initiated by Culler-Dunfield \cite{cullerdunfield} to prove that Dehn fillings of one-cusped 3-manifolds are left-orderable by exhibiting representations of fundamental groups of such Dehn fillings into $\widetilde{PSL}_2(\mathbb{R})$, particularly synthesizing the ideal point surface detection theory of Culler-Shalen \cite{cullershalen} with the holonomy extension locus techniques of Gao \cite{gao}. The idea of deforming real ideal points detecting well-understood essential surfaces originated in the author's previous work on detecting Seifert surfaces \cite{WangDetectedSeifertSurfaces}. Two new difficulties arise when carrying out this approach for Conway spheres. First, real representation data at ideal points are often degenerate, which results in singularities which obstruct real deformations from such ideal points. Second, even after a deformed real character arc has been constructed, left-orderability requires precise control of longitudinal translation numbers of the associated $\widetilde{PSL}_2(\mathbb R)$ representations. 

\medskip

Our first result (Theorem \ref{mainthm:a}) resolves the former difficulty near a class of $SL_2(\mathbb R)$ representations on tangle exteriors. Under explicit local deformation-theoretic hypotheses on the character varieties of the component tangle exteriors, we prove that the associated singularity is an ordinary real quadratic node. We prove that exactly one of the two associated real branches consists of a genuine deformation arc of $SL_2(\mathbb R)$ representations of $\pi_1(M)$. The degeneracies of our chosen tangle representations force this arc to escape to an ideal point of the character curve which detects the essential Conway sphere in the sense of Culler-Shalen \cite{cullershalen}. 

\medskip

Our second result (Theorem \ref{mainthm:b}) controls the longitudinal translation number of such arcs after lifting them to $\widetilde{PSL}_2(\mathbb R)$ representations. In contrast to the Seifert surface case \cite{WangDetectedSeifertSurfaces}, the homological longitude has nontrivial intersection with an essential Conway sphere, and thus must be broken into pieces lying in each tangle complement. A suitable decomposition of the longitude, combined with a careful analysis of trace asymptotics at the ideal point, reduces this translation number to a finite sum of limiting heights. If this sum is zero, the deformed arcs of $SL_2(\mathbb R)$ representations lift to two arcs with vertical asymptotes on opposite sides in the translation-zero sheet of Gao's holonomy extension locus \cite{gao}. Results from \cite{gao} then imply left-orderability for sufficiently large positive and negative rational Dehn fillings of $M$. 

\medskip

Our third main result (Theorem \ref{mainthm:c}) verifies these two mechanisms completely for explicit tangle models and produces three infinite families of knots for which all sufficiently large positive and negative rational surgeries are left-orderable. We expect that Theorem \ref{mainthm:a} and \ref{mainthm:b} apply to more general families of tangles and knots; we leave this as a direction of future research. 

\medskip

Paoluzzi-Porti \cite{paoluzziporti} construct non-realized limiting gluing data from orbifold holonomy representations, but their construction does not provide the real parabolic limiting representations required for our methods. The limiting representations in this paper instead factor through two-bridge knot quotients; this factoring creates the singularities which are not present in the ideal-point construction of \cite{paoluzziporti}. Our results are also complementary to the results of Culler-Dunfield \cite{cullerdunfield}, Gao \cite{gao}, and the author's previous work \cite{WangDetectedSeifertSurfaces}, whose results typically proved left-orderability in a neighborhood of 0 or a one-sided infinite interval. In contrast, the results of this paper establish left-orderability simultaneously for all sufficiently large positive and negative slopes. To our knowledge, these are the first results of this form through the $\widetilde{PSL}_2(\mathbb R)$-representation program.

\medskip

For the rest of the paper, a \emph{tangle} is a pair $(B^3, t)$ consisting of a 3-ball and two properly embedded disjoint arcs; its \emph{exterior} is the complement of an open tubular neighborhood of $t$ in $B^3$. Throughout, $K \subset S^3$ is a knot with exterior $M = S^3 \setminus K$, and $C$ is an essential Conway sphere for $K$. We write $C' = C \cap M$ for the associated four-holed sphere, and $M = T_1 \cup_{C'}T_2$ for the decomposition of $M$ into the two tangle exteriors it bounds; by a common abuse of notation, we also write $T_1 \cup_CT_2$. Meridians of $K$ and of tangle strands are oriented consistently with a fixed orientation of $K$.

\begin{mainthm}\label{mainthm:a}
Let $M=T_1\cup_C T_2$ be a knot exterior split along an essential Conway sphere. Suppose that there are irreducible representations $\rho_i: \pi_1(T_i) \to SL_2(\mathbb{R})$, $i = 1, 2$, of the distinguished form of Definition \ref{def:admissible}. Suppose that each pair $(T_i, \rho_i)$ satisfies the hypotheses of Theorem \ref{thm:main}. Then the $SL_2(\mathbb C)$ character variety of $M$ has an associated real ideal point detecting $C$, with limiting characters $\chi_1$ and $\chi_2$. At this ideal point, the real normalization consists of two half-arcs of irreducible $SL_2(\mathbb{R})$-characters.
\end{mainthm}

The precise statement appears as Theorem \ref{thm:main}. The tangle hypotheses consist of smoothness assumptions on the two tangle character germs, a local coordinate assumption on each quotient character germ, and two explicit endpoint nondegeneracy conditions. We show that these hypotheses hold for two simple families of tangles: the \emph{torus-trivial} and \emph{twist-trivial} tangles. The simplest tangle satisfying the hypotheses of Theorem \ref{mainthm:a} is the \emph{trefoil-trivial tangle}.

\begin{figure}[h]
	\centering
	\includegraphics[scale=.5]{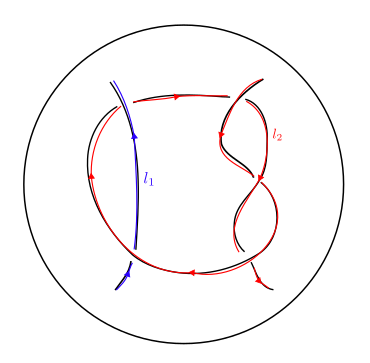}
	\caption{The trefoil-trivial tangle. Notice how the strand with endpoints on the right closes up to a trefoil, and the other strand is trivial. Formally, quotienting the tangle group by the normal closure of the meridian of the trivial strand gives the trefoil knot group, as verified in Section \ref{sec:twobridge}.}\label{fig:trefoiltrivial}
\end{figure}

\begin{mainthm}\label{mainthm:b}
Let $x$ be a real ideal point supplied by Theorem \ref{mainthm:a}. Suppose that the homological longitude admits a four-longitude decomposition which is in positive parabolic normal form along both real half-branches and satisfies height cancellation on both. Then the augmented holonomy images of both half-arcs near $x$ lie in $H_{0,0}(M)$, and their associated Dehn filling slopes tend to $+\infty$ and $-\infty$. Consequently, there exist real numbers $R^-<R^+$ such that $\pi_1(M(r))$ is left-orderable for every rational slope $r\notin [R^-,R^+]$.
\end{mainthm}

The terms \emph{positive parabolic normal form} and \emph{height cancellation} are formally defined in Section \ref{sec:trans}. These two hypotheses have complementary roles. Positive parabolic normal form is a specific sign condition on how the representations of Theorem \ref{mainthm:a} are glued which allows limiting translation numbers of the four longitude pieces to be directly summed. Height cancellation forces this sum to vanish, placing the $SL_2(\mathbb R)$ arcs supplied by Theorem \ref{mainthm:a} into the zero-sheet of Gao's holonomy extension locus. For applications, Section \ref{sec:trans} reduces positive parabolic normal form to a single sign test on independently chosen limiting representatives of the representations on each tangle.

\medskip

For each $n, m \geq 1$, let $T^{tor}_n$ and $T^{tw}_m$ denote the torus-trivial and twist-trivial tangles, constructed in Section \ref{sec:twobridge}. Section \ref{sec:twobridge} also constructs their \emph{flipped}, \emph{mirrored}, and \emph{mirror-flipped} variants and records the relevant data for all of these variants. For $n, m \geq 1$, Section \ref{sec:families} defines the $(n, m)$ \emph{torus-torus}, \emph{twist-torus}, and \emph{twist-twist} knots, each with a base class and a flipped-pair class. Their diagrams are designed so that the four exact local traversals concatenate in the order required by their homological longitude. For these knots, we verify the sign criterion to prove positive parabolic normal form and compute their summed longitudinal translation number directly. 

\begin{mainthm}\label{mainthm:c}
Let $K$ be an $(n, m)$ torus-torus, twist-torus, or twist-twist knot in either the base or the flipped-pair class, and let $M$ be its exterior.
\begin{enumerate}
	\item Theorem \ref{mainthm:a} applies to exactly $nm$, $m$, and $1$ admissible real limiting data respectively, producing that number of distinct real ideal points of $X(M)$ detecting the Conway sphere, each with a two-sided real analytic branch of irreducible $SL_2(\mathbb R)$-characters.
	\item Along every such branch, $\operatorname{tr}\rho_u(\mu) > 2$ while $|\operatorname{tr}\rho_u(\lambda)| \to \infty$; in particular the meridian-normalized lifts define a pair of unbounded augmented holonomy arcs.
	\item Theorem \ref{mainthm:b} applies to exactly $\min(n, m)$, $1$, and $1$ of these branches respectively; the corresponding arcs lie in $H_{0,0}(M)$ and are unbounded with vertical asymptotes on opposite sides. Consequently, there exist $R^- < R^+$ such that $\pi_1(M(r))$ is left-orderable for every rational $r \notin [R^-, R^+]$.
\end{enumerate}
\end{mainthm}

More precisely, the proof of Theorem \ref{mainthm:c} determines the magnitude of the longitudinal translation number on every branch produced by these data. For a torus-torus datum indexed by $(r_{n,j},r_{m,k})$, it is $2|j-k|$. For a twist-torus datum whose torus root is $r_{m,k}$, it is $2k$. For the twist-twist datum it is zero. 

\medskip

The simplest knot to which Theorems \ref{mainthm:b} and \ref{mainthm:c} apply is $K8a14$, obtained by gluing together two trefoil-trivial tangles. 

\begin{figure}[h]
	\centering
	\includegraphics[scale=.5]{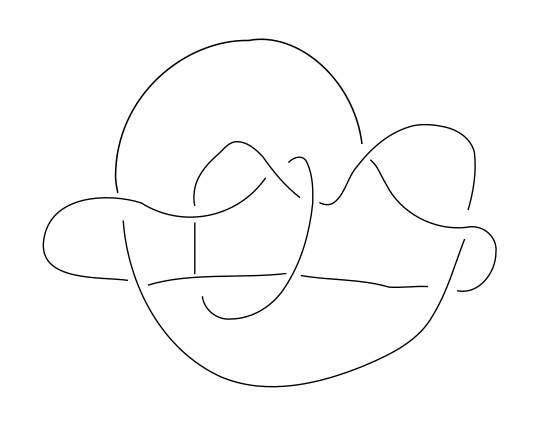}
	\caption{The knot $K8a14$, obtained by gluing together two trefoil-trivial tangles. By Theorem \ref{mainthm:c}, for some compact interval $I = [R^-, R^+]$, every Dehn filling at a rational slope outside $I$ is left-orderable, verifying Conjecture \ref{conj:main} for these fillings.}
\end{figure}

Our results are corroborated by data gathered using the program PE of Culler and Dunfield \cite{CullerDunfieldPE}. For 28 diagrams of torus-torus, twist-torus, and twist-twist knots, with diagrams and identifications performed by SnapPy \cite{SnapPy}, the holonomy extension loci contain the data predicted by Theorem \ref{mainthm:c}. Specifically, there are the predicted numbers of pairs of arcs tending toward slope $\pm\infty$, and the numbers of pairs of arcs with prescribed translation numbers coincide with the above predictions. 

\begin{figure}[h]
	\centering
	\includegraphics[scale=.4]{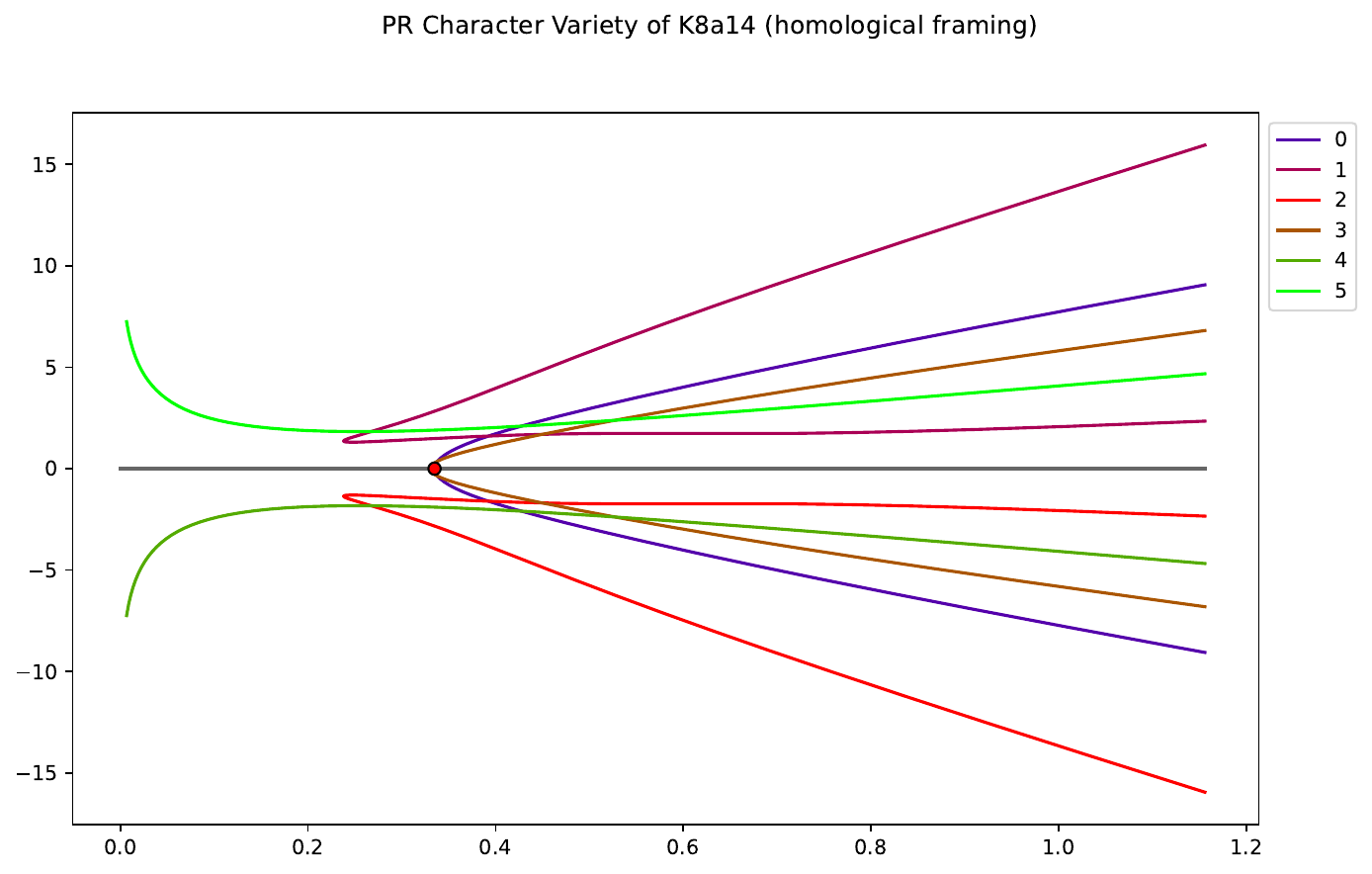}
	\caption{The holonomy extension locus of $K8a14$, as computed by PE \cite{CullerDunfieldPE}. The light- and dark-green arcs approaching slope $\pm\infty$ on the left are the arcs generated by Theorem \ref{mainthm:a}. PE computes the longitudinal translation number at these arcs to be 0, corroborating the results of Theorems \ref{mainthm:b} and \ref{mainthm:c}.}
\end{figure}

We view this paper as the start of a general program aimed at attacking Conjecture \ref{conj:main}. Section \ref{sec:future} describes the first natural questions stemming from the techniques in this paper.

\subsection{Outline of the paper}

Section \ref{sec:cullershalen} recalls Culler--Shalen theory, limiting characters, and the non-realized gluing criterion of Paoluzzi-Porti. Section \ref{sec:h00} reviews translation numbers and Gao's holonomy extension locus, presents the computed loci of the family knots in census range, and compares the four gluings of two trefoil-trivial tangles. In Section \ref{sec:idealpoints} we prove the local node theorem, identify the unique genuine gluing branch, and construct the corresponding real ideal point and two-sided real $SL_2(\mathbb R)$-branch, which completes Theorem \ref{mainthm:a}. Section \ref{sec:trans} introduces four-longitude decompositions and positive parabolic normal form, proves the trace asymptotics, computes the longitudinal translation number from the four limiting heights, and derives the large-slope left-orderability conclusion of Theorem \ref{mainthm:b}. It concludes with a practical sign-and-height criterion stated using chosen limiting representatives. Section \ref{sec:twobridge} certifies the torus-trivial and twist-trivial tangles for applications, computes data to input into the theorems of the previous sections, and transfers these data under mirrors and flips. Section \ref{sec:families} defines the three knot family types, reads the canonical four-longitude decompositions from their diagrams, and proves Theorem \ref{mainthm:c}. Section \ref{sec:future} discusses future directions of research stemming from this paper.

\subsection{Acknowledgements}

The author would like to thank Marc Culler and Nathan Dunfield for providing the PE program to perform extremely informative computations, and for helpful and inspiring discussions about this work. 

\paragraph{AI disclosure statement}

The author used OpenAI's ChatGPT and Anthropic's Claude during the preparation of this manuscript to assist with algebraic calculations, proof verification, literature searches, revisions to exposition, and LaTeX formatting. In particular, ChatGPT was used in working out preliminary free-group reductions in Propositions \ref{prop:torusgroup} and \ref{prop:twistgroup} and intermediate computations in Lemmas \ref{lma:discriminantterms}, \ref{lma:crossjet}, and \ref{lma:commutatorasymptotics}. All mathematical arguments, computations and citations appearing in the final manuscript were checked by the author, who takes full responsibility for the contents. 

%% file: cullershalen.tex
Let $M$ be a compact orientable 3-manifold. Its \emph{$SL_2(\mathbb{C})$ character variety}, denoted $X(M)$, is the affine algebraic set of characters
\begin{equation*}
	\chi_\rho(\gamma) = \mathrm{tr}(\rho(\gamma)) \ \ \ \ \ \rho: \pi_1(M) \to SL_2(\mathbb{C})
\end{equation*}
The natural affine algebraic structure was constructed by Culler-Shalen \cite{cullershalen}.

\subsection{Deformation theory}

We use two standard facts. First, an irreducible $SL_2(\mathbb{C})$-representation is determined by its character up to conjugacy; see Proposition 1.5.2 of \cite{cullershalen}. Second, if $R(\Gamma) = \hom(\Gamma, SL_2(\mathbb{C}))$ is reduced at an irreducible representation $\rho$, then $T_{\chi_\rho}X(\Gamma) \cong H^1(\Gamma; \sltwo(\mathbb{C})_{\Ad\rho})$; see Theorems 53(2), 54(2), and Corollary 55 of \cite{sikora}. Indeed, differentiating an analytic path $\rho_t$ through $\rho$ gives the \emph{cocycle} $u(g) = \frac{d}{dt}|_{t=0}\rho_t(g)\rho(g)^{-1}$, while derivatives of conjugation paths are coboundaries; see Sections 2-3 of \cite{weil}.

\medskip

When $\Gamma$ is free on $x_1, \dots, x_k$, one has $R(\Gamma) = SL_2(\mathbb{C})^k$, and every cocycle integrates: prescribing $\rho_t(x_i) = \exp(tu(x_i))\rho(x_i)$ defines a path of representations with derivative $u$. Every tangle group arising in the applications of this paper is free, so every tangent cocycle appearing in Section \ref{sec:idealpoints} is integrable.

\subsection{Character varieties and limiting characters}

Let $\mathcal C \subset X(M)$ be an irreducible affine curve, let $\overline{\mathcal{C}}$ be its projective closure, and let $\nu: \widetilde{\mathcal C} \to \overline{\mathcal C}$ be the normalization. An \emph{ideal point} of $\mathcal C$ is a point $x \in \widetilde{\mathcal C}$ lying over $\overline{\mathcal C} \setminus \mathcal C$. A properly embedded surface is \emph{incompressible} if it has no compressing disk and no sphere component bounding a ball. It is \emph{essential} if it is incompressible, boundary-incompressible, and has no boundary-parallel component. Culler-Shalen \cite{cullershalen} established a method which associates ideal points on a curve in $X(M)$ to incompressible surfaces in 3-manifolds. The theory is developed in full detail in \cite{cullershalen} and \cite{shalen}. We say that an incompressible surface $S$ associated to an ideal point $x$ of an irreducible curve in $X(M)$ is \emph{detected by $x$}. Classifying which surfaces are detected by ideal points is a well-established problem in 3-manifold topology which has spurred well-known developments such as the \emph{$A$-polynomial}; see \cite{ccgls, chesebrotillmann, ohtsuki}. This paper uses and further develops the methodology of \emph{limiting characters} to study ideal points detecting essential Conway spheres. We describe this methodology here.

\begin{definition}
Let $M$ be a compact orientable 3-manifold, and let $\mathcal{C} \subset X(M)$ be an irreducible algebraic curve. Let $\widetilde{\mathcal{C}}$ be the smooth projective model of $\mathcal{C}$, and let $x \in \widetilde{\mathcal{C}} \setminus \mathcal{C}$ be an ideal point. Let $H \subset \pi_1(M)$ be a finitely generated subgroup, and let $r: X(M) \to X(H)$ be the restriction map. Following Tillmann \cite{tillmann}, we say that $x$ \emph{has a limiting character on $H$} if the induced rational map $\widetilde{\mathcal{C}} \to X(H)$ is defined at $x$. Since the character ring of the finitely generated group $H$ is generated by finitely many trace functions, this is equivalent to requiring that every $I_\gamma$, $\gamma \in H$, have a finite value at $x$. In this case, the value $\chi_{x, H} \in X(H)$ is called the \emph{limiting character of $x$ on $H$}. Equivalently, if $\nu: (\Delta, 0) \to (\widetilde{\mathcal{C}}, x)$ is a local analytic parametrization with $\nu(0) = x$ and $\nu(t) \in \mathcal{C}$ for $t \neq 0$, then the limiting character on $H$ is the character
\begin{equation*}
	\chi_{x, H}(\gamma) = \lim_{t \to 0}I_\gamma(\nu(t)) \ \ \ \ \ \gamma \in H
\end{equation*}
where the limit exists and is finite for all $\gamma \in H$.
\end{definition}

Note that the limiting character can exist even though $\nu(t)$ has no limit in $X(M)$: the traces of elements of $H$ converge while traces of other elements blow up. This is exactly the situation at the ideal points constructed in this paper.

\medskip

By Theorem 2.2.1 and Proposition 2.3.1 of \cite{cullershalen}, the fundamental group of each component of $M$ cut along an associated surface lies in a vertex group of the splitting at $x$, so every trace function on it is pole-free at $x$; hence $x$ has a limiting character on each component. In addition to studying detected surfaces, determining limiting characters on complementary regions of detected surfaces at ideal points has emerged as a research direction in Culler-Shalen theory. Limiting characters have been studied by Dunfield \cite{dunfield}, Tillmann \cite{tillmann}, Paoluzzi-Porti \cite{paoluzziporti}, and the author \cite{WangLimitingCharacters, WangPuncturedJSJTori}. In \cite{WangDetectedSeifertSurfaces}, the author used real limiting characters to construct arcs in the holonomy extension locus and obtain left-orderable fillings. The present paper continues that direction for ideal points detecting Conway spheres.

\subsection{Non-realized gluing data}

We next recall the non-realized gluing mechanism of Paoluzzi-Porti \cite{paoluzziporti} for constructing ideal points with prescribed limiting characters. The first step is to construct a \emph{limiting gluing datum}. Let $M = M_1 \cup_S M_2$, where $S$ is connected, essential, and separating. Let $s_i: X(M_i) \to X(S)$ be the restriction maps. Define the \emph{gluing subvariety}
\begin{equation*}
	V_S(M_1, M_2) = X(M_1) \times_{X(S)}X(M_2) = \{(\chi_1, \chi_2) \mid s_1(\chi_1) = s_2(\chi_2)\}
\end{equation*}
If $r: X(M) \to X(M_1) \times X(M_2)$ is the product restriction map, then $r(X(M)) \subset V_S(M_1, M_2)$. We refer to a point $\xi = (\chi_1, \chi_2) \in V_S(M_1, M_2)$ as a \emph{limiting gluing datum}.

\medskip

Let $\nu_Z: \widetilde Z \to Z$ denote the normalization of an irreducible algebraic curve $Z \subset V_S(M_1, M_2)$ through $\xi$. We say that the limiting gluing datum $\xi$ is \emph{branch-admissible} if there are a curve $Z$ and a point $\widetilde\xi \in \nu_Z^{-1}(\xi)$ for which a punctured analytic neighborhood of $\widetilde\xi$ maps into $r(X(M))$. 

\medskip

We say that $\xi$ is \emph{non-realized} if $\xi \notin r(X(M))$. When the common restriction to $S$ is reducible, such a tuple need not come from an actual representation of $\pi_1(M)$, since equality of boundary characters does not always imply compatibility of boundary representations. When a non-realized tuple is approached along a curve of genuine characters of $\pi_1(M)$, the approaching characters have no limit in $X(M)$: any limit would restrict to the tuple, contradicting non-realization. They therefore converge to an ideal point of the normalization of the character curve, and under a nonconjugacy hypothesis on the boundary restrictions, this ideal point detects the cutting surface. This argument is the content of Lemma 7 of \cite{paoluzziporti}. We restate and prove the version we need for this paper.

\begin{lemma}\label{lma:detection}
Let $\xi = (\chi_1, \chi_2) \in V_S(M_1, M_2)$ be a branch-admissible non-realized limiting gluing datum, where each $\chi_i$ is the character of an irreducible representation $\rho_i: \pi_1(M_i) \to SL_2(\mathbb{C})$. Then there exists an irreducible algebraic curve $\mathcal X \subset X(M)$ and an ideal point of its smooth projective model such that the limiting characters on $\pi_1(M_i)$ are $\chi_i$. If the restrictions $\rho_1|_{\pi_1(S)}$ and $\rho_2|_{\pi_1(S)}$ are not conjugate, then this ideal point detects $S$.
\end{lemma}

\begin{proof}
Let $Z$ be an irreducible algebraic curve representing the branch-admissible germ. By Chevalley's theorem (see Ex. 3.19 of \cite{hartshorne}), $r(X(M)) \cap Z$ is constructible. Since it contains the image of a punctured analytic neighborhood of $\widetilde\xi$, it is infinite and hence Zariski dense in the irreducible curve $Z$. A dense constructible subset of an irreducible variety contains a nonempty Zariski-open subset, so there is a nonempty open $U \subset Z \setminus \{\xi\}$ contained in $r(X(M))$. Since $U$ is contained in the image of $r^{-1}(Z) \to Z$, some irreducible component $Y \subset r^{-1}(Z)$ dominates $Z$. Since $Y \to Z$ is dominant, its generic fiber is nonempty. Choose a closed point of the generic fiber and let $\mathcal X \subset Y$ be its closure. Then $\mathcal X$ is an irreducible algebraic curve and $r|_{\mathcal X}: \mathcal X \to Z$ is dominant. Let $\widetilde{\mathcal{X}}$ and $\widetilde Z$ be their smooth projective models. The dominant map $\mathcal X \to Z$ extends to a nonconstant, hence surjective, morphism $\widetilde r: \widetilde{\mathcal X} \to \widetilde Z$. Let $\widetilde\xi \in \widetilde Z$ be the point over $\xi$ corresponding to the chosen germ, and choose $x \in \widetilde r^{-1}(\widetilde\xi)$. If $x$ lay over an affine point $x_0 \in \mathcal X$, then $r(x_0) = \xi$, contradicting $\xi \notin r(X(M))$. Hence $x$ is an ideal point.

\medskip

For $\gamma \in \pi_1(M_i)$, the trace function $I_\gamma$ on $\mathcal{X}$ is the pullback of the corresponding regular function on $X(M_i)$. Since $\widetilde r(x) = \widetilde\xi$, it follows that $\lim_{\chi \to x}I_\gamma(\chi) = \chi_i(\gamma)$. Thus the limiting character on $\pi_1(M_i)$ is $\chi_i$. Choose a sequence $\chi_n \in \mathcal X$ converging to $x$. Its restrictions to $\pi_1(M_i)$ converge to the irreducible characters $\chi_i$. Under the nonconjugacy hypothesis, Lemma 7 of \cite{paoluzziporti} implies that a subsequence converges to an ideal point detecting $S$. The same subsequence already converges to $x$ in $\widetilde{\mathcal{X}}$, so that ideal point is $x$. 
\end{proof}

The problem is therefore to construct real, branch-admissible, non-realized gluing data. In \cite{paoluzziporti}, such data came from holonomy representations of orbifold components in the Bonahon-Siebenmann decomposition of the $(2, 0)$ Dehn filling. In the author's previous work \cite{WangPuncturedJSJTori, WangLimitingCharacters}, they came from holonomy representations of JSJ components of toroidal Dehn fillings. In both settings, geometry supplied the limiting characters, while non-realization arose because matching boundary characters did not come from compatible boundary representations.

\medskip

The limiting data used here are instead required to be real. At these data, one of the trace-matching equations is first-order degenerate, so the local gluing problem is singular. Section \ref{sec:idealpoints} gives explicit hypotheses under which the pseudo-gluing germ is an ordinary real quadratic node and identifies its unique genuine, hence branch-admissible, branch.

%% file: h00.tex
\subsection{Facts about translation numbers}

We collect here facts about translation numbers. We primarily base our exposition on the relevant sections of Culler-Dunfield \cite{cullerdunfield}; the expository article by Ghys \cite{ghys} also serves as a source. 

\medskip

We use the following conventions from \cite{cullerdunfield}. Let $G = PSL_2(\mathbb{R})$ and $\widetilde{G} = \widetilde{PSL}_2(\mathbb{R})$. We use the successive covering homomorphisms
\begin{equation*}
	\widetilde{G} \xrightarrow{q} SL_2(\mathbb{R}) \xrightarrow{p} PSL_2(\mathbb{R})
\end{equation*}
Let $z$ denote the positive generator of the center of $\widetilde{G}$, normalized by $\mathrm{trans}(z) = 1$. Then $q(z) = -I$, and $\ker q = \langle z^2 \rangle$.

\medskip

For $\widetilde{g} \in \widetilde{G}$, acting on $\mathbb{R}$ as the universal cover of the $G$-action on $\mathbb{RP}^1$, the \emph{translation number} is
\begin{equation*}
	\mathrm{trans}(\widetilde{g}) = \lim_{n \to \infty}\frac{\widetilde{g}^n(x) - x}{n}
\end{equation*}
which exists and is independent of $x \in \mathbb{R}$. We use the following standard properties; see Section 6 of \cite{ghys}. Translation number is invariant under conjugation, satisfies $\mathrm{trans}(\widetilde{g}^n) = n\,\mathrm{trans}(\widetilde{g})$ for all $n \in \mathbb{Z}$ and $\mathrm{trans}(z^k\widetilde{g}) = k + \mathrm{trans}(\widetilde{g})$, and is a quasimorphism with defect one:
\begin{equation*}
	|\mathrm{trans}(\widetilde{g}\widetilde{h}) - \mathrm{trans}(\widetilde{g}) - \mathrm{trans}(\widetilde{h})| \leq 1
\end{equation*}
for all $\widetilde{g}, \widetilde{h} \in \widetilde{G}$. In addition, $\mathrm{trans}: \widetilde G \to \mathbb{R}$ is continuous. Finally, if $\widetilde g$, $\widetilde h \in \widetilde G$ commute, then $\mathrm{trans}(\widetilde g\widetilde h) = \mathrm{trans}(\widetilde g) + \mathrm{trans}(\widetilde h)$. (To see this, combine the defect one property and homogeneity, then divide by $n$ and take the limit as $n \to \infty$.)

\begin{lemma}\label{lma:commutator}
For all $\widetilde{g}, \widetilde{h} \in \widetilde{G}$, we have $|\mathrm{trans}([\widetilde{g}, \widetilde{h}])| \leq 1$.
\end{lemma}

\begin{proof}
By conjugation invariance, $\mathrm{trans}(\widetilde{g}\widetilde{h}\widetilde{g}^{-1}) = \mathrm{trans}(\widetilde{h})$ and $\mathrm{trans}(\widetilde{h}^{-1}) = -\mathrm{trans}(\widetilde{h})$. Applying the defect-one bound to the product $(\widetilde{g}\widetilde{h}\widetilde{g}^{-1})\cdot\widetilde{h}^{-1}$ gives $|\mathrm{trans}([\widetilde{g},\widetilde{h}])| \leq |\mathrm{trans}(\widetilde{h}) - \mathrm{trans}(\widetilde{h})| + 1 = 1$.
\end{proof}

We fix once and for all the following identification, used in every rotation-sense and translation-number computation in this paper (notably in Propositions \ref{prop:torus-heights} and \ref{prop:twist-height}). Identify $\mathbb{RP}^1 = \mathbb{R}/\mathbb{Z}$ via $x \mapsto [\cos \pi x : \sin \pi x]$. Under this identification the positive central generator $z$ acts by $x \mapsto x + 1$, and $\mathrm{trans}$ is the translation number of the induced homeomorphism of $\mathbb{R}$. A hyperbolic or parabolic element of $G$ fixes the points of $\mathbb R / \mathbb Z$ corresponding to its eigenlines, and the rotation $R_\theta = \left(\begin{smallmatrix}\cos\theta & -\sin\theta\\ \sin\theta & \cos\theta\end{smallmatrix}\right)$ with $0 < \theta < \pi$ acts by $x \mapsto x + \theta/\pi$. Reversing this identification negates all translation numbers simultaneously and changes no statement in this paper; what matters is that a single choice is fixed throughout. More generally, for every $\widetilde g \in \widetilde G$, the reduction of $\mathrm{trans}(\widetilde g)$ modulo $\mathbb Z$ is the rotation number of the projected element $p(q(\widetilde g)) \in PSL_2(\mathbb{R})$. In particular, if $\widetilde R_\theta$ is any lift of $R_\theta$, then $\mathrm{trans}(\widetilde R_\theta) \equiv \frac{\theta}{\pi} \pmod{\mathbb Z}$. 

\medskip

The following lemma extends Claim 8.5 of \cite{cullerdunfield} from parabolic elements to hyperbolic elements. 

\begin{lemma}\label{lma:parity}
Let $\widetilde{g} \in \widetilde{PSL}_2(\mathbb{R})$, and let $A = q(\widetilde{g}) \in SL_2(\mathbb{R})$. Suppose $A$ is noncentral hyperbolic or parabolic. Then $\mathrm{trans}(\widetilde{g}) \in \mathbb{Z}$. Moreover, 
\begin{equation*}
	\mathrm{trans}(\widetilde{g}) \equiv 0 \mod 2 \iff \tr(A) > 0
\end{equation*}
Equivalently, $\mathrm{trans}(\widetilde{g}) \equiv 1 \mod 2 \iff \tr(A) < 0$.
\end{lemma}

\begin{proof}
Every hyperbolic or parabolic element of $\widetilde{PSL}_2(\mathbb{R})$ has integral translation number; see Section 3.1 of \cite{cullerdunfield} or \cite{ghys}. If $A$ is parabolic, the parity assertion is precisely Claim 8.5 of \cite{cullerdunfield}. Suppose now that $A$ is hyperbolic. If $\tr(A) > 0$, then, after $SL_2(\mathbb{R})$-conjugation, $A = \left(\begin{smallmatrix}e^s & 0 \\ 0 & e^{-s}\end{smallmatrix}\right)$, $s > 0$. This matrix admits a distinguished lift $\widetilde{A}_0$ with $\mathrm{trans}(\widetilde{A}_0) = 0$. The elements lying over the same $PSL_2(\mathbb{R})$ element are $z^k\widetilde{A}_0$, $k \in \mathbb{Z}$, and satisfy $q(z^k\widetilde{A}_0) = (-1)^kA$, with $\mathrm{trans}(z^k\widetilde{A}_0) = k$. Consequently, the elements whose image under $q$ is $A$ have even translation number, whereas those whose image under $q$ is $-A$ have odd translation number. Since the latter matrices have negative trace, the result follows. 
\end{proof}

\begin{corollary}\label{cor:translationconstant}
Let $I$ be connected and let $u \mapsto \widetilde g_u \in \widetilde G$ be continuous. Suppose that for every $u \in I$, $q(\widetilde g_u)$ is noncentral hyperbolic, noncentral parabolic, or central. Then $u \mapsto \mathrm{trans}(\widetilde g_u)$ is constant.
\end{corollary}

\begin{proof}
In the noncentral hyperbolic and parabolic cases, the translation number is integral by Lemma \ref{lma:parity}. In the central case, $\widetilde g_u$ lies in $\langle z \rangle$, so its translation number is again integral. The result follows because a continuous integer-valued function on a connected space is constant. 
\end{proof}

\begin{corollary}\label{cor:zerotranslift}
Let $A \in SL_2(\mathbb{R})$ be noncentral hyperbolic or parabolic with $\mathrm{tr}(A) > 0$. Then there is a unique $\widetilde A \in \widetilde G$ such that $q(\widetilde A) = A$, $\mathrm{trans}(\widetilde{A}) = 0$. Under the angular convention fixed above, this lift fixes the chosen angular representatives of every eigenline on which $A$ acts by a positive eigenvalue.
\end{corollary}

\begin{proof}
By Lemma \ref{lma:parity}, the translation number of every lift of $A$ through $q$ is an even integer. Multiplication by a unique power of $z^2$ therefore produces a unique lift of translation number zero. If $Av = \lambda v$ with $\lambda > 0$, the lift of the projection action induced by the matrix $A$ fixes the angular representative of $[v]$. It consequently has translation number zero and is the lift just constructed. 
\end{proof}

We have the following standard facts about lifts of representations to $\widetilde{PSL}_2(\mathbb{R})$, stemming from Section 3.4 of \cite{cullerdunfield}. 

\begin{lemma}\label{lma:changinglifts}
Let $\rho: \Gamma \to PSL_2(\mathbb{R})$ be a representation which lifts to $\widetilde{PSL}_2(\mathbb{R})$. If $\widetilde\rho, \widetilde\rho': \Gamma \to \widetilde{PSL}_2(\mathbb{R})$ are two such lifts of $\rho$, then there is a unique homomorphism $\phi \in H^1(\Gamma; \mathbb{Z}) = \mathrm{Hom}(\Gamma, \mathbb{Z})$ such that $\widetilde\rho'(g) = z^{\phi(g)}\widetilde\rho(g)$ for every $g \in \Gamma$. Consequently, 
\begin{equation*}
	\mathrm{trans}(\widetilde\rho'(g)) - \mathrm{trans}(\widetilde\rho(g)) = \phi(g)
\end{equation*}
\end{lemma}

\begin{proof}
This is the parameterization of lifts from Section 3.4 of \cite{cullerdunfield}. Indeed, two lifts of the same $PSL_2(\mathbb R)$ representation differ pointwise by an element of the central subgroup $\langle z\rangle$. Thus there is a function $\phi:\Gamma\to\mathbb Z$ satisfying the displayed equation. Since both lifts are homomorphisms and $z$ is central, $\phi$ is a homomorphism. The translation-number formula follows from $\operatorname{trans}(z^k\widetilde g) = k+\operatorname{trans}(\widetilde g)$. 
\end{proof}

\begin{remark}\label{rmk:evenliftchange}
If $\rho: \Gamma \to SL_2(\mathbb R)$ is fixed and $\widetilde\rho$, $\widetilde\rho'$ both satisfy $q \circ \widetilde\rho = q \circ \widetilde\rho' = \rho$, then the homomorphism of Lemma \ref{lma:changinglifts} is even: there is a unique $\phi \in H^1(\Gamma; \mathbb Z)$ such that $\widetilde\rho'(g) = z^{2\phi(g)}\widetilde\rho(g)$. 
\end{remark}

\begin{corollary}\label{cor:nullhomologous}
If $[g] = 0 \in H_1(\Gamma; \mathbb{Z})$, then $\mathrm{trans}(\widetilde\rho(g))$ is independent of the chosen lift $\widetilde\rho$ of $\rho$. 
\end{corollary}

\begin{proof}
Every $\phi \in \mathrm{Hom}(\Gamma, \mathbb{Z})$ factors through $H_1(\Gamma; \mathbb{Z})$ and hence $\phi(g) = 0$. 
\end{proof}

\begin{lemma}\label{lma:normalizedlift}
Let $M$ be a knot exterior with meridian $\mu$, and let $\rho_u:\pi_1(M)\longrightarrow SL_2(\mathbb R)$, $u>0$ be a continuous family on a one-sided branch, and suppose that $\operatorname{tr}(\rho_u(\mu))>2$ for all sufficiently small $u>0$. Then, after shrinking the branch, there is a unique continuous family of lifts $\widetilde\rho_u: \pi_1(M)\longrightarrow\widetilde{PSL}_2(\mathbb R)$ such that $q\circ\widetilde\rho_u=\rho_u$ and $\operatorname{trans}(\widetilde\rho_u(\mu))=0$. 
\end{lemma}

\begin{proof}
We first construct a continuous family of lifts through $q: \widetilde{PSL}_2(\mathbb{R}) \to SL_2(\mathbb{R})$. Fix a finite generating set $g_1, \dots, g_k$ of $\pi_1(M)$. For each $i$, the path $u \mapsto \rho_u(g_i)$ lifts continuously through the covering $q$, uniquely once a lift is chosen at a single parameter. For each relator $w$ of a presentation on the $g_i$, the corresponding word in the lifted paths is a continuous path in the discrete subgroup $\ker q = \langle z^2 \rangle$, hence constant; choosing the initial lifts at one parameter so that these words are trivial there (this is possible because the obstruction to lifting the $SL_2(\mathbb R)$-representation through its universal cover is the corresponding Euler class in $H^2(\pi_1(M); \mathbb{Z})$, which vanishes for a knot exterior; see \cite{cullerdunfield}) makes them trivial for all $u$. The lifted paths therefore assemble to a continuous family of homomorphisms $\widetilde{\rho}_u$ with $q \circ \widetilde{\rho}_u = \rho_u$. Two lifts through $q$ differ by $z^{2\phi}$, $\phi\in H^1(M;\mathbb Z)$, because $\ker q=\langle z^2\rangle$. By Lemma \ref{lma:parity}, $\operatorname{trans}(\widetilde\rho_u(\mu))$ is an even integer. It is also continuous in $u$, and hence is constant on the connected half-branch. Since $H^1(M;\mathbb Z)\cong\mathbb Z$ is generated by the class dual to the meridian, a fixed central change by $z^{2\phi}$ makes the meridional translation number zero. If two such lifts are meridian-normalized, their difference is $z^{2\phi}$ with $\phi(\mu)=0$. Since $[\mu]$ generates $H_1(M;\mathbb Z)$, this forces $\phi=0$.
\end{proof}

We call $\widetilde\rho_u$ the \emph{meridian-normalized lift}. 

\medskip

Let $H = \left(\begin{smallmatrix}1 & 0 \\ 0 & -1\end{smallmatrix}\right) \in GL_2(\mathbb{R})$ and define $\Omega: SL_2(\mathbb{R}) \to SL_2(\mathbb{R})$ with $\Omega(A) = HAH^{-1}$. Thus, $\Omega\left(\begin{smallmatrix}a & b \\ c & d\end{smallmatrix}\right) = \left(\begin{smallmatrix}a & -b \\ -c & d\end{smallmatrix}\right)$. For $A = \left(\begin{smallmatrix}a & b \\ c & d\end{smallmatrix}\right)$, define  $\mathfrak{d}(A) = b - c$.  We have the following standard facts; these were also used in the results of \cite{WangDetectedSeifertSurfaces}.

\begin{lemma}\label{lma:realformsign}
Let $R_\theta = \left(\begin{smallmatrix}\cos\theta & -\sin\theta \\ \sin\theta & \cos\theta\end{smallmatrix}\right)$.
\begin{enumerate}
	\item If $A = G(\varepsilon \left(\begin{smallmatrix}1&s\\0&1\end{smallmatrix}\right))G^{-1}$ with $\varepsilon = \pm 1$, $s \neq 0$, and $G = \left(\begin{smallmatrix} g_{11} & g_{12} \\ g_{21} & g_{22} \end{smallmatrix}\right) \in GL_2(\mathbb{R})$, then $\mathfrak{d}(A) = \varepsilon s(g_{11}^2 + g_{21}^2)/\det G$. In particular, $\mathfrak{d} \neq 0$ on noncentral parabolic matrices, and two noncentral parabolic matrices of the same trace are $SL_2(\mathbb{R})$-conjugate if and only if their $\mathfrak{d}$-signs agree.
	\item If $A = GR_\theta G^{-1}$ with $0 < \theta < \pi$, then $\mathfrak{d}(A) = -\sin\theta\,(g_{11}^2 + g_{12}^2 + g_{21}^2 + g_{22}^2)/\det G$. In particular $\mathfrak{d} \neq 0$ on elliptic matrices, and two elliptic matrices of the same trace are $SL_2(\mathbb{R})$-conjugate if and only if their $\mathfrak{d}$-signs agree.
	\item The sign of $\mathfrak{d}$ on an elliptic matrix depends only on its rotation direction: if $\operatorname{tr} A = 2\cos\theta_A$ with $0 < \theta_A < \pi$, then $\mathfrak{d}(A) < 0$ if and only if the projective action of $A$ has rotation number $\theta_A/\pi$, and $\mathfrak{d}(A) > 0$ if and only if it has rotation number $-\theta_A/\pi$.
	\item For $\Omega = \mathrm{Ad}_{\operatorname{diag}(1,-1)}$, we have $\mathfrak{d}(\Omega(A)) = -\mathfrak{d}(A)$.
\end{enumerate}
\end{lemma}

\begin{proof}
The formulas in (1) and (2) are direct computations; note in (1) that the ambiguity in the pair $(G, s)$ leaves the sign of $\varepsilon s \det G$, hence of $\mathfrak{d}(A)$, well-defined. For the conjugacy criteria: conjugating by $\operatorname{diag}(1,-1)$ carries $\varepsilon\left(\begin{smallmatrix}1&s\\0&1\end{smallmatrix}\right)$ to $\varepsilon\left(\begin{smallmatrix}1&-s\\0&1\end{smallmatrix}\right)$ and $R_\theta$ to $R_{-\theta}$, so the two $SL_2(\mathbb{R})$-classes of each fixed trace are distinguished by the sign of $\det G$, which the formulas record. For (3): if $\det G > 0$, conjugation preserves the rotation number, which for $R_{\theta_A}$ is $\theta_A/\pi$ by the identification above, and (2) gives $\mathfrak{d}(A) < 0$; if $\det G < 0$, replacing $G$ by $G\operatorname{diag}(1,-1)$ reduces to the previous case with $\theta_A$ replaced by $-\theta_A$, giving rotation number $-\theta_A/\pi$ and $\mathfrak{d}(A) > 0$. Claim (4) is immediate, as $\Omega$ negates both off-diagonal entries.
\end{proof}

The following is another standard translation number fact, shown in the proof of Lemma 6.1 of \cite{cullerdunfield}:

\begin{lemma}\label{lma:outerreverses}
The automorphism $\Omega$ lifts to an automorphism $\widetilde\Omega: \widetilde{PSL}_2(\mathbb R) \to \widetilde{PSL}_2(\mathbb R)$ satisfying $\operatorname{trans}(\widetilde\Omega(\widetilde g)) = -\operatorname{trans}(\widetilde g)$.
\end{lemma}

\subsection{Holonomy extension loci and the method of Gao}

We recall the portion of Gao's holonomy extension locus construction used in this paper. Fix a meridian-longitude basis $\pi_1(\partial M)=\langle\mu,\lambda\rangle$. An \emph{augmented representation} is a pair $(\widetilde\rho,v)$, where $\widetilde\rho:\pi_1(M)\longrightarrow\widetilde G$ is a representation whose peripheral image is hyperbolic, parabolic, or central, and $v\in\mathbb RP^1$ is a common fixed point of the projected peripheral representation. For $\gamma\in\pi_1(\partial M)$, let $A_\gamma\in SL_2(\mathbb R)$ be the image of $\widetilde\rho(\gamma)$ under the intermediate covering $q: \widetilde G\longrightarrow SL_2(\mathbb R)$. Let $a_\gamma\in\mathbb R^\times$ be the eigenvalue of $A_\gamma$ on the line $v$. Following Definitions 3.1 and 3.3 of \cite{gao}, define the \emph{eigenvalue-translation map} 
\begin{equation*}
	\operatorname{EV}(\widetilde\rho,v)(\gamma) = \left(\log|a_\gamma|, \operatorname{trans}(\widetilde\rho(\gamma))\right)
\end{equation*}
Since the peripheral subgroup is abelian, both components are homomorphisms on $\pi_1(\partial M)$. Thus $\operatorname{EV}(\widetilde\rho,v) \in H^1(\partial M;\mathbb R) \times H^1(\partial M;\mathbb Z)$. 

\begin{definition}
The \emph{holonomy extension locus} of $M$, denoted $HL_{\widetilde{G}}(M)$, is the closure in $H^1(\partial M; \mathbb{R}) \times H^1(\partial M; \mathbb{Z})$ of the image of the augmented representations under $\operatorname{EV}$.
\end{definition}

Let $\mu^*,\lambda^*$ be the basis of $H^1(\partial M)$ dual to $\mu,\lambda$. We write $\operatorname{EV}(\widetilde\rho,v) = \left(x\mu^*+y\lambda^*, i\mu^*+j\lambda^*\right)$ as $(x,y;i,j)$. For $i,j\in\mathbb Z$, define the corresponding sheet by
\begin{equation*}
	H_{i,j}(M) = \left\{(x,y)\in\mathbb R^2: (x,y;i,j)\in HL_{\widetilde G}(M)\right\}
\end{equation*}
In particular, $H_{0,0}(M)$ is the sheet on which both the meridian and longitude have translation number zero. If $u\longmapsto(\widetilde\rho_u,v_u)$ is a continuous family of augmented representations, we call its image under $\operatorname{EV}$ an \emph{augmented holonomy arc}. Thus, on the sheet $H_{0,0}(M)$, its coordinates are
\begin{equation*}
	x(u)=\log|a_\mu(u)| \ \ \ \ \ y(u)=\log|a_\lambda(u)|
\end{equation*}
By Corollary \ref{cor:translationconstant}, the translation coordinates are constant along connected families of augmented representations. Consequently, every augmented holonomy arc is contained in a single sheet $H_{i,j}(M)$. For a rational slope $r=\frac pq$, written in the basis $(\mu,\lambda)$ so that the filling slope is $\mu^p\lambda^q$, let $L_r = \left\{ (x,y)\in\mathbb R^2: px+qy=0 \right\}$. Thus, when $q\neq0$, $L_r$ is the line through the origin having ordinary slope $-\frac pq=-r$. Equivalently, at a point with $x\neq0$, $r=-\frac yx$. Following Definition 3.4 of \cite{gao}, a point of the holonomy extension locus is an \emph{ideal point of the locus} if it lies in the closure of the image of the eigenvalue-translation map but not in the image itself.

\medskip

We will also use the following description of the asymptotic behavior of unbounded arcs.

\begin{lemma}[Lemma 3.6 of {\cite{gao}}]\label{lma:gaoasymptotic}	
Suppose that, for some $i,j$, the sheet $H_{i,j}(M)$ contains an unbounded arc. Then the arc approaches an asymptote $y = -rx$ in $\mathbb{R}^2$, where $r$ is the boundary slope of an incompressible surface associated to some ideal point of $X(M)$. 
\end{lemma}

In our applications, $r$ will be the meridional slope $1/0$, and the unbounded arcs in $H_{0,0}$ approach the vertical line $x = 0$.

\begin{remark}\label{rmk:gaogenerality}
Although the results of \cite{gao} are stated with a view toward homology spheres, Lemma 3.8 of \cite{gao} applies to an arbitrary rational filling slope: its hypotheses are that the slope line meets $H_{0,0}(M)$ at a nonzero point that is not an ideal point of the locus, and that the filled manifold is irreducible. We use it in this generality in Corollary \ref{cor:slopes}.
\end{remark}

\subsection{Computed holonomy extension loci}

The results of this paper originated in the computations presented in this subsection. We list the 28 distinct census exteriors obtained from the verified small-parameter family presentations of Table \ref{tab:instances}. For each of them, Table \ref{tab:instances} compares the exact ideal-point counts and absolute longitudinal translation numbers determined later in the paper with the PE computations. We also display the four gluing knots of a pair of trefoil-trivial tangles as a compact illustration of the general behavior. 

\medskip

All loci displayed in this section were computed with PE \cite{CullerDunfieldPE}. Each knot was identified by constructing the family diagram, computing its exterior in SnapPy \cite{SnapPy}, and verifying an isometry with the census exterior. Census names do not distinguish a knot from its mirror; since a global mirror negates every signed translation number, all comparisons are between multisets of absolute values. Translation numbers are those of the meridian-normalized lift of Lemma \ref{lma:normalizedlift}; each is a constant integer along its arc. Several exteriors admit more than one verified family presentation; Table \ref{tab:instances} is indexed by gluing, so such knots appear in more than one cell, and Remark \ref{rmk:priority} fixes the primary presentation. For every knot whose locus we display, the computed vertical-asymptote arcs, their count, and their translation numbers agree with the values that Proposition \ref{prop:matching-template} and Lemma \ref{lma:trans} determine, and no additional vertical-asymptote arcs were found in the computed output.

\medskip

\begin{figure}[h]
	\centering
	\begin{subfigure}{0.48\textwidth}
		\includegraphics[scale=.4]{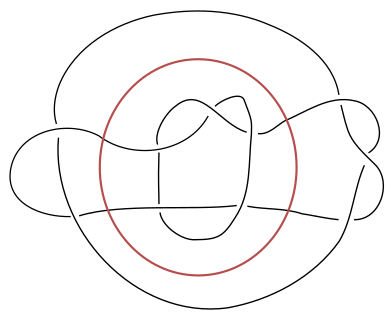}
		\caption{$K8a14$}\label{fig:K8a14}
	\end{subfigure}\hfill
	\begin{subfigure}{0.48\textwidth}
		\includegraphics[scale=.4]{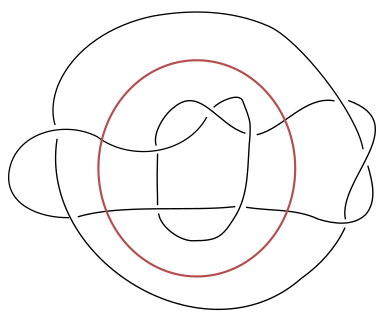}
		\caption{$K10n11$}\label{fig:K10n11}
	\end{subfigure}
	\medskip
	\begin{subfigure}{0.48\textwidth}
		\includegraphics[scale=.4]{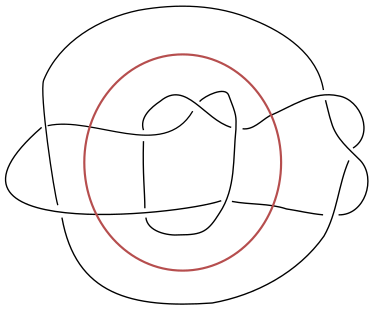}
		\caption{$K10n12$}\label{fig:K10n12}
	\end{subfigure}\hfill
	\begin{subfigure}{0.48\textwidth}
		\includegraphics[scale=.4]{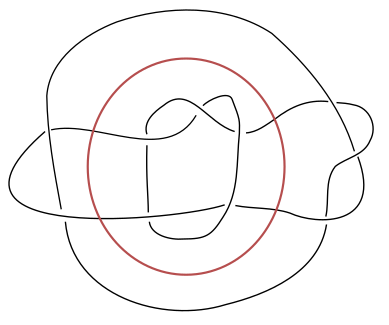}
		\caption{$K8a15$}\label{fig:K8a15}
	\end{subfigure}
	\caption{Four ways of gluing together two trefoil-trivial tangles. As shown in Figure \ref{fig:fourclasses}, the orderability results of this paper will only apply to $K8a14$.}
	\label{fig:fourclassesknots}
\end{figure}

\begin{figure}[H]
	\centering
	\begin{subfigure}{0.48\textwidth}
		\includegraphics[width=\linewidth]{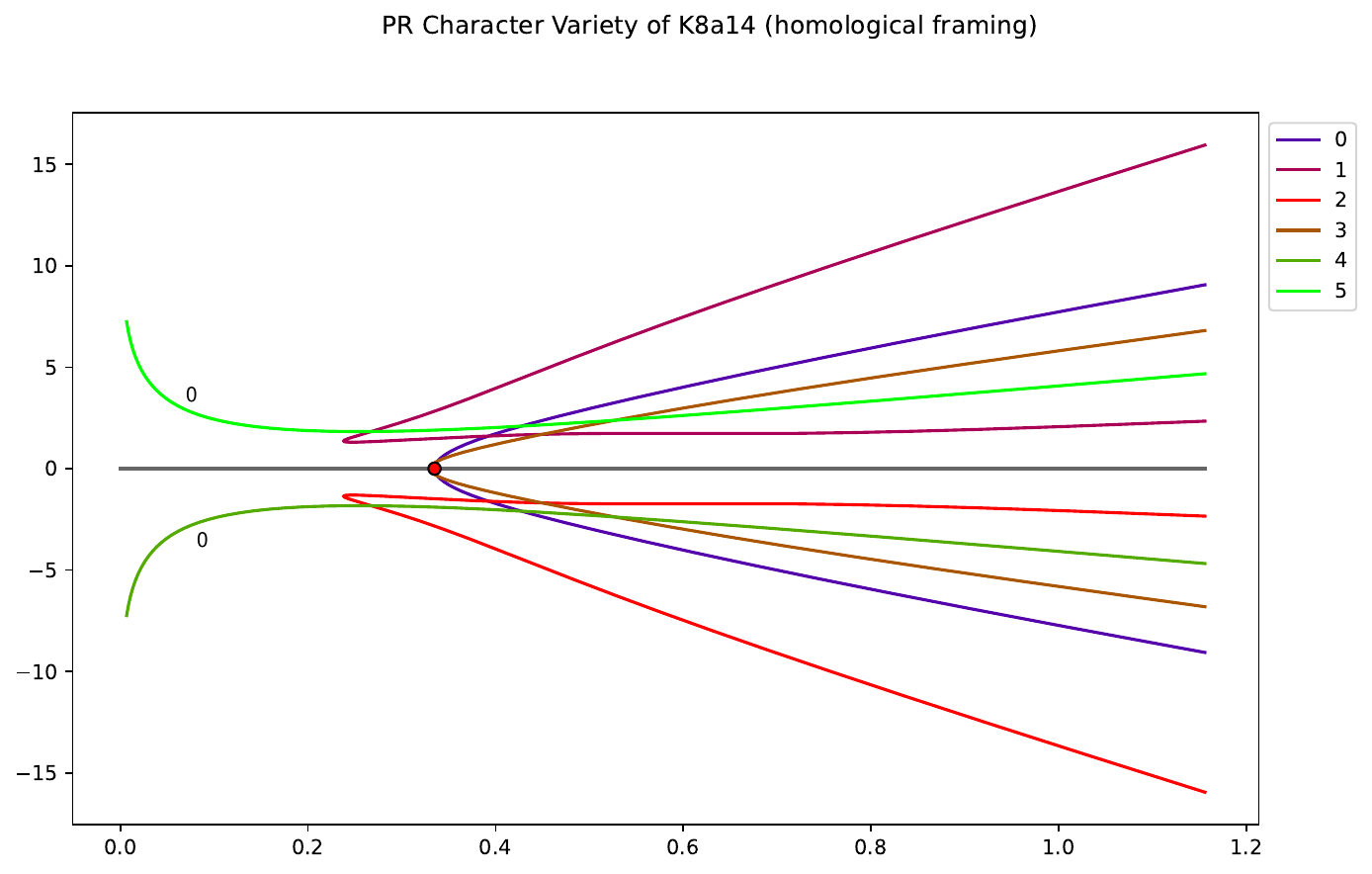}
		\caption{}\label{fig:fourclasses-a}
	\end{subfigure}\hfill
	\begin{subfigure}{0.48\textwidth}
		\includegraphics[width=\linewidth]{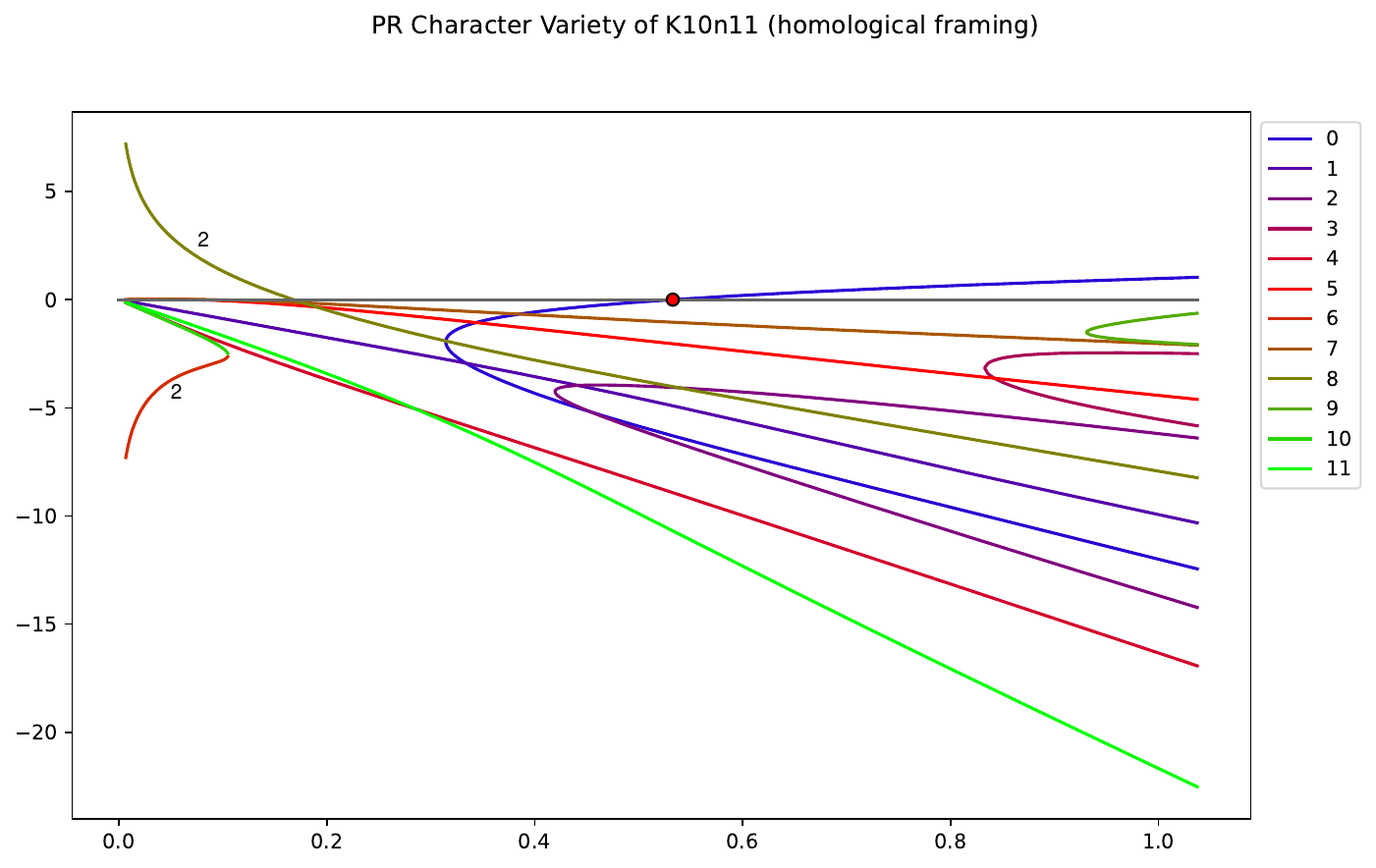}
		\caption{}\label{fig:fourclasses-b}
	\end{subfigure}
	\medskip
	\begin{subfigure}{0.48\textwidth}
		\includegraphics[width=\linewidth]{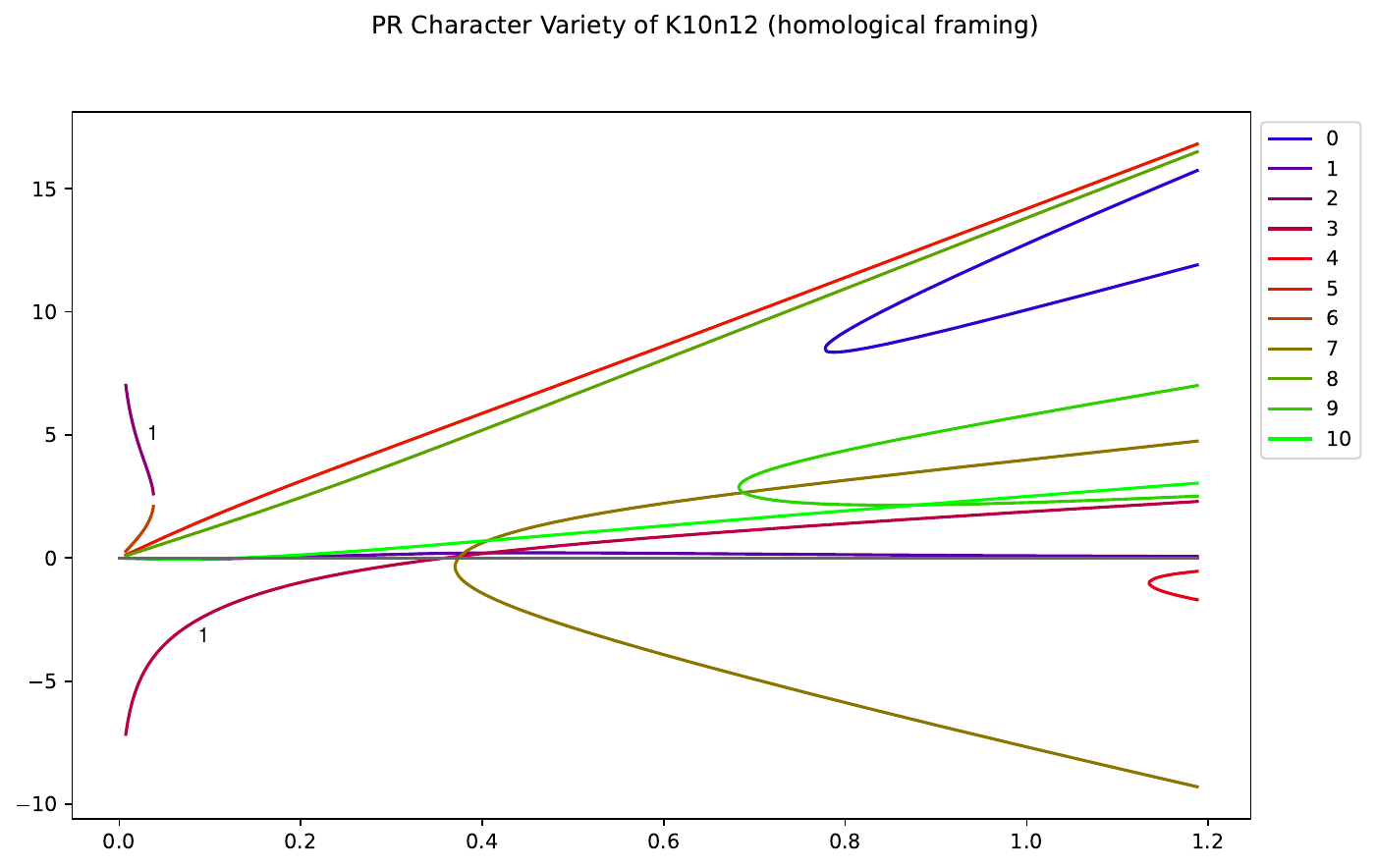}
		\caption{}\label{fig:fourclasses-c}
	\end{subfigure}\hfill
	\begin{subfigure}{0.48\textwidth}
		\includegraphics[width=\linewidth]{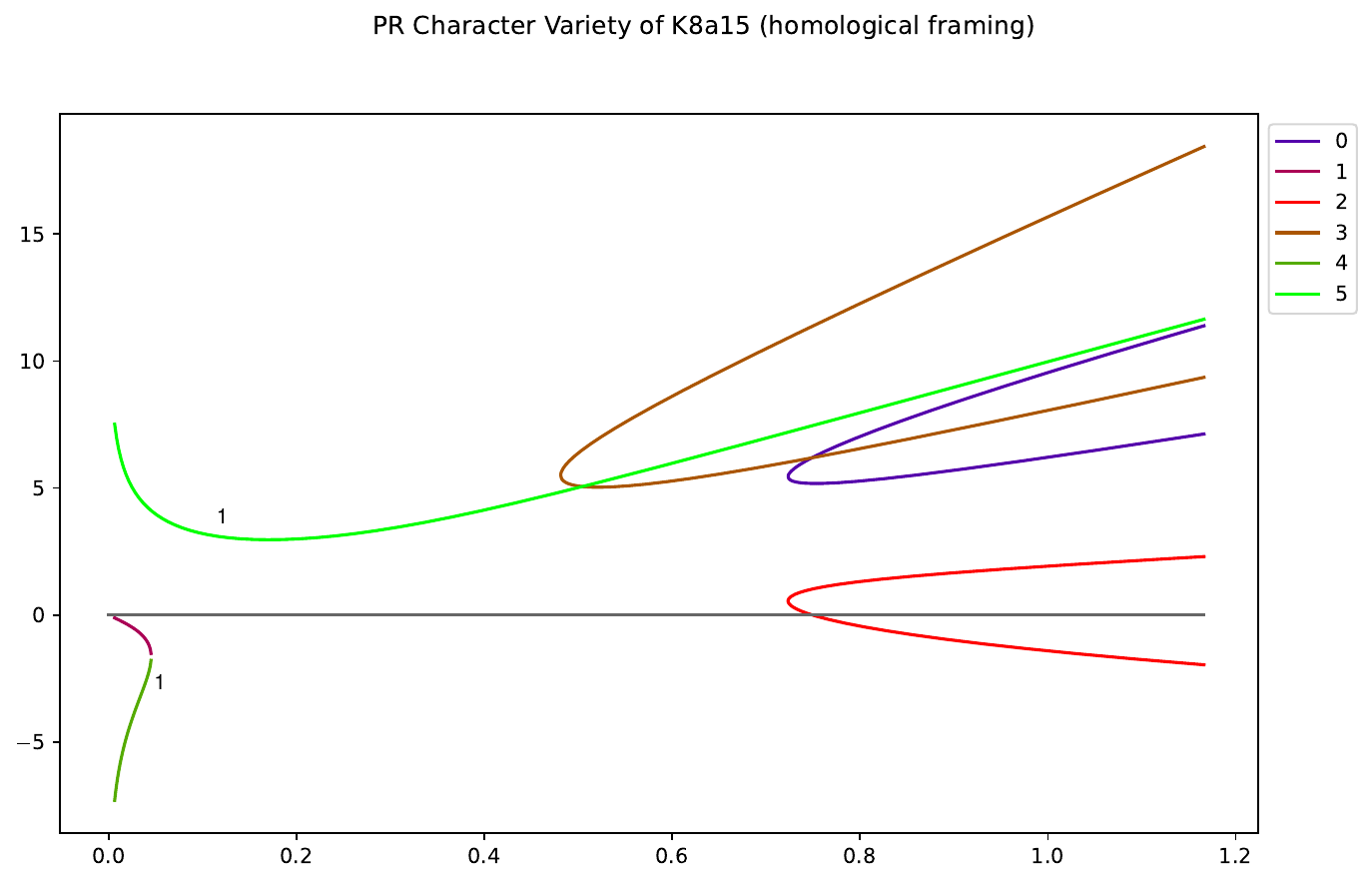}
		\caption{}\label{fig:fourclasses-d}
	\end{subfigure}
	\caption{The computed loci of the four knot-forming gluing classes of a pair of trefoil-trivial tangles. Only in (a) do such arcs lie in $H_{0,0}$; these produce the left-orderable fillings of Theorem \ref{mainthm:c}.} \label{fig:fourclasses}
\end{figure}

There are four ways to glue together two trefoil-trivial tangles, shown in Figure \ref{fig:fourclassesknots}; the knots are $K8a14$, $K10n11$, $K10n12$, and $K8a15$. Studying their holonomy extension loci highlights the main questions this paper addresses. Theorem \ref{mainthm:a} applies to all four of these knots, since the hypotheses of the theorem are only based on the tangle exteriors. However, the hypotheses of Theorems \ref{mainthm:b} and \ref{mainthm:c} apply to only $K8a14$. As can be seen in Figure \ref{fig:fourclasses}, all four knots admit a pair of opposite-signed arcs with vertical asymptotes, as predicted by Theorem \ref{mainthm:a}. However, the longitudinal translation numbers of the arcs are 0 for $K8a14$, 2 for $K10n11$, and 1 for both $K10n12$ and $K8a15$. 

\medskip

Theorems \ref{mainthm:a}-\ref{mainthm:c} prove the ideal-point counts and the orderability result for three families of knots which include $K8a14$. These families are depicted in Figure \ref{fig:families}; they are called the \emph{$(n,m)$ torus-torus, twist-torus, and twist-twist knots}. For these knots, Proposition \ref{prop:matching-template} and Lemma \ref{lma:trans} also determine the longitudinal translation numbers of all of the augmented holonomy arcs stemming from Theorem \ref{mainthm:a}. The PE computations in Table \ref{tab:instances} agree with these values in every row, and no additional vertical-asymptote arcs were found in the computed output.

\begin{corollary}\label{cor:instances}
Let $K$ be any of the 28 knots of Table \ref{tab:instances} and $M$ its exterior. There exist $R^- < R^+$ such that $\pi_1(M(r))$ is left-orderable for every rational slope $r \notin [R^-, R^+]$. 
\end{corollary}

\begin{proof}
Each knot belongs to one of the three families in Definition \ref{def:families}; apply Theorem \ref{mainthm:c}. 
\end{proof}

\begin{figure}[h]
	\centering
	\includegraphics[width=0.7\textwidth]{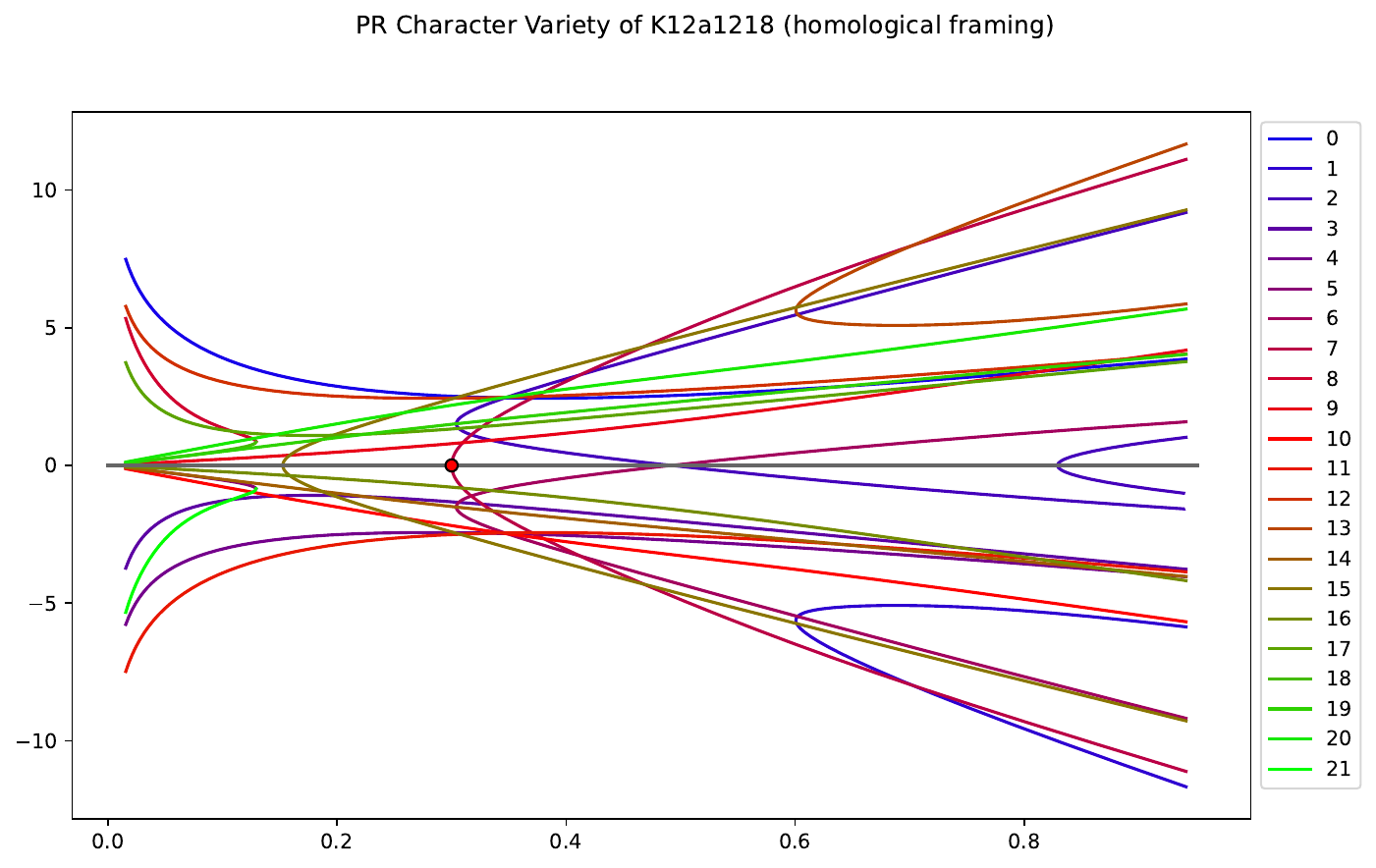}
	\caption{The computed locus of the $(2, 2)$ torus-torus knot, which is $K12a1218$. Arcs unbounded with vertical asymptotes are labeled by their longitudinal translation numbers: four pairs, with longitudinal translation numbers $\{0,0,2,2\}$. The two pairs with translation number $0$ lie in $H_{0,0}$, with exactly $\min\{2,2\} = 2$ of the four, as predicted by Theorem \ref{mainthm:c}.}
	\label{fig:selectivity}
\end{figure}

The $(2, 2)$ torus-torus knot $K12a1218$ (Figure \ref{fig:selectivity}) demonstrates the full extent of what the results in this paper prove: four pairs of arcs coming from ideal points detecting the essential Conway sphere, of which exactly $\min\{2, 2\} = 2$ lie in $H_{0,0}$. 

\begin{longtable}{c|cc}
	\caption{Census identifications of family knots at small parameters. Each entry is the name, in the Hoste-Thistlethwaite census as distributed with SnapPy, of the knot obtained from the indicated gluing of Definition \ref{def:families}. Every identification is a SnapPy verified isometry of the constructed diagram against the census exterior. Census names do not distinguish mirrors. The constructed exteriors at the blank entries admit no verified isometry with any census exterior. The named entries yield 28 distinct verified census exteriors.}
	\label{tab:instances} \\
		\toprule
		$(n, m)$ & base class & flipped-pair class \\
		\midrule
		\endfirsthead
		\toprule
		$(n, m)$ & base class & flipped-pair class \\
		\midrule
		\endhead
		\midrule
		\multicolumn{3}{r}{\emph{continued on next page}}
		\endfoot
		\bottomrule
		\endlastfoot
		\multicolumn{3}{l}{\emph{torus-torus}} \\
		$(1,1)$ & K8a14 & K10a78 \\
		$(2,1)$ & K10a83 & K12a815 \\
		$(3,1)$ & K12a850 & K14a12838 \\
		$(4,1)$ & K14a13425 & \\
		$(2,2)$ & K12a1218 & K14a18218 \\
		$(3,2)$ & K14a18259 & \\
		\midrule
		\multicolumn{3}{l}{\emph{twist-torus}} \\
		$(1,1)$ & K10n10 & K8a14 \\
		$(2,1)$ & K12n437 & K10a39 \\
		$(3,1)$ & K14n16405 & K12a640 \\
		$(4,1)$ & & K14a11421 \\
		$(1,2)$ & K12n115 & K10a83 \\
		$(2,2)$ & K14n9964 & K12a365 \\
		$(3,2)$ & & K14a8859 \\
		$(1,3)$ & K14n1970 & K12a850 \\
		$(2,3)$ & & K14a5197 \\
		\midrule
		\multicolumn{3}{l}{\emph{twist-twist}} \\
		$(1,1)$ & K10a7 & K8a14 \\
		$(2,1)$ & K12a46 & K10a39 \\
		$(3,1)$ & K14a498 & K12a640 \\
		$(2,2)$ & K14a4128 & K12a273 \\
		$(3,2)$ & & K14a3441 \\
		\bottomrule
\end{longtable}

\begin{remark}\label{rmk:priority}
Several knots admit more than one family presentation: K8a14 admits three (torus-torus base $(1,1)$, and the flipped-pair $(1,1)$ gluings of both other types), and K10a39, K12a640, K10a83, and K12a850 admit two each. Each presentation carries its own certified branches, and Theorem \ref{mainthm:c} applies to each independently. When we refer to a knot by a single presentation, we use the first available in the order: torus-torus before twist-torus before twist-twist; base class before flipped-pair; then lexicographically least $(n, m)$. Under this convention K8a14 is the $(1,1)$ torus-torus knot, as in Figure \ref{fig:fourlongitude}.
\end{remark}

%% file: smoothness.tex
We prove Theorem \ref{mainthm:a} in this section. Let $K$ be a knot with exterior $M = S^3 \setminus N(K)$, and let $C$ be an essential Conway sphere giving a decomposition $M = T_1 \cup_C T_2$. Write $C' = C \cap M$ for the corresponding four-holed sphere. Since $C'$ is separating, van Kampen's theorem gives $\pi_1(M) \cong \pi_1(T_1) *_{\pi_1(C')}\pi_1(T_2)$ with respect to the basepoints fixed below.

\medskip

The construction proceeds in four steps, which structure this section. 
\begin{enumerate}
	\item Restriction to the tangle exteriors maps the relevant characters of $\pi_1(M)$ into the fiber product of the two tangle character varieties over the character variety of the four-holed sphere; near the limiting datum this is the \emph{gluing variety} $V_C(T_1, T_2)$.
	\item Define the \emph{pseudo-gluing variety} $V'_C(T_1, T_2)$ using coordinates on the character variety of $C'$; this contains $V_C(T_1, T_2)$ as a subvariety. The regularity and nondegeneracy hypotheses force the germ of $V'_C(T_1, T_2)$ at the limiting character $\xi$ to be an ordinary real quadratic node.
	\item Exactly one of the two branches lies in the full gluing variety $V_C(T_1, T_2)$, and its generic points are realized by irreducible characters of $\pi_1(M)$, while the limiting datum itself is non-realized.
	\item Realization along the branch together with non-realization at its limit produces a real ideal point of the character curve detecting $C$, together with a two-sided real analytic branch of irreducible $SL_2(\mathbb{R})$-characters.
\end{enumerate} 

\subsection{Boundary and gluing conventions}

In this paper, oriented Wirtinger generators create a \emph{left-crossing} with the oriented strand. This results in the \emph{right-handed} Wirtinger convention: if $a$ is the in-undercrossing, $a'$ the out-undercrossing and $z$ is the overcrossing, then $a' = zaz^{-1}$. Accordingly, the strand following the under-crossing, say $t$, is the inverse of $z$. Our traversal convention is thus $a' = t^{-1}at$. 

\begin{figure}[h]
	\centering
	\includegraphics[scale=.7]{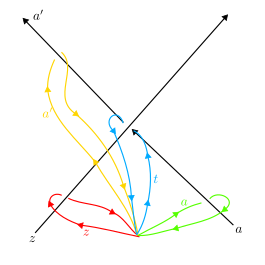}
	\caption{The right-handed Wirtinger convention. The over-arc is $z$, the incoming and outgoing under-arcs are $a$ and $a'$, and the oriented strand traversal is $t=z^{-1}$. Thus $a' = zaz^{-1}$ and $a' = t^{-1}at$.}
\end{figure}

Let $C_i' = C \cap T_i \subset \partial T_i$, $i = 1, 2$, denote the two copies of the boundary Conway spheres in the two tangles. Let $\phi: C_1' \to C_2'$ be the gluing homeomorphism. Choose a basepoint $*$ in the interior of $C_1'$, and use $\phi(*)$ as the basepoint of $C_2'$. We use these basepoints for $\pi_1(C_i')$, $\pi_1(T_i)$, and after gluing, $\pi_1(M)$. 

\medskip

Orient the knot $K$ and the embedded four-holed sphere $C'$. All puncture loops below are oriented using this fixed orientation of $C'$. In particular, the loops on the two boundary copies are not oriented independently using the boundary orientations inherited from $T_1$ and $T_2$. We have
\begin{equation}\label{eq:boundaryconventions1}
	\pi_1(C_1') = \langle \mu_{1, 1}, \mu_{1, 2}, \mu_{1, 3}, \mu_{1, 4} \mid \mu_{1, 1}\mu_{1, 2}\mu_{1, 3}\mu_{1, 4} = 1 \rangle
\end{equation}
\begin{equation}\label{eq:boundaryconventions2}
	\pi_1(C_2') = \langle \mu_{2, 4}, \mu_{2, 1}, \mu_{2, 2}, \mu_{2, 3} \mid \mu_{2, 4}\mu_{2, 1}\mu_{2, 2}\mu_{2, 3} = 1 \rangle
\end{equation}
For notational purposes, we fix the map 
\begin{equation*}
	\phi_*(\mu_{1,1}) = \mu_{2,4} \ \ \ \ \ \phi_*(\mu_{1,2}) = \mu_{2,1} \ \ \ \ \ \phi_*(\mu_{1, 3}) = \mu_{2,2} \ \ \ \ \ \phi_*(\mu_{1, 4}) = \mu_{2,3}
\end{equation*}
For each tangle, let $m_{i, 1}, \dots, m_{i, 4} \in \pi_1(T_i, *)$ be based Wirtinger meridians oriented consistently with the corresponding oriented tangle strands. From now on, we use the following convention. The symbols $m_{i, j}$ always denote same-oriented strand meridians. Strand conjugacy relations, quotient groups, and the hypotheses imposed on a tangle representation are written in the $m$-variables. Four-holed-sphere trace coordinates and all gluing equations are written in the $\mu$-variables. To summarize, we include a sample image of the $K8a14$ knot, with all $m_{i, j}$ and $\mu_{i,j}$ labeled with arrows.
\begin{figure}[h]
	\centering
	\includegraphics[scale=1]{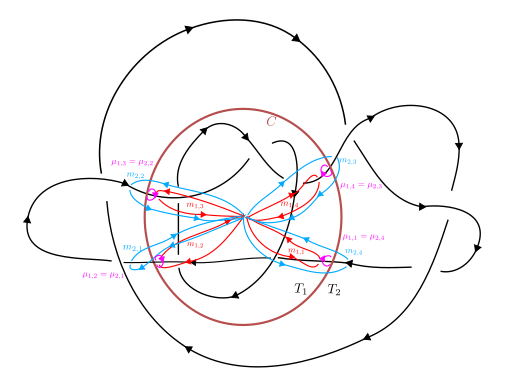}
	\caption{The boundary-oriented puncture loops $\mu_{i,j}$ and the same-oriented strand meridians $m_{i,j}$ for the knot $K8a14$. The arrows record the global knot orientation, and the labels show the crossed boundary identification $\mu_{1,j} \sim \mu_{2,j-1}$ used throughout the paper, where $j-1$ is read modulo $4$.}\label{fig:orientation}
\end{figure}
We then have the equivalences 
\begin{equation*}
	\mu_{1, 1} = m_{1, 1}^{-1} \ \ \ \ \ \mu_{1, 2} = m_{1, 2} \ \ \ \ \ \mu_{1, 3} = m_{1, 3}^{-1} \ \ \ \ \ \mu_{1, 4} = m_{1, 4}
\end{equation*}
\begin{equation*}
	\mu_{2, 1} = m_{2, 1} \ \ \ \ \ \mu_{2, 2} = m_{2, 2}^{-1} \ \ \ \ \ \mu_{2, 3} = m_{2, 3} \ \ \ \ \ \mu_{2, 4} = m_{2, 4}^{-1}
\end{equation*}
Crucially, we label the punctures so that the chosen orientation of $K$ encounters them in the order $\mu_{1, 1} \to \mu_{1, 2} \to \mu_{1, 3} \to \mu_{1, 4}$. These identifications with this chosen orientation will be used throughout the remainder of this paper.

\subsection{Admissible limiting gluing data}

By Lemma \ref{lma:detection}, it remains to prove that the admissible limiting datum is branch-admissible. We will consider the following special forms of limiting data.

\begin{definition}\label{def:admissible}
An \emph{admissible real limiting datum} is a limiting gluing datum $\xi = (\chi_1, \chi_2)$ represented by irreducible homomorphisms $\rho_i: \pi_1(T_i) \to SL_2(\mathbb{R})$ such that, in the same-oriented strand variables,
\begin{equation}\label{eq:admissible1}
	\rho_i(m_{i, 1}) = \rho_i(m_{i, 2}) = I \ \ \ \ \ \rho_i(m_{i, 3}) = \rho_i(m_{i, 4}) = P_i
\end{equation}
where $P_i$ is a noncentral parabolic of trace 2. Equivalently, in boundary variables, 
\begin{equation}\label{eq:admissible2} \begin{array}{c|cccc} &\mu_{i,1}&\mu_{i,2}&\mu_{i,3}&\mu_{i,4}\\ \hline T_1&I&I&P_1^{-1}&P_1\\ T_2&I&I&P_2&P_2^{-1}. \end{array} \end{equation}
\end{definition}

Informally, an admissible real limiting datum kills the strand with endpoint meridians $m_{i, 1}$ and $m_{i, 2}$, while it acts as a noncentral parabolic on the meridian of the other strand. We refer to the former strand as the \emph{identity strand}, and the latter as the \emph{parabolic strand}. We refer to the meridian of the parabolic strand, i.e. $m_{i, 3}$, as the \emph{surviving meridian}. Every word in the boundary subgroup has trace 2 at such a datum, so the two boundary characters agree and $\xi \in V_C(T_1, T_2)$. 

\begin{lemma}\label{lma:unrealized}
An admissible real limiting datum is not in the image of $X(M) \to X(T_1) \times X(T_2)$. Moreover, the two limiting boundary representations are not conjugate.
\end{lemma}

\begin{proof}
Suppose a representation of $\pi_1(M)$ realizes $\xi$. Since the two tangle characters are irreducible, its restrictions would be conjugate to $\rho_1$ and $\rho_2$. According to our gluing convention, we have $\mu_{1, 3} = \mu_{2, 2}$. At the limiting datum, we have $\rho_1(\mu_{1, 3}) = P_1^{-1}$ and $\rho_2(\mu_{2, 2}) = I$. These matrices cannot be conjugate, so the two boundary representations are not conjugate. Thus, the datum cannot be realized by a representation of $\pi_1(M)$.
\end{proof}

The following is Lemma 5 of \cite{paoluzziporti}.

\begin{lemma}\label{lma:irreduciblegluing}
Let $(\eta_1,\eta_2)\in V_C(T_1,T_2)$. If the common restriction $\eta_1|_{\pi_1(C')} = \eta_2|_{\pi_1(C')}$ is an irreducible character, then there exists a unique character $\eta\in X(M)$ such that $\eta|_{\pi_1(T_i)}=\eta_i$, $i = 1, 2$. 
\end{lemma}

\subsection{Infinitesimal trace calculus}

Fix an admissible real limiting datum $\xi$. We begin our analysis by proving deformation-theoretic formulas for trace functions. 

\medskip

Over the irreducible locus, the representation-to-character map is locally a holomorphic $PGL_2(\mathbb C)$ fiber bundle (Lemma 4 of \cite{paoluzziporti}). In particular, character paths admit local analytic representation lifts. Over the real locus, these lifts may be chosen real analytic. For a real analytic path of representations $\rho_t$, use the left logarithmic convention
\begin{equation*}
	\rho_t(g) = \bigl(I + tu(g) + O(t^2)\bigr)\rho_0(g) \ \ \ \ \ 
	u(g) = \left.\frac{d}{dt}\right|_{t=0}\rho_t(g)\rho_0(g)^{-1}
\end{equation*}
Then $u \in Z^1(\pi_1(T); \sltwo(\mathbb{R})_{\mathrm{Ad}\rho_0})$ and $d\mathrm{tr}_g(u) = \mathrm{tr}\bigl(u(g)\rho_0(g)\bigr)$. For $g, h \in \pi_1(T)$, the \emph{cocycle rule} states that $u(gh) = u(g) + \Ad_{\rho_0}(g)u(h)$, $u(g^{-1}) = -\Ad{\rho_0}(g)^{-1}u(g)$. In particular, if $\rho_0(g) = I$, then $u(gh) = u(g) + u(h)$ for every $h$. 

\medskip

Two facts that we use this section: the \emph{cyclicity of trace}, i.e. for $A, B \in M_2$, $\tr(AB) = \tr(BA)$, and the standard $SL_2(\mathbb{C})$ identity $\tr(AB) = \tr(A)\tr(B) - \tr(AB^{-1})$. 

\medskip

With our conventions, equality of the two Conway-boundary characters is expressed by
\begin{equation*}
	\mathrm{tr}(\mu_{1, 1}) = \mathrm{tr}(\mu_{2, 4}) \ \ \ \ \ \mathrm{tr}(\mu_{1, 2}) = \mathrm{tr}(\mu_{2, 1}) \ \ \ \ \ \mathrm{tr}(\mu_{1, 3}) = \mathrm{tr}(\mu_{2, 2}) \ \ \ \ \ \mathrm{tr}(\mu_{1, 4}) = \mathrm{tr}(\mu_{2, 3})
\end{equation*}
together with
\begin{equation*}
	\mathrm{tr}(\mu_{1, 1}\mu_{1, 2}) = \mathrm{tr}(\mu_{2, 4}\mu_{2, 1}) \ \ \ \ \ \mathrm{tr}(\mu_{1, 2}\mu_{1, 3}) = \mathrm{tr}(\mu_{2, 1}\mu_{2, 2}) \ \ \ \ \ \mathrm{tr}(\mu_{1, 1}\mu_{1, 3}) = \mathrm{tr}(\mu_{2, 4}\mu_{2, 2})
\end{equation*}

\begin{lemma}\label{lma:identitytoidentity}
Let $f_0 = \mathrm{tr}(\mu_{1, 2}) - \mathrm{tr}(\mu_{2, 1})$. At an admissible real limiting datum $\xi$, one has $d(f_0)_\xi = 0$.
\end{lemma}

\begin{proof}
If $(u_1, u_2) \in T_\xi(X(T_1) \times X(T_2))$ is a tangent cocycle pair, then
\begin{equation*}
	d(f_0)_\xi(u_1, u_2) = \mathrm{tr}(u_1(\mu_{1, 2})\rho_1(\mu_{1, 2})) - \mathrm{tr}(u_2(\mu_{2, 1})\rho_2(\mu_{2, 1})) = \mathrm{tr}(u_1(\mu_{1, 2})) - \mathrm{tr}(u_2(\mu_{2, 1})) = 0
\end{equation*}
because both cocycle values lie in $\sltwo(\mathbb{R})$ and both base matrices are $I$.
\end{proof}

We thus have a first-order degeneracy at an admissible real limiting datum $\xi \in V_C(T_1, T_2)$. The central goal of this section is to identify sufficient conditions on $T_1, T_2$ to conclude that there is a real quadratic node at $\xi$ in an affine variety containing $V_C(T_1, T_2)$; we then analyze which of the two resulting real arcs genuinely lies in $V_C(T_1, T_2)$ and glue to genuine $SL_2(\mathbb{R})$-characters of $\pi_1(M)$.

\medskip

We first prove some preliminary lemmas. 

\begin{lemma}\label{lma:linalg}
If $Q \in SL_2(\mathbb{R})$ is noncentral parabolic with fixed line $\mathcal Q$, then $\{X \in \sltwo(\mathbb{R}) \mid \mathrm{tr}(XQ) = 0\}$ is the Borel subalgebra preserving $\mathcal Q$.
\end{lemma}

\begin{proof}
After conjugation, write $Q = \epsilon\left(\begin{smallmatrix} 1 & s \\ 0 & 1 \end{smallmatrix}\right)$ with $\epsilon \in \{\pm 1\}$, $s \neq 0$, and fixed line $\mathcal Q = \langle e_1 \rangle$. Write $X = \left(\begin{smallmatrix}a & b \\ c & -a\end{smallmatrix}\right) \in \sltwo(\mathbb{R})$. Then we have $XQ = \epsilon\left(\begin{smallmatrix}a&as+b\\c&cs-a\end{smallmatrix}\right)$ and thus $\mathrm{tr}(XQ) = \epsilon cs$. Since $s \neq 0$, we have $\mathrm{tr}(XQ) = 0$ implies $c = 0$, i.e. $X\langle e_1 \rangle \subset \langle e_1 \rangle$. Thus, $\{X \in \sltwo(\mathbb{R}) \mid \mathrm{tr}(XQ) = 0\} = \left\{\begin{smallmatrix}a & b \\ 0 & -a\end{smallmatrix}\right\}$ which is the Borel subalgebra preserving $\mathcal Q = \langle e_1 \rangle$. The conclusion holds up to conjugation and hence for the original parabolic $Q$. 
\end{proof}

\begin{lemma}\label{lma:mixedtrace}
Let $T$ be a tangle, and let $m_i \in \pi_1(T)$ denote the four Conway meridians, and suppose $\rho:\pi_1(T)\to SL_2(\mathbb R)$ satisfies $\rho(m_1) = \rho(m_2) = I$ and $\rho(m_3) = \rho(m_4) = P$. For $j \in \{1, 2\}, k \in \{3, 4\}$ and any tangent cocycle $u$, $d\bigl(\mathrm{tr}(m_jm_k) - \mathrm{tr}(m_k)\bigr)(u) = \mathrm{tr}(u(m_j)P)$. 
\end{lemma}

\begin{proof}
The cocycle rule gives $u(m_jm_k)=u(m_j)+u(m_k)$, and the base image of the product is $P$. Subtracting the two trace differentials gives the formula.
\end{proof}

Keep the setting of Lemma \ref{lma:mixedtrace}, where $m_1, m_2$ are the meridians at the two ends of the identity strand and $m_2 = \tau m_1\tau^{-1}$. We compare the four differentials $d\tr(m_jm_k)$, $j \in \{1, 2\}$, $k \in \{3, 4\}$, obtained by pairing $P = \rho(m_3)$ against a meridian from each end of an identity strand. Since $\rho(m_1) = I$, the cocycle rule gives $u(m_2) = \Ad_{\rho(\tau)}u(m_1)$. Applying Lemma \ref{lma:mixedtrace} with $k = 3$ to both $j = 1$ and $j = 2$ and subtracting, we get 
\begin{equation*}
	d\tr(m_2m_3)(u) - d\tr(m_1m_3)(u) = \tr((u(m_2)-u(m_1))P) = \tr((\rho(\tau)u(m_1)\rho(\tau)^{-1} - u(m_1))P)
\end{equation*}
Cyclicity of trace turns the rightmost term into $\tr(u(m_1)(P^\tau - P))$, where $P^\tau = \rho(\tau)^{-1}P\rho(\tau)$. So, passing from $m_1$ to $m_2$ leaves $u(m_1)$ fixed and replaces $P$ by $P^\tau$; in summary, we have $d\tr(m_2m_3)(u) - d\tr(m_1m_3)(u) = \tr(u(m_1)(P^\tau - P))$. Thus, if the functional on the right vanishes on the tangent space, then $d\tr(m_2m_3) = d\tr(m_1m_3)$. The next lemma shows that if this functional vanishes on the tangent space, all four differentials $d\tr(m_1m_3)$, $d\tr(m_2m_3)$, $d\tr(m_1m_4)$, $d\tr(m_2m_4)$ coincide. This is crucial for computing the rank of the Jacobian at an admissible real limiting datum.

\begin{lemma}\label{lma:beta}
Let $T$ be a tangle with Conway holes $m_1$ connected to $m_2$ and $m_3$ connected to $m_4$. Let $\rho: \pi_1(T) \to SL_2(\mathbb{R})$ satisfy $\rho(m_1) = \rho(m_2) = I$, and $\rho(m_3) = \rho(m_4) = P$, with $P$ noncentral parabolic. Suppose $m_2 = \tau m_1\tau^{-1}$, and let $P^\tau = \rho(\tau)^{-1}P\rho(\tau)$, and suppose that for every tangent cocycle $u$ on the chosen character germ, one has $\mathrm{tr}(u(m_1)(P^\tau - P)) = 0$. Then we have 
\begin{equation*}
	d\mathrm{tr}(m_1m_3) = d\mathrm{tr}(m_2m_3) = d\mathrm{tr}(m_1m_4) = d\mathrm{tr}(m_2m_4)
\end{equation*}
\end{lemma}

\begin{proof}
The above computation shows that under the hypothesis, $d\mathrm{tr}(m_2m_3)(u) = d\mathrm{tr}(m_1m_3)(u)$. The same computation with $m_4$ in place of $m_3$ gives $d\mathrm{tr}(m_2m_4) = d\mathrm{tr}(m_1m_4)$. Finally, the trace functions of $m_3$ and $m_4$ agree because these elements are conjugate; differentiating that equality and using Lemma \ref{lma:mixedtrace} shows that the products with $m_3$ and $m_4$ have the same differential.
\end{proof}

\subsection{Proof of Theorem \ref{mainthm:a}}

We now construct the nodal pseudo-gluing germ, identify its genuine branch, and complete the proof of Theorem \ref{mainthm:a}. This establishes sufficient criteria for an ideal point to detect an essential Conway sphere with an admissible real limiting datum, and for two arcs of real traces to approach this ideal point. 

\medskip

Near an admissible datum, let $V_C'(T_1, T_2)$ be the five-equation locus
\begin{equation}\label{eq:fiveeq1}
	\mathrm{tr}(\mu_{1, 1}) = \mathrm{tr}(\mu_{2, 4}) \ \ \ \ \ \mathrm{tr}(\mu_{1, 2}) = \mathrm{tr}(\mu_{2, 1}) \ \ \ \ \ \mathrm{tr}(\mu_{1, 3}) = \mathrm{tr}(\mu_{2, 2})
\end{equation}
\begin{equation}\label{eq:fiveeq2}
	\mathrm{tr}(\mu_{1, 1}\mu_{1, 2}) = \mathrm{tr}(\mu_{2, 4}\mu_{2, 1}) \ \ \ \ \ \mathrm{tr}(\mu_{1, 2}\mu_{1, 3}) = \mathrm{tr}(\mu_{2, 1}\mu_{2, 2}) 
\end{equation}
We refer to this as the \emph{pseudo-gluing subvariety}. In the proof of Theorem \ref{mainthm:a}, we first show $V_C'(T_1, T_2)$ has an ordinary quadratic node at $\xi$. The full gluing variety $V_C(T_1, T_2)$ is obtained by imposing the additional equation $\mathrm{tr}(\mu_{1, 1}\mu_{1, 3}) = \mathrm{tr}(\mu_{2, 4}\mu_{2, 2})$. We will show that exactly one of the two local node branches satisfies this additional equation.

\medskip

For $T_1$, set
\begin{equation*}
	a_1 = \mathrm{tr}(\mu_{1, 1}) \ \ \ \ \ b_1 = \mathrm{tr}(\mu_{1, 2}) \ \ \ \ \ c_1 = \mathrm{tr}(\mu_{1, 3}) \ \ \ \ \ d_1 = \mathrm{tr}(\mu_{1, 4})
\end{equation*}
\begin{equation*}
	x_1 = \mathrm{tr}(\mu_{1, 1}\mu_{1, 2}) \ \ \ \ \ y_1 = \mathrm{tr}(\mu_{1, 2}\mu_{1, 3}) \ \ \ \ \ z_1 = \mathrm{tr}(\mu_{1, 1}\mu_{1, 3})
\end{equation*}
For $T_2$, use the boundary order adapted to the identification:
\begin{equation*}
	a_2 = \mathrm{tr}(\mu_{2, 4}) \ \ \ \ \ b_2 = \mathrm{tr}(\mu_{2, 1}) \ \ \ \ \ c_2 = \mathrm{tr}(\mu_{2, 2}) \ \ \ \ \ d_2 = \mathrm{tr}(\mu_{2, 3})
\end{equation*}
\begin{equation*}
	x_2 = \mathrm{tr}(\mu_{2, 4}\mu_{2, 1}) \ \ \ \ \ y_2 = \mathrm{tr}(\mu_{2, 1}\mu_{2, 2}) \ \ \ \ \ z_2 = \mathrm{tr}(\mu_{2, 4}\mu_{2, 2})
\end{equation*}
Thus, $V_C'(T_1, T_2)$ is defined by $a_1 = a_2$, $b_1 = b_2$, $c_1 = c_2$, $x_1 = x_2$, $y_1 = y_2$. Since the two endpoints of each strand are conjugate in the tangle group and trace is invariant under inversion, we have $a_1 = b_1$, $c_1 = d_1$, $a_2 = d_2$, $b_2 = c_2$, meaning that $d_1 = d_2$ is already implied by the original five equations. Thus, on $V_C'(T_1, T_2)$, the two boundary restrictions have the same six coordinates $a_i, b_i, c_i, d_i, x_i, y_i$. When considering such variables on $V_C'(T_1, T_2)$, denote them $a, b, c, d, x, y$. The only possible unmatched coordinate is $z_1, z_2$. 

\begin{lemma}\label{lma:branches}
On $V_C'(T_1, T_2)$, we have $(z_1 - z_2)(z_1 + z_2 + xy - (ac + bd)) = 0$. Put $\mathcal D = \{z_1 = z_2\}$, and $\mathcal H = \{z_1 + z_2 + xy - (ac + bd) = 0\}$. Then $V_C'(T_1, T_2) \subset \mathcal D \cup \mathcal H$, and the full gluing variety is $V_C(T_1, T_2) = V_C'(T_1, T_2) \cap \mathcal D$. 
\end{lemma}

\begin{proof}
For an ordered boundary tuple $(A, B, C, D)$ satisfying $ABCD = 1$, put $a = \tr A$, $b = \tr B$, $c = \tr C$, $d = \tr D$, $x = \tr(AB)$, $y = \tr(BC)$, $z = \tr(AC)$. The four-holed-sphere relation, stated as 5.2.2 of \cite{goldman}, viewed as a quadratic in $z$, is
\begin{equation*}
	0 = z^2 + (xy-ac-bd)z+x^2+y^2 - (ab+cd)x-(ad+bc)y + a^2 + b^2 + c^2 + d^2 + abcd - 4
\end{equation*}
On $V_C'(T_1, T_2)$, the two boundary tuples have the same coefficients in this quadratic. Thus $z_1$ and $z_2$ are roots of one quadratic. Subtracting their two equations gives $(z_1 - z_2)(z_1 + z_2 + xy - (ac+bd)) = 0$. The final statement follows because equality of the remaining coordinate $z$ is precisely the omitted full-gluing equation. 
\end{proof}

\begin{lemma}\label{lma:boundary-discriminant}
Let $\chi\in X(C')$ have Fricke coordinates $(a, b, c, d, x, y, z)$ in the order above, and define $\omega(\chi) = 2z+xy-(ac+bd)$. If $\chi$ is reducible, then $\omega(\chi)=0$. Consequently, $\omega(\chi)\neq0$ implies that $\chi$ is irreducible.
\end{lemma}

\begin{proof}
Represent a reducible character by a diagonal representation and write the first three eigenvalues as $\alpha, \beta, \gamma$. Then
\begin{equation*}
	a = \alpha+\alpha^{-1} \ \ \ \ \ b=\beta+\beta^{-1} \ \ \ \ \ c=\gamma+\gamma^{-1} \ \ \ \ \ d=\alpha\beta\gamma+(\alpha\beta\gamma)^{-1}
\end{equation*}
\begin{equation*}
	x=\alpha\beta+(\alpha\beta)^{-1} \ \ \ \ \ y=\beta\gamma+(\beta\gamma)^{-1} \ \ \ \ \ z=\alpha\gamma+(\alpha\gamma)^{-1}
\end{equation*}
A direct expansion gives $2z + xy = ac + bd$. Hence $\omega(\chi) = 0$. 
\end{proof}

We now describe, in outline, how these two lemmas combine at $\xi$. Every boundary trace at $\xi$ equals $2$, so the limiting datum is reducible on $\pi_1(C')$. In particular, Lemma \ref{lma:boundary-discriminant} forces $\omega_1(\xi) = \omega_2(\xi) = 0$. Step 2 in the below proof of Theorem \ref{thm:main} shows that $V_C'(T_1, T_2)$ is the zero locus, in a smooth surface $S \ni \xi$, of a single function whose Hessian at $\xi$ is indefinite. This becomes an ordinary node with two real analytic branches. Lemma \ref{lma:branches} shows that every nearby point of $V_C'(T_1, T_2)$ satisfies $\omega_1 = \omega_2$ or $\omega_1 = -\omega_2$. On $S$, the former is equivalent to the omitted gluing equation $z_1 = z_2$, so each branch must lie entirely in either $\mathcal D$ or $\mathcal H$. Because both $\omega_i$ vanish at $\xi$, distinguishing the branches requires tracking them to their first nonvanishing order, which is the content of Step 3. To make this precise, we introduce coordinates adapted to $\xi$. Since the two ends of a tangle strand are conjugate, $a_1 = b_1$, $c_1 = d_1$ on $X(T_1)$ and $a_2 = d_2$, $b_2 = c_2$ on $X(T_2)$. Combined with the matching equations $a_1 = a_2$, $c_1 = c_2$, these force 
\begin{equation*}
	p := a_1 = b_1 = a_2 = d_2 \ \ \ \ \ r := c_1 = d_1 = b_2 = c_2
\end{equation*}
on $S$. This pairing is crossed in the following sense: $T_1$'s identity strand matches $T_2$'s parabolic strand trace, and vice versa, which resulted in the nonconjugacy in Lemma \ref{lma:unrealized}. Writing $x = x_1 = x_2$, $y = y_1 = y_2$ for the remaining matched coordinates, the two boundary tuples become $(p, p, r, r)$ and $(p, r, r, p)$, with $\xi$ itself at $p = r = x = y = 2$. 

\medskip

Given a noncentral parabolic $P$, we say that $\mathrm{Fix}(P)$ is the fixed line of $P$, viewed as an action on $\mathbb{R}^2$. In addition, given $\gamma \in G$ and a representation $\rho: G \to SL_2(\mathbb C)$, write $I_\gamma: X(G) \to \mathbb C$ to be the trace function on $X(G)$.

\begin{theorem}\label{thm:main}
Let $\xi = (\chi_1, \chi_2) \in V_C'(T_1, T_2)$ be an admissible real limiting datum, represented by $\rho_i: \pi_1(T_i) \to SL_2(\mathbb{R})$, satisfying \eqref{eq:admissible1}-\eqref{eq:admissible2}. Let $\Gamma_i = \pi_1(T_i) / \langle \langle m_{i,1} \rangle \rangle$ be the quotient with representation $\bar\rho_i: \Gamma_i \to SL_2(\mathbb R)$, $q_i: \pi_1(T_i) \to \Gamma_i$ the quotient map, and $\bar m_i = q_i(m_{i, 3})$. Let $\mathcal P_i = \mathrm{Fix}(P_i)$. Choose $\tau_i \in \pi_1(T_i)$ such that $m_{i, 2} = \tau_im_{i, 1}\tau_i^{-1}$, and let $P_i^{\tau_i} = \rho_i(\tau_i)^{-1}P_i\rho_i(\tau_i)$. Let $\mathcal Q_i = \mathrm{Fix}(P_i^{\tau_i})$. Assume, for each $i$:
\begin{enumerate}
	\item $X(T_i)_{\chi_i}$ is smooth of complex dimension 3.
	\item The character germ $X(\Gamma_i)_{\chi_{\bar\rho_i}}$ is a smooth complex curve, and the trace function $I_{\bar m_i}$ is a local analytic coordinate on it.
	\item every cocycle $u_i$ representing a vector in $T_{\chi_i}X(T_i)$ satisfies $\mathrm{tr}(u_i(m_{i, 1})(P_i^{\tau_i} - P_i)) = 0$.
	\item $\mathcal P_i \neq \mathcal Q_i$.
\end{enumerate}
Then the germ of $V_C'(T_1,T_2)$ at $\xi$ is an ordinary real quadratic node, and exactly one of its two complex analytic branches, namely the branch contained in $\mathcal D = \{z_1=z_2\}$, lies in $V_C(T_1,T_2)$. This branch is the germ of an irreducible algebraic curve $Z\subset V_C(T_1,T_2)$ whose generic Conway-boundary character is irreducible and whose generic points lie in $r(X(M))$. Consequently, $\xi$ is branch-admissible, and the Conway sphere is detected by an ideal point with limiting datum $\xi$.
\end{theorem}

\begin{proof}
\emph{Step 1: rank of the Jacobian.} Let $u_i$ be a cocycle representing a vector in $T_{\chi_i}X(T_i)$, and set $X_i = u_i(m_{i, 1}) \in \sltwo(\mathbb C)$ and $X_i^\tau = \Ad_{\rho_i(\tau_i)}X_i$. Define the strand-variable functionals $\alpha_i = dI_{m_{i,3}}$ and $\beta_i = dI_{m_{i,1}m_{i,3}}$. By Lemma \ref{lma:mixedtrace}, we have 
\begin{equation*}
	\beta_i(u_i) = \alpha_i(u_i) + \tr(X_iP_i) \ \ \ \ \ \tr(X_iP_i) = \beta_i(u_i) - \alpha_i(u_i)
\end{equation*}
Cyclicity of trace and hypothesis (3) give $\tr(X_i^\tau P_i) = \tr(X_iP_i^{\tau_i}) = \tr(X_iP_i)$. Using \eqref{eq:boundaryconventions1} and \eqref{eq:boundaryconventions2}, the trace identity $\tr(AB^{-1}) = \tr(A)\tr(B)-\tr(AB)$, and Lemma \ref{lma:beta}, we obtain the differential table
\begin{equation}\label{eq:boundarydifferentials} \begin{array}{c|ccccccc} &da_i&db_i&dc_i&dd_i&dx_i&dy_i&dz_i\\ \hline T_1 &0&0&\alpha_1&\alpha_1 &0&2\alpha_1-\beta_1&\beta_1\\ T_2 &\alpha_2&0&0&\alpha_2 &2\alpha_2-\beta_2&0&\beta_2. \end{array} \end{equation}
The two identities involving a mixed product with $\mu_{i, 2}$ use hypothesis (3), via Lemma \ref{lma:beta}. (Thus, \eqref{eq:boundarydifferentials} is a dictionary on $T_{\chi_i}X(T_i)$ under hypothesis (3), not a formal cocycle identity independent of the theorem's assumptions.)

\medskip

We first work on the complex tangent spaces. Let $L_i^{\mathbb C} = \ker\alpha_i \cap \ker\beta_i \subset T_{\chi_i}X(T_i)$. For $u_i \in L_i^{\mathbb C}$, choose a cocycle representative and set $X_i = u_i(m_{i,1}) \in \sltwo(\mathbb C)$. The preceding identities give $\tr(X_iP_i) = \tr(X_iP_i^{\tau_i}) = 0$. The complexification of Lemma \ref{lma:linalg} therefore shows that $X_i$ preserves both $\mathcal P_i \otimes_{\mathbb R}\mathbb C$ and $\mathcal Q_i \otimes_{\mathbb R}\mathbb C$. By hypothesis (4), these two lines are distinct, so in the basis ($\mathcal P_i$, $\mathcal Q_i$), one has $X_i = \left(\begin{smallmatrix}\lambda & 0 \\ 0 & -\lambda\end{smallmatrix}\right)$ for $\lambda \in \mathbb C$. 

\medskip
	
We now show that evaluation at $m_{i,1}$ is injective on $L_i^{\mathbb C}$. Since $\rho_i(m_{i,1}) = I$, every coboundary vanishes at $m_{i, 1}$; hence $L_i^{\mathbb C} \to \sltwo(\mathbb C)$, $u_i \mapsto u_i(m_{i, 1})$ is well-defined on character tangent classes. Suppose for contradiction that $0 \neq u_i \in L_i^{\mathbb C}$ and $X_i = u_i(m_{i, 1}) = 0$. The cocycle rule shows that $u_i$ vanishes on the normal closure of $m_{i, 1}$, so it descends to a cocycle $\bar u_i \in Z^1(\Gamma_i; \sltwo(\mathbb C)_{\Ad\bar\rho_i})$. Since $\bar\rho_i$ is irreducible, the local holomorphic $PGL_2(\mathbb C)$-bundle over the irreducible locus identifies $T_{\chi_{\bar\rho_i}}X(\Gamma_i) \cong Z^1(\Gamma_i; \sltwo(\mathbb C)_{\Ad\bar\rho_i})/B_1(\Gamma_i; \sltwo(\mathbb C)_{\Ad\bar\rho_i})$. Thus $\bar u_i$ determines a tangent vector $[\bar u_i]$, which is zero precisely when $\bar u_i$ is a coboundary. Since $dI_{\bar m_i}([\bar u_i]) = \alpha_i(u_i) = 0$ and $I_{\bar m_i}$ is a local analytic coordinate by hypothesis (2), its differential is injective. Hence $[\bar u_i] = 0$, so $\bar u_i$, and therefore its pullback $u_i$is a coboundary. This contradicts $u_i \neq 0$ in $T_{\chi_i}X(T_i)$. Hence, every nonzero $u_i \in L_i^{\mathbb C}$ has $u_i(m_{i, 1}) \neq 0$, and evaluation at $m_{i, 1}$ is injective.

\medskip

Its image lies in the one-dimensional subspace of $\sltwo(\mathbb C)$ preserving the two distinct complexified lines $\mathcal P_i$ and $\mathcal Q_i$. Hence $\dim_{\mathbb C}L_i^{\mathbb C} \leq 1$. On the other hand, $T_{\chi_i}X(T_i)$ has dimension three, so the intersection of the kernels of the two linear functionals $\alpha_i$ and $\beta_i$ has dimension at least one. Therefore, $\dim_{\mathbb C}L_i^{\mathbb C} = 1$, and $\alpha_i$, $\beta_i$ are linearly independent.

\medskip

All of the defining data are real, so $L_i^{\mathbb C}$ is the complexification of a one-dimensional real line $L_i^{\mathbb R} = L_i$. Choose a nonzero vector $v_i \in L_i$. Then the diagonal coefficient above is real and nonzero, and therefore $Q_i := \frac12\mathrm{tr}(v_i(m_{i, 1})^2) > 0$.

\medskip

Differentiating the five pseudo-gluing equations \eqref{eq:fiveeq1}-\eqref{eq:fiveeq2} gives the Jacobian
\begin{equation}\label{eq:gluingjacobian} \begin{pmatrix} 0&-\alpha_2\\ 0&0\\ \alpha_1&0\\ 0&-(2\alpha_2-\beta_2)\\ 2\alpha_1-\beta_1&0 \end{pmatrix}. \end{equation}
The zero row is the identity-to-identity equation $f_0 = b_1 - b_2$. Since $\alpha_i$ and  $\beta_i$ are independent, the other four rows have rank four. The four nonzero rows have complex rank four and cut out a smooth complex analytic surface $S_{\mathbb C}$, defined over $\mathbb R$, in the six-dimensional complex product germ. Its real locus $S$ is a smooth real analytic surface, and $T_\xi S_{\mathbb C} = L_1^{\mathbb C} \oplus L_2^{\mathbb C}$, and $T_\xi S = L_1 \oplus L_2$. For comparison, the implied equation $d_1=d_2$ has differential $(\alpha_1, -\alpha_2)$, while the omitted full-gluing equation $z_1=z_2$ has differential $(\beta_1, -\beta_2)$.
	
\medskip

\emph{Step 2: the quadratic term and the node.} We now compute the quadratic term of $f_0=b_1-b_2 = \tr(\mu_{1,2})-\tr(\mu_{2,1})$ on $T_\xi S$. Since $\rho_i(\mu_{i,1}) = I$, we may write $\rho_{i,t}(m_{i,1}) = I + tX_i + t^2Y_i + O(t^3)$. Since $\rho_{i,t}(m_{i,1}) \in SL_2(\mathbb{R})$, using the general identity $\det(I + tM) = 1 + t\tr(M) + t^2\det(M)$, we have
\begin{equation*}
	1 = \det\rho_{i,t}(m_{i,1}) = 1 + t(\mathrm{tr}(X_i)) + t^2(\mathrm{tr}(Y_i) + \mathrm{det}(X_i)) + O(t^3)
\end{equation*} 
Taking the $t^2$ term to be zero and using $\det(X_i) = -\frac12\mathrm{tr}(X_i^2)$, we get $\mathrm{tr}(Y_i)=\frac12\mathrm{tr}(X_i^2)$. Therefore $\mathrm{tr}\rho_{i,t}(m_{i,1}) = 2+\frac{t^2}{2}\operatorname{tr}(X_i^2)+O(t^3)$. The quadratic term of $f_0$ is hence $Q_\xi = \frac12\mathrm{tr}(X_1^2) - \frac12\mathrm{tr}(X_2^2)$. (Here we also use the fact that $m_{1,2}$ is conjugate to $m_{1,1}$, so they have the same trace differentials.)

\medskip
	
Choose nonzero real generators $v_i \in L_i$ and put $Q_i = \frac12\mathrm{tr}(v_i(m_{i,1})^2) > 0$. In the real coordinates $(s_1, s_2)$ tangent to $(v_1, 0)$ and $(0, v_2)$, the quadratic part of $f_0|_{S(\mathbb R)}$ is $Q_1s_1^2 - Q_2s_2^2$. Its complexification is nondegenerate. By the holomorphic Morse lemma (see 1.6 of \cite{morse}), the complex germ is analytically equivalent to a nondegenerate quadratic form, hence in two variables to $f_0 = u^2 - v^2$. Thus the complex germ is an ordinary quadratic node with two complex analytic branches. By the real-analytic Morse lemma (1.7 of \cite{morse}), since the real Hessian has signature $(1, 1)$, the same real normal form $u^2 - v^2$ may be chosen on the real locus. Consequently, the two complex branches are defined over $\mathbb R$, and their real loci are smooth and transverse.  

\medskip

\emph{Step 3: branch selection.} It remains to decide which of the two node branches satisfies the omitted full-gluing equation $z_1=z_2$. Recall the coordinates $p, r, x, y$ described in the discussion preceding the theorem statement. We have $f_0 = p - r$. The ordered boundary traces on the two sides are respectively $(a_1, b_1, c_1, d_1) = (p, p, r, r)$, and $(a_2, b_2, c_2, d_2) = (p, r, r, p)$. For $i = 1, 2$, the Fricke relation is a quadratic in the unmatched coordinate $z_i$, i.e. $z_i^2+(xy-2pr)z_i+G_i=0$, where
\begin{equation*}
	G_1 = x^2 + y^2 - (p^2 + r^2)x - 2pry + 2p^2 + 2r^2 + p^2r^2 - 4 \ \ \ \ \ G_2 = x^2 + y^2 - 2prx - (p^2 + r^2)y + 2p^2 + 2r^2 + p^2r^2 - 4
\end{equation*} 
Define $\omega_i = 2z_i - (2pr - xy) = 2z_i + xy - 2pr$ for $i = 1, 2$, and $\Delta_i = (2pr - xy)^2 - 4G_i$. Then $\omega_i^2 = \Delta_i$, and on $S$, 
\begin{equation*}
	z_1 - z_2 = \frac{\omega_1 - \omega_2}{2} \ \ \ \ \ z_1 + z_2 + xy - 2pr = \frac{\omega_1 + \omega_2}{2}
\end{equation*}
The goal of this computation is to show that precisely one node branch satisfies $\omega_1 = \omega_2$, while the other satisfies $\omega_1 = -\omega_2$. Because $ac+bd = 2pr$ on $S$, these are precisely the two factors from Lemma \ref{lma:branches}. Let
\begin{equation*}
	P_0 = p - 2 \ \ \ \ \ R_0 = r - 2 \ \ \ \ \ X_0 = x - 2 \ \ \ \ \ Y_0 = y - 2
\end{equation*}
The following computation appears in the appendix. By Lemma \ref{lma:discriminantterms}, the constant, linear, and quadratic homogeneous parts of $\Delta_1, \Delta_2$ vanish, and the cubic homogeneous parts are 
\begin{equation*}
	\Delta_1^{(3)} = 4X_0\bigl(P_0^2 - 2P_0R_0 - 2P_0Y_0 + R_0^2 - 2R_0Y_0 + X_0Y_0 + Y_0^2\bigr)
\end{equation*}
\begin{equation*}
	\Delta_2^{(3)} = 4Y_0\bigl(P_0^2 - 2P_0R_0 - 2P_0X_0 + R_0^2 - 2R_0X_0 + X_0^2 + X_0Y_0\bigr)
\end{equation*}
Choose local coordinates $s_1, s_2$ on $S$, tangent respectively to nonzero vectors $(v_1, 0)$ and $(0, v_2)$ in $L_1 \oplus L_2$, and put $U = Q_1s_1^2$, and $V = Q_2s_2^2$. The differential table \eqref{eq:gluingjacobian} shows that $dP_0|_{T_\xi S} = dR_0|_{T_\xi S} = dX_0|_{T_\xi S} = dY_0|_{T_\xi S} = 0$. Since $P_0(\xi) = R_0(\xi) = X_0(\xi) = Y_0(\xi) = 0$, their restrictions to $S$ vanish to second order. Since $\rho_i(m_{i, 1}) = I$, Step 2's determinant computation applies to $g = m_{i, 1}$ for any tangent direction $u_i$, giving $\tr\rho_{i,t}(m_{i,1}) = 2 + \frac{t^2}{2}\tr(u_i(m_{i,1})^2) + O(t^3)$. As $a_1$ and $b_2$ each depend on only one of $\chi_1$, $\chi_2$, restricting to $u = s_1v_1 + s_2v_2 \in L_1 \oplus L_2$ (where $v_2$ has no $T_1$-component and $v_1$ has no $T_2$-component), and polarizing gives $j_2(P_0) = Q_1s_1^2 = U$ and $j_2(R_0) = Q_2s_2^2 = V$. To compute the remaining two jets, write along the $i$th tangent direction
\begin{equation*}
	\rho_{i, s}(m_{i, 1}) = I + sX_i + s^2Y_i + O(s^3) \ \ \ \ \ \rho_{i, s}(m_{i, 2}) = I + sX_i^\tau + s^2Y_i^\tau + O(s^3)
\end{equation*}
The Step-2 calculation gives $\tr Y_i = \tr Y_i^\tau = Q_i$, and $\tr(X_i^2) = \tr((X_i^\tau)^2) = 2Q_i$. In the ordered basis $(\mathcal P_i, \mathcal Q_i)$, write $X_i = \left(\begin{smallmatrix}\lambda_i&0\\0&-\lambda_i\end{smallmatrix}\right)$. Put $R_i = \rho_i(\tau_i)$. Since $\mathcal Q_i = R_i^{-1}\mathcal P_i$, the line $R_i\mathcal Q_i = \mathcal P_i$ is an eigenline of $X_i^\tau = R_iX_iR_i^{-1}$ with eigenvalue $-\lambda_i$. Hence, in the same basis, $X_i^\tau$ has diagonal entries $-\lambda_i$, $\lambda_i$. Therefore $\tr(X_iX_i^\tau) = -2\lambda_i^2 = -2Q_i$. By Lemma \ref{lma:crossjet}, $x_1 = 2 + 4Q_1s^2 + O(s^3)$, and $y_2 = 2 + 4Q_2s^2 + O(s^3)$. Consequently, $j_2(X_0) = 4U$ and $j_2(Y_0) = 4V$. Since $P_0$, $R_0$, $X_0$, $Y_0$ vanish to the second order, every homogeneous term of degree at least four in $\Delta_i$ has order at least eight on $S$. Moreover, replacing one argument of a cubic homogeneous polynomial by its order-three remainder changes the result only in order at least seven. Hence the order-six jets of $\Delta_i|_S$ are obtained by substituting the quadratic jets of $P_0$, $R_0$, $X_0$, $Y_0$ into $\Delta_i^{(3)}$. This substitution gives
\begin{equation*}
	j_6(\omega_1^2|_S) = 16U(U + 3V)^2 \ \ \ \ \ j_6(\omega_2^2|_S) = 16V(3U+V)^2
\end{equation*}
Therefore, there are fixed signs $\sigma_1, \sigma_2 \in \{\pm 1\}$ such that
\begin{equation*}
	W_1 = j_3(\omega_1|_S) = 4\sigma_1 \sqrt{Q_1}s_1\bigl(Q_1s_1^2 + 3Q_2s_2^2\bigr)
\end{equation*}
\begin{equation*}
	W_2 = j_3(\omega_2|_S) = 4\sigma_2 \sqrt{Q_2}s_2\bigl(3Q_1s_1^2 + Q_2s_2^2\bigr)
\end{equation*}
The node quadratic is defined by $U - V = Q_1s_1^2 - Q_2s_2^2$. The two tangent lines of the node are $\sqrt{Q_1}s_1 = \sqrt{Q_2}s_2$ and $\sqrt{Q_1}s_1 = -\sqrt{Q_2}s_2$. On the first line, with both sides equal to $s$, we have $W_1 = 16\sigma_1s^3$ and $W_2 = 16\sigma_2s^3$. On the second line, with $\sqrt{Q_1}s_1 = s$ and $\sqrt{Q_2}s_2 = -s$, we have $W_1 = 16\sigma_1s^3$ and $W_2 = -16\sigma_2s^3$

\medskip

Thus, on each tangent line, exactly one of the cubic initial forms $W_1 - W_2$, $W_1 + W_2$ is nonzero. If a holomorphic function vanishes identically on a smooth curve germ, then its initial homogeneous form vanishes on the tangent line of that germ. Hence a nonzero restriction of $W_1 \pm W_2$ to a node tangent line rules out the corresponding factor on the branch tangent to that line. Each smooth node branch is irreducible, and $(\omega_1-\omega_2)(\omega_1+\omega_2) = 0$ identically on it. Therefore, one node branch must lie in $\mathcal D = \{z_1 = z_2\} = \{\omega_1 = \omega_2\}$, while the other must lie in $\mathcal H = \{z_1 + z_2 + xy - (ac + bd) = 0\} = \{\omega_1 = -\omega_2\}$. Since the full gluing variety is obtained by imposing $z_1 = z_2$, its reduced germ at $\xi$ is precisely the smooth branch contained in $\mathcal D$. 

\medskip

\emph{Step 4: realization and the ideal point.} Let $\mathcal B$ denote the complex analytic branch contained in $\mathcal D$. Since $V_C(T_1,T_2)$ is an affine algebraic set and its reduced analytic germ at $\xi$ is the smooth one-dimensional germ $\mathcal B$, there is a unique irreducible algebraic component $Z\subset V_C(T_1,T_2)$ through $\xi$ whose analytic germ is $\mathcal B$. The local dimension at $\xi$ is one, and hence $Z$ is an irreducible algebraic curve. 

\medskip

We claim that the generic Conway-boundary character on $Z$ is irreducible. Choose a holomorphic local parameter $s$ on $\mathcal B$, compatible with its real structure. Recall that, on the two tangle sides, $\omega_i = 2z_i-(2pr-xy) = 2z_i + xy - 2pr$. The cubic calculation gives $\omega_i(s) = \kappa_is^3 + O(s^4)$, $\kappa_i \neq 0$. After shrinking the complex disk, $\omega_i(s) \neq 0$ for $0 < |s| < \epsilon$. Since the restrictions $\omega_i|_Z$ have nonzero cubic initial terms, neither is identically zero on $Z$. Therefore $Z^\circ = Z \setminus (V(\omega_1) \cup V(\omega_2))$ is a nonempty Zariski-open subset. Every point of $Z^\circ$ has irreducible Conway boundary character by Lemma \ref{lma:boundary-discriminant}, and lies in $r(X(M))$ by Lemma \ref{lma:irreduciblegluing}. So $\xi$ is branch-admissible. 

\medskip

Finally, Lemma \ref{lma:unrealized} shows that $\xi$ is nonrealized and that the two limiting boundary representations are not conjugate. Lemma \ref{lma:detection} now produces the irreducible algebraic curve $\mathcal{C}\subset X(M)$ and ideal point $x$ whose limiting characters are $\chi_1$ and $\chi_2$ and which detects the Conway sphere $C$.
\end{proof}

\begin{corollary}\label{cor:realbranches}
Let $\mathcal C\subset X(M)$ be the irreducible algebraic curve obtained from Theorem \ref{thm:main}, and let $x\in\widetilde{\mathcal C}\setminus\mathcal C$ be the corresponding ideal point. Then $x$ is a real point of $\widetilde{\mathcal C}$, and there is a holomorphic local coordinate at $x$, defined over $\mathbb R$, whose restriction gives a signed real analytic coordinate
\begin{equation*}
	t:(\widetilde{\mathcal C}(\mathbb R),x)\longrightarrow(\mathbb R,0)
\end{equation*}
such that, after shrinking the parameter interval, the following hold for every $t\neq0$:
\begin{enumerate}
	\item the corresponding character $\chi_t\in X(M)$ is real-valued and irreducible;
	\item if $\mu$ is a knot meridian, then $I_\mu(t) = \chi_t(\mu) = 2+Qt^2+O(t^3)$ for some $Q>0$;
	\item $\chi_t$ is the character of an irreducible representation $\rho_t:\pi_1(M)\longrightarrow SL_2(\mathbb R)$
\end{enumerate}
\end{corollary}

\begin{proof}
The genuine branch $\mathcal B$ is defined over $\mathbb R$, and hence so is the unique algebraic curve $Z$ having germ $\mathcal B$ at $\xi$. By step 4 of the proof of Theorem \ref{thm:main}, the punctured germ of $Z$ has irreducible common Conway boundary character. For every sufficiently small nonzero real point $\eta \in Z$, Lemma \ref{lma:irreduciblegluing} gives a unique global character $\chi \in X(M)$. Since the two restrictions of $\chi$ are real-valued, the conjugate character $\overline\chi$ has the same restrictions; uniqueness therefore gives $\overline\chi = \chi$. The restriction map from the curve $\mathcal C$ supplied by Lemma \ref{lma:detection} to $Z$ is consequently generically one-to-one. It is therefore birational, and the induced morphism between the smooth projective normalizations is an isomorphism. The unique generic character above $Z$ is fixed by complex conjugation, so the corresponding rational section and its closure $\mathcal C$ are defined over $\mathbb R$, as is the birational map $\mathcal C \dashrightarrow Z$. Hence the point $x$ lying over $\xi$ is real, and the real coordinate on the normalization of $Z$ gives a signed real analytic coordinate at $x$. 

\medskip
	
We next compute the meridian trace. Let $u_i\in L_i=\ker\alpha_i\cap\ker\beta_i$ be as in the proof of
Theorem \ref{thm:main}, and write the tangent vector to the genuine branch as $a u_1+b u_2$. The quadratic node equation is $Q_1a^2-Q_2b^2=0$, where $Q_i = \frac12\operatorname{tr}\bigl(u_i(m_{i,1})^2\bigr)>0$. Since $t$ is a local parameter, the tangent vector is nonzero. Since $Q_1a^2 = Q_2b^2$ with $Q_1, Q_2 > 0$, we have $a \neq 0$ and $b \neq 0$. Choose a local complex analytic representation slice on the first tangle side and write
\begin{equation*}
	\rho_{1,t}(m_{1,1}) = I+tX+t^2Y+O(t^3), \ \ \ \ \ X=a\,u_1(m_{1,1})
\end{equation*}
Since the matrices have determinant one, $\operatorname{tr}(X)=0$ and $\operatorname{tr}(Y) = -\det X = \frac12\operatorname{tr}(X^2)$ (the computation in Step 2 of Theorem \ref{thm:main}). Therefore
\begin{equation*}
	I_\mu(t) = \operatorname{tr}\rho_{1,t}(m_{1,1}) = 2+ \frac12\operatorname{tr}(X^2)t^2 + O(t^3)
\end{equation*}
The element $m_{1,1}$ represents the knot meridian up to conjugacy and inversion, neither of which changes trace. Thus $I_\mu(t) = 2+Qt^2+O(t^3)$, where $Q = \frac12\operatorname{tr}(X^2) = a^2Q_1 >0$. 

\medskip
	
Finally, an irreducible real-valued $SL_2(\mathbb C)$-character is represented either in $SU(2)$ or $SL_2(\mathbb R)$; see \cite{goldman} or \cite{dunfieldrasmussen}. But every element of $SU(2)$ has trace in $[-2,2]$, whereas $I_\mu(t)>2$ for all sufficiently small $t\neq0$. Hence the $SU(2)$ real form is impossible. Therefore each $\chi_t$ is the character of an irreducible representation $\rho_t:\pi_1(M)\longrightarrow SL_2(\mathbb R)$.
\end{proof}

%% file: trans.tex
Let $M=T_1\cup_C T_2$ be a knot exterior with a chosen essential Conway sphere, and assume throughout that the hypotheses of Theorem \ref{thm:main} hold. Fix one of the two real half-branches supplied by Corollary \ref{cor:realbranches} and denote it by $B$.

\medskip

In this section, we prove Theorem \ref{mainthm:b} in the following steps.
\begin{enumerate}
	\item Algebraize $B$ and choose integral models of its two tangle restrictions.
	\item Break the longitude into four pieces, each contained in a tangle exterior group. This is called a \emph{four-longitude decomposition}. 
	\item Define \emph{positive parabolic normal form}, which is a behavior which allows us to bypass the defect-one ambiguity of the quasimorphism and directly sum translation numbers of the four pieces.
	\item The two half-branches produce unbounded augmented holonomy arcs with vertical asymptotes on opposite sides; if the translation numbers of the four-longitude decomposition sum to 0, $B$ lifts to an augmented branch in $H_{0,0}$ which intersects every sufficiently steep positive and negative filling line.
	\item Use the results of \cite{gao} to conclude orderability for Dehn filling at sufficiently large positive and negative slopes.   
\end{enumerate}
We finally give a practical criterion for verifying positive parabolic normal form. 

\subsection{Algebraization and integral models}

The estimates below require uniform control as the half-branch approaches the ideal point. We therefore replace the analytic family of characters by an algebraic representation over a discretely valued field, where reduction via the maximal ideal records limiting matrices, positive valuation records vanishing, and negative valuation records poles. The resulting limiting information will be used to show that under certain conditions, translation numbers can be directly summed. 

\medskip

Let $\mathcal C\subset X(M)$ be the irreducible character curve obtained from Theorem \ref{thm:main}, and let $x\in\widetilde{\mathcal C}\setminus\mathcal C$ be its ideal point. By construction, the extended restriction map takes $x$ to the admissible real limiting datum $\xi=(\chi_1,\chi_2)\in V_C(T_1,T_2)$. By Corollary \ref{cor:realbranches}, there is a signed real analytic normalization
\begin{equation*}
	\nu:(-\epsilon,\epsilon) \longrightarrow\widetilde{\mathcal C}(\mathbb R) \ \ \ \ \ \nu(0)=x
\end{equation*} 
such that, for every sufficiently small $t\neq0$, the point $\nu(t)$ represents an irreducible real character $\chi_t$ which is the character of an irreducible representation $\rho_t:\pi_1(M)\longrightarrow SL_2(\mathbb R)$. Moreover, $\lim_{t\to0}\chi_t|_{\pi_1(T_i)}=\chi_i$ and $I_\mu(t)=2+Qt^2+O(t^3)$, with $Q > 0$.

\medskip

Write $B^+=\nu((0,\epsilon))$, $B^-=\nu((-\epsilon,0))$. For the one-sided arguments below, fix either half-branch and denote it by $B$.

\begin{prop}\label{prop:branchmodel}
Let $B$ be either real half-branch. After reparametrizing $B$ and passing to a finite real Puiseux extension $K/\mathbb{R}((u))_{\mathrm{conv}}$ with residue field $\mathbb R$, there exist $\rho_B: \pi_1(M) \to SL_2(K)$ and $C_i \in SL_2(K)$, $i = 1, 2$ such that:
\begin{enumerate}
	\item For sufficiently small $u > 0$, specialization of $\rho_B$ gives a representation with character on $B$.
	\item $C_i\rho_B(\pi_1(T_i))C_i^{-1} \subset SL_2(\mathcal O_K)$.
	\item Reduction mod the maximal ideal of $K$ gives irreducible representations $\rho_i^B: \pi_1(T_i) \to SL_2(\mathbb{R})$ with characters $\chi_i$. 
\end{enumerate}
\end{prop}

A full proof is provided in Appendix \ref{sec:algebra}. For expositional purposes, we only provide a rough sketch here.

\begin{proof}[Sketch of proof]
The proof realizes $B$ as a real algebraic arc via a Nash-approximation argument on the semialgebraic set of representations restricting to $B$; this uses standard results in real algebraic geometry, see e.g. \cite{bochnak-coste-roy}. Then pass to the normalization of its Zariski closure to obtain a discrete valuation, and the tautological representation of \cite{cullershalen}. Finally, adjoin square roots of two determinants to fix the tangle-restricted representations into $SL_2(K)$ rather than merely $GL_2(K)$. 
\end{proof}

We fix the reductions $\rho_i^B$ as the \emph{branch-adapted limiting representatives}. 

\begin{definition}\label{def:conjugator}
Choose $C_i \in SL_2(K)$ as in Proposition \ref{prop:branchmodel}, and set $R = C_1C_2^{-1} \in SL_2(K)$. Let $\varpi$ be a uniformizer of $K$, positive for small $u > 0$, and put
\begin{equation*}
	m_R = \min_{j,k}v(R_{jk}) \ \ \ \ \ \widehat R = \varpi^{-m_R}R \ \ \ \ \ \widehat\Delta = \det(\widehat R) = \varpi^{-2m_R}	
\end{equation*}
Then every entry of $\widehat{R}$ lies in $\mathcal{O}_K$, at least one entry is a unit, and $\widehat\Delta > 0$ at every small $u > 0$. We refer to $R$ as the \emph{relative conjugator}, and $\widehat R$ as the \emph{projective relative conjugator}. 
\end{definition}

\begin{lemma}\label{lma:unbounded}
We have $R \notin M_2(\mathcal O_K)$. Consequently, $\widehat\Delta \in \mathfrak m_K$, and the reduction $J=\operatorname{red}(\widehat R)$ is nonzero and singular.
\end{lemma}

\begin{proof}
For $\eta\in\pi_1(C_1')$, let $\alpha_{1,B}(\eta) = C_1\rho_B(\iota_{1*}\eta)C_1^{-1}$, and $\alpha_{2,B}(\phi_*(\eta)) = C_2\rho_B(\iota_{2*}\phi_*(\eta))C_2^{-1}$. By our labeling conventions, the two boundary elements represent the same element of $\pi_1(M)$. Since $R = C_1C_2^{-1}$, 
\begin{equation}\label{eq:intertwiner}
	\alpha_{1,B}(\eta)R = R\alpha_{2,B}(\phi_*(\eta)) \ \ \ \ \ \eta \in \pi_1(C_1')
\end{equation}
Suppose that $R \in M_2(\mathcal O_K)$, and put $\overline R = \mathrm{red}(R)$. Since $\det R = 1$, $\overline R \in SL_2(\mathbb{R})$. Reducing \eqref{eq:intertwiner} gives $\rho_1^B(\eta)\overline R = \overline R\rho_2^B(\phi_*(\eta))$. After conjugating $\rho_2^B$ by $\overline R$, the two limiting tangle representations agree on the Conway subgroup. The universal property of the amalgamated product therefore gives a representation $\rho_0:\pi_1(M)\longrightarrow SL_2(\mathbb R)$ whose restrictions to the two tangle groups have characters $\chi_1$ and $\chi_2$. Thus $\xi\in r(X(M))$ contrary to Lemma \ref{lma:unrealized}. Therefore $R \notin M_2(\mathcal O_K)$; equivalently $m_R < 0$, so $v(\widehat\Delta) = -2m_Rv(\varpi) > 0$, i.e. $\widehat\Delta \in \mathfrak m_K$. Finally, at least one entry of $\widehat R$ is a unit, so $J\neq0$, while $\det J = \operatorname{red}(\widehat\Delta) = 0$, so $J$ is singular.
\end{proof}

Let $\alpha_i^B = \rho_i^B|_{\pi_1(C_i')}$, and let $\mathcal P_i \subset \mathbb R^2$ be the fixed line of the nontrivial Conway parabolics in the image of $\alpha_i^B$. 

\begin{lemma}\label{lma:rankone}
Let $J = \mathrm{red}(\widehat R)$. Then $\ker J = \mathcal P_2$ and $\im J = \mathcal P_1$. In particular, $J$ has rank one.
\end{lemma}

\begin{proof}
Reducing \eqref{eq:intertwiner} gives $\alpha_1(\eta)J = J\alpha_2(\phi_*(\eta))$. Let $\eta=\mu_{1,1}$. By the admissible limiting datum conventions, we have $\alpha_1^B(\mu_{1,1}) = I$, and $\alpha_2^B(\phi_*(\mu_{1,1})) = \alpha_2^B(\mu_{2,4}) = P_2^{-1}$. Hence $J(P_2^{-1}-I) = 0$. The image of $P_2^{-1}-I$ is $\mathcal P_2$, so $\mathcal P_2 \subseteq \ker J$. Similarly, take $\eta = \mu_{1,3}$. Then $\alpha_1^B(\mu_{1,3}) = P_1^{-1}$ and $\alpha_2^B(\phi_*(\mu_{1,3})) = \alpha_2^B(\mu_{2,2}) = I$, so $(P_1^{-1}-I)J = 0 $, so $\im J \subseteq \mathcal P_1$. By Lemma \ref{lma:unbounded}, $J$ is nonzero and singular, so it has rank one. It follows that $\ker J$ is one-dimensional, so the two inclusions are equalities. 
\end{proof}

\begin{lemma}\label{lma:analyticgauge}
After shrinking $B$ and enlarging $K$ by a finite real extension, the normalizers $C_i$ in Proposition \ref{prop:branchmodel} may be chosen so that $\Sigma_i(t) = C_i\rho_B|_{\pi_1(T_i)}C_i^{-1}$ is real analytic in the signed parameter $t$ and extends to $t = 0$ with $\Sigma_i(0) = \rho_i^B$. In addition, the projective relative conjugator admits a meromorphic representative in $t$. The resulting change of normalizers does not alter $J = \mathrm{red}(\widehat R)$ up to multiplication by a nonzero real scalar. 
\end{lemma}

\begin{proof}
Since $\rho_i^B$ is irreducible, the character map over the irreducible locus is a holomorphic fiber bundle (Lemma 4 of \cite{paoluzziporti}) and hence admits a local real-analytic section. Choose a real analytic family $\sigma_i(t): \pi_1(T_i) \to SL_2(\mathbb R)$ with $\sigma_i(0) = \rho_i^B$, $\chi_{\sigma(t)} = \chi_t|_{\pi_1(T_i)}$. For $t \neq 0$, the representations $C_i\rho_B|_{\pi_1(T_i)}C_i^{-1}$ and $\sigma_i(t)$ have the same irreducible character. Hence they are conjugate by a projectively unique matrix $H_i \in GL_2(K)$. Scale $H_i$ so that its entries are integral and its reduction is nonzero. Reducing $\sigma_i(t)H_i(t) = H_i(t)\Sigma_i(t)$ shows that $\mathrm{red}(H_i)$ centralizes the irreducible representation $\rho_i^B$, and is therefore a nonzero scalar. After a further scalar normalization and a finite real extension, we may take $H_i \in SL_2(\mathcal O_K)$, $\mathrm{red}(H_i) = I$. Replacing $C_i$ by $H_iC_i$ makes the normalized restriction equal to $\sigma_i(t)$. The projective relative conjugator is determined by linear intertwining equations with analytic coefficients, so its projective class is meromorphic in $t$. Since the $H_i$ reduce to $I$, replacing $C_i$ by $H_iC_i$ changes $\operatorname{red}(\widehat R)$ only by a nonzero scalar.
\end{proof}

\subsection{Four-longitude decompositions and trace asymptotics}

Let $\lambda$ denote the homological longitude of $M$. We break $\lambda$ into four pieces lying alternately in the two tangle groups and use this decomposition to formulate \emph{positive parabolic normal form}. Once a four-longitude decomposition is put into parabolic normal form, we then analyze the asymptotics of the ascending terms of the four-longitude decomposition. These asymptotics, combined with the defect-one quasimorphism property, will enable direct summation of the translation numbers of the four pieces.

\medskip

Let $\xi = (\chi_1, \chi_2) \in V_C(T_1, T_2)$ be an admissible real limiting datum, and fix branch-adapted limiting representatives $\rho_i^B: \pi_1(T_i) \to SL_2(\mathbb R)$. A \emph{longitudinal piece} $\ell \in \pi_1(T_i)$ is a group element satisfying the traversal identity $m_{i, 2} = \ell^{-1}m_{i, 1}\ell$ or $m_{i, 4} = \ell^{-1}m_{i, 3}\ell$. Visually, this means that $\ell$ can be seen as a strand traversal from $m_{i, 1} \to m_{i, 2}$ or $m_{i, 3} \to m_{i, 4}$. Longitudinal pieces are not unique; one may twist by meridians $m_{i, k}$ or their conjugates. 

\begin{definition}\label{def:fourlongitude}
A \emph{four-longitude decomposition along $B$} is a decomposition of $\lambda$ into four longitudinal pieces $\lambda = \ell_1\ell_2\ell_3\ell_4$, where:
\begin{enumerate}
	\item $\ell_i, \ell_{i+2} \in \pi_1(T_i)$ for $i = 1, 2$.
	\item $m_{i, 2} = \ell_i^{-1}m_{i,1}\ell_i$ and $m_{i, 4} = \ell_{i+2}^{-1}m_{i, 3}\ell_{i+2}$.
\end{enumerate}
We say that a four-longitude decomposition is in \emph{parabolic form} if, in addition:
\begin{enumerate}
	\item $\rho_i^B(\ell_i)$ and $\rho_i^B(\ell_{i+2})$ are noncentral parabolic for $i = 1, 2$.
	\item $\mathrm{Fix}(\rho_i^B(\ell_i)) \neq \mathrm{Fix}(\rho_i^B(\ell_{i+2}))$.
\end{enumerate}
\end{definition}

\begin{figure}[h]
	\centering
	\includegraphics[scale=.6]{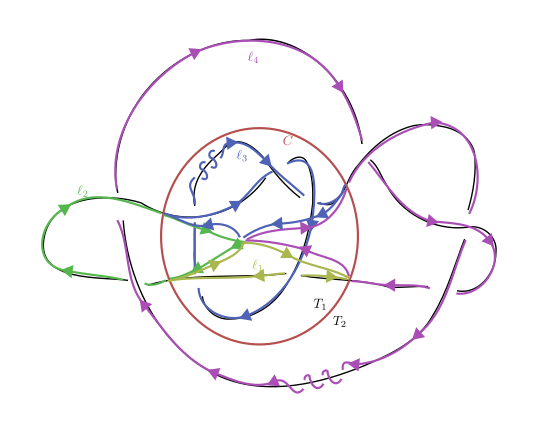}
	\caption{A four-longitude decomposition for $K8a14$.}\label{fig:fourlongitude}
\end{figure}

Let $\rho_B: \pi_1(M) \to SL_2(K)$ with normalizations $C_1, C_2$ and let $\widehat R$ be the projective relative conjugator. Let $\lambda = \ell_1\ell_2\ell_3\ell_4$ be a four-longitude decomposition in parabolic form. Put
\begin{equation*}
	A_1 = C_1\rho_B(\ell_1)C_1^{-1} \ \ \ \ \ B_2 = C_2\rho_B(\ell_2)C_2^{-1}
\end{equation*}
\begin{equation*}
	A_3 = C_1\rho_B(\ell_3)C_1^{-1} \ \ \ \ \ B_4 = C_2\rho_B(\ell_4)C_2^{-1}
\end{equation*}
We will work with the asymptotics of the traces
\begin{equation}\label{eq:firstasymptotic}
	\tr\rho_B(\ell_1\ell_2) = \tr(A_1\widehat RB_2\widehat R^{-1}) \ \ \ \ \ \tr\rho_B(\ell_1\ell_2\ell_3) = \tr(A_1\widehat RB_2\widehat R^{-1}A_3)
\end{equation}
\begin{equation}\label{eq:thirdasymptotic}
	\tr\rho_B(\lambda) = \tr(A_1\widehat RB_2\widehat R^{-1}A_3\widehat RB_4\widehat R^{-1})
\end{equation}
as one approaches the ideal point constructed in Section \ref{sec:idealpoints}. For each $i$, put $\mathcal L_i = \mathrm{Fix}\rho_i^B(\ell_i)$. Let $P_i = \rho_i^B(m_{i, 3}) = \rho_i^B(m_{i, 4})$ with fixed line $\mathcal P_i$; condition (2) of the definition of a four-longitude decomposition forces $\rho_i^B(\ell_{i+2})$ to centralize $P_i$, and so $\mathrm{Fix}(\rho_i^B(\ell_{i+2})) = \mathcal P_i$. By the definition of parabolic form, $(\mathcal P_i, \mathcal L_i)$ is an ordered pair of distinct real lines. 

\begin{definition}\label{def:normalform}
A four-longitude decomposition is in \emph{parabolic normal form along $B$} if, after independent $SL_2(\mathbb{R})$-changes of basis in the two tangle sides adapted to the pair of fixed lines $(\mathcal P_i, \mathcal L_i)$, one has
\begin{equation*}
	\rho_1^B(\ell_1) = \mathrm{red}(A_1) = \epsilon_1\begin{pmatrix}1 & 0 \\ a_1 & 1\end{pmatrix} \ \ \ \ \ \rho_2^B(\ell_2) = \mathrm{red}(B_2) = \epsilon_2\begin{pmatrix}1 & 0 \\ a_2 & 1\end{pmatrix}
\end{equation*}
\begin{equation*}
	\rho_1^B(\ell_3) = \mathrm{red}(A_3) = \epsilon_3\begin{pmatrix}1 & b_1 \\ 0 & 1\end{pmatrix} \ \ \ \ \ \rho_2^B(\ell_4) = \mathrm{red}(B_4) = \epsilon_4\begin{pmatrix}1 & b_2 \\ 0 & 1\end{pmatrix}
\end{equation*}
where $a_1, a_2, b_1, b_2 \in \mathbb R^\times$ and $\epsilon_j \in \{\pm 1\}$. It is in \emph{positive parabolic normal form} if, in addition, $a_1a_2 < 0$. 
\end{definition}

The fact that $J$ is singular is an artifact of the non-realization of $\xi$; this is reframed into a statement about which lines the reduced conjugator collapses. We replace the relative conjugator with a matrix exhibiting this single obstruction as one vanishing parameter $\theta$, which will be helpful in computing trace asymptotics. 

\begin{lemma}\label{lma:antidiag}
In the ordered bases $(\mathcal P_i, \mathcal L_i)$, we have $J = \mathrm{red}(\widehat R) = \omega E_{12}$, for some $\omega \in \mathbb{R}^\times$. After multiplying the normalizers $C_i$ by elements reducing to the identity and rescaling, the relative conjugator can be represented by 
\begin{equation*}
	S_\theta = \begin{pmatrix}0&1 \\ -\theta & 0\end{pmatrix} \ \ \ \ \ \theta \in \mathfrak m_K
\end{equation*}
where $\theta(u) > 0$ for every sufficiently small $u$. 
\end{lemma}

\begin{proof}
By Lemma \ref{lma:rankone}, $\ker J = \mathcal P_2$ and $\im J = \mathcal P_1$. In the ordered bases $(\mathcal P_i, \mathcal L_i)$, this is precisely $J = \omega E_{12}$ for some $\omega \neq 0$, since $\mathcal P_2$, the first basis vector of the domain, is killed, and the image lies along $\mathcal P_1$, the first basis vector of the codomain. Write $\widehat R = \left(\begin{smallmatrix}x&y\\z&w\end{smallmatrix}\right)$. Then $y$ is a unit and $x, z, w \in \mathfrak m_K$. The computation
\begin{equation*}
	\begin{pmatrix}1 & 0 \\ -w/y & 1\end{pmatrix}\widehat R\begin{pmatrix}1 & 0 \\ -x/y & 1\end{pmatrix} = \begin{pmatrix}0 & y \\ -\widehat\Delta/y&0\end{pmatrix}
\end{equation*}
uses matrices which reduce to the identity. Dividing by $y$ gives $S_\theta$, where $\theta = \frac{\widehat\Delta}{y^2}$. Since $\widehat\Delta \in \mathfrak m_K$, one has $\theta \in \mathfrak m_K$, $\mathrm{red}(\theta) = 0$. Moreover, $\widehat\Delta(u)$ and $y(u)^2$ are positive for sufficiently small $u$, and therefore $\theta(u) > 0$ for sufficiently small $u$. 
\end{proof}

Retain the notation $\Sigma_i(t) = C_i\rho_B|_{\pi_1(T_i)}C_i^{-1}$ from Lemma \ref{lma:analyticgauge}, where we now use $S_\theta$ as the relative conjugator. One useful identity, which informs the following asymptotic computations, is 
\begin{equation*}
	S_\theta \begin{pmatrix}1 & 0 \\ a_2 & 1\end{pmatrix} S_\theta^{-1} = \begin{pmatrix}1 & -a_2/\theta \\ 0 & 1\end{pmatrix}
\end{equation*}
which implies the trace identity
\begin{equation*}
	\tr\left(\begin{pmatrix}1&0\\a_1&1\end{pmatrix}S_\theta\begin{pmatrix}1&0\\a_2&1\end{pmatrix}S_\theta^{-1}\right) = 2 - \frac{a_1a_2}{\theta}
\end{equation*}
This is where the sign of $a_1a_2$ enters. 

\begin{lemma}\label{lma:residuehalf}
Write $P_i = \rho_i^B(m_{i,3}) = \left(\begin{smallmatrix}1&p_i\\0&1\end{smallmatrix}\right)$. Set
\begin{equation*}
	K_i(t) = \epsilon_{i+2}\Sigma_i(t)(\ell_{i+2}) \ \ \ \ \ q_i(t) := (K_i(t))_{21}
\end{equation*}
Then there is a nonzero real number $s$ such that $q_i(t) = \frac{2s}{p_i}t + O(t^2)$, and $\theta(t) = \frac{8s^2}{p_1p_2}t^2 + O(t^3)$. Consequently, $p_1p_2 > 0$, and $\frac{q_1(t)q_2(t)}{\theta(t)} = \frac12 + O(t)$. 
\end{lemma}

\begin{proof}
Let $u_i \in T_{\chi_i}X(T_i)$ be the cocycle $u_i(g) = \frac{d}{dt}|_{t=0}\Sigma_i(t)(g)\Sigma_i(0)(g)^{-1}$. By Step 1 of the proof of Theorem \ref{thm:main}, $u_i \in L_i$, and by Corollary \ref{cor:realbranches}, $u_i \neq 0$. Write $\Sigma_i(t)(m_{i, 1}) = I + tX_i + O(t^2)$, so $X_i = u_i(m_{i, 1})$. Thus, from Step 1 of the proof of Theorem \ref{thm:main}, $X_i$ preserves $\mathcal P_i$ and $\mathcal Q_i = \mathrm{Fix}(P_i^{\ell_i^{-1}})$, which are the same lines as in Theorem \ref{thm:main}, since $\tau_i = \ell_i^{-1}$. In the basis where $\rho_i^B(\ell_i) = \left(\begin{smallmatrix}1&0\\a_i&1\end{smallmatrix}\right)$, direct computation shows $\mathcal Q_i = \langle \binom{1}{a_i}\rangle$, a line distinct from $\mathcal L_i$. Tracelessness together with fixing $(1, 0)$ and $(1, a_i)$ forces
\begin{equation*}
	X_i = s_i\begin{pmatrix}1&-2/a_i\\0&-1\end{pmatrix} \ \ \ \ \ X_i^\tau = \Ad_{\rho_i^B(\ell_i^{-1})}X_i = s_i\begin{pmatrix}-1 & -2/a_i \\ 0 & 1\end{pmatrix}
\end{equation*}
for some nonzero $s_i$ (Corollary \ref{cor:realbranches}). The common identity puncture is $m_{1, 2} = m_{2, 1}$. Since $S_\theta$ is the relative conjugator, we have the equation $\Sigma_1(t)(m_{1,2})S_\theta = S_\theta\Sigma_2(t)(m_{2, 1})$. Put
\begin{equation*}
	A(t) = \Sigma_1(t)(m_{1, 2}) \ \ \ \ \ B(t) = \Sigma_2(t)(m_{2, 1}) \ \ \ \ \ S(t) = S_{\theta(t)}
\end{equation*} 
Then we have $A(t)S(t) = S(t)B(t)$. Differentiating at $t = 0$ gives $A'(0)S(0) + A(0)S'(0) = S'(0)B(0) + S(0)B'(0)$. Since $A(0) = B(0) = I$, the two $S'(0)$ terms cancel; using $S(0) = E_{12}$ for what remains gives $X_1^\tau E_{12} = E_{12}X_2$. Comparing the $(1, 2)$ entries yields $s_1 = s_2$. Denote their common nonzero value by $s$. Write
\begin{equation*}
	M_i(t) = \Sigma_i(t)(m_{i, 3}) = (I + tY_i + O(t^2))\begin{pmatrix}1&p_i\\0&1\end{pmatrix} \ \ \ \ \ K_i(t)  = \epsilon_{i+2}\Sigma_i(t)(\ell_{i+2}) = (I + tZ_i + O(t^2))\begin{pmatrix}1&b_i\\0&1\end{pmatrix}
\end{equation*}
(Recall that $m_{i,3}$ and $\ell_{i+2}$ have the same fixed line.) The tangent $u_i$ lies in $L_i \subset \ker\alpha_i = \ker dI_{m_{i, 3}}$, so $\alpha_i(u_i) = \tr(Y_i\left(\begin{smallmatrix}1&p_i\\0&1\end{smallmatrix}\right)) = p_i(Y_i)_{21} = 0$. Since $p_i \neq 0$, it follows that $(Y_i)_{21} = 0$. 

\medskip

The Conway sphere group relations \eqref{eq:boundaryconventions1}-\eqref{eq:boundaryconventions2}, translated to $m$-coordinates, are $m_{1, 1}^{-1}m_{1, 2}m_{1, 3}^{-1}m_{1, 4} = 1$ and $m_{2,1}m_{2,2}^{-1}m_{2,3}m_{2,4}^{-1} = 1$. Taking inverses and substituting $m_{1, 4} = \ell_3^{-1}m_{1,3}\ell_3$, $m_{2,4} = \ell_4^{-1}m_{2,3}\ell_4$, we get
\begin{equation*}
	m_{1,2}^{-1}m_{1,1} = m_{1,3}^{-1}\ell_3^{-1}m_{1,3}\ell_3 \ \ \ \ \ m_{2,2}m_{2,1}^{-1} = m_{2,3}\ell_4^{-1}m_{2,3}^{-1}\ell_4
\end{equation*}
Applying $\Sigma_i(t)$ and using $\Sigma_i(t)(\ell_{i+2}) = \epsilon_{i+2}K_i(t)$:
\begin{equation*}
	\Sigma_1(t)(m_{1,2}^{-1}m_{1,1}) = M_1^{-1}K_1^{-1}M_1K_1 \ \ \ \ \ \Sigma_2(t)(m_{2,2}m_{2,1}^{-1}) = M_2K_2^{-1}M_2^{-1}K_2
\end{equation*}
Write $z_i = (Z_i)_{21}$. Lemma \ref{lma:commutatorasymptotics} gives
\begin{equation*}
	M_i^{-1}K_i^{-1}M_iK_i = I + t\begin{pmatrix}p_iz_i & * \\ 0 & -p_iz_i\end{pmatrix} + O(t^2) \ \ \ \ \ M_iK_i^{-1}M_i^{-1}K_i = I + t\begin{pmatrix}-p_iz_i&*\\0&p_iz_i\end{pmatrix} + O(t^2)
\end{equation*} 
On the other hand, directly from $\Sigma_1(t)(m_{1,1}) = I + t(X_1) + O(t^2)$ and $\Sigma_1(t)(m_{1,2}) = I + tX_1^\tau + O(t^2)$, we get
\begin{equation*}
	\Sigma_1(t)(m_{1,2}^{-1}m_{1,1}) = I + t(X_1 - X_1^\tau) + O(t^2) \ \ \ \ \ \Sigma_2(t)(m_{2,2}m_{2,1}^{-1}) = I + t(X_2^\tau - X_2) + O(t^2)
\end{equation*}
Comparing $(1,1)$ entries with $X_1 - X_1^\tau = \left(\begin{smallmatrix}2s&0\\0&-2s\end{smallmatrix}\right)$ and $X_2^\tau-X_2=\left(\begin{smallmatrix}-2s&0\\0&2s\end{smallmatrix}\right)$ gives $p_i(Z_i)_{21} = 2s$. Thus $q_i(t) = t(Z_i)_{21} + O(t^2) = \frac{2s}{p_i}t+ O(t^2)$, as desired. 

\medskip

For the second-order asymptotic of $\theta$, put 
\begin{equation*}
	M(t) = (I + tY + O(t^2))\begin{pmatrix}1&p\\0&1\end{pmatrix} \ \ \ \ \ K(t) = (I + tZ + O(t^2))\begin{pmatrix}1&b\\0&1\end{pmatrix}
\end{equation*}
with $Y_{21} = 0$. Recall $m_{1, 1} = m_{1,2}(m_{1,3}^{-1}\ell_3^{-1}m_{1,3}\ell_3)$, so 
\begin{equation*}
	\Sigma_1(t)(m_{1,1}) - \Sigma_1(t)(m_{1, 2}) = \Sigma_1(t)(m_{1, 2})(M_1^{-1}K_1^{-1}M_1K_1 - I)
\end{equation*}
Since $X_1^\tau$ has zero $(2, 1)$ entry, $(\Sigma_1(t)(m_{1, 2}))_{21} = O(t^2)$. As $M_1^{-1}K_1^{-1}M_1K_1 - I = O(t)$, the contribution of this entry to the $(2, 1)$ entry of the difference is $O(t^3)$, so only $(M_1^{-1}K_1^{-1}M_1K_1 - I)_{21}$ survives at order $t^2$:
\begin{equation*}
	(\Sigma_1(m_{1,1}))_{21} - (\Sigma_1(m_{1,2}))_{21} = (M_1^{-1}K_1^{-1}M_1K_1)_{21} + O(t^3)
\end{equation*}
Write $A = \Sigma_1(t)(m_{1,3})$ (whose linear term is $Y_1$), $B = \Sigma_2(t)(m_{2, 2})$ (whose linear term is $X^\tau$). The $(1, 2)$ entry of $AS_\theta = S_\theta B$ gives $A_{11} = B_{22}$, hence $(Y_1)_{11} = (X_2^\tau)_{22} = s$. Combined with $(Z_1)_{21} = 2s/p_1$ established above, Lemma \ref{lma:commutatorasymptotics} with $\upsilon = (Y_1)_{11} = s$, $\zeta = (Z_1)_{21} = 2s/p_1$ gives
\begin{equation}\label{eq:first}
	(\Sigma_1(m_{1,1}))_{21} - (\Sigma_1(m_{1,2}))_{21} = -\frac{8s^2}{p_1}t^2 + O(t^3)
\end{equation}
Applying the relative conjugator identity to $m_{1, 1} = m_{2, 4}$ and $m_{1, 2} = m_{2, 1}$ gives 
\begin{equation*}
	\Sigma_1(t)(m_{1, 1}) = S_\theta\Sigma_2(t)(m_{2, 4})S_\theta^{-1} \ \ \ \ \ \Sigma_1(t)(m_{1, 2}) = S_\theta\Sigma_2(t)(m_{2, 1})S_\theta^{-1}
\end{equation*}
For any $B = (b_{jk})$, the lower left entry of $S_\theta BS_\theta^{-1}$ is $-\theta b_{12}$, so
\begin{equation*}
	(\Sigma_1(m_{1,1}))_{21} - (\Sigma_1(m_{1,2}))_{21} = -\theta(t)((\Sigma_2(t)(m_{2, 4}))_{12} - (\Sigma_2(t)(m_{2, 1}))_{12})
\end{equation*}
At $t = 0$, $\rho_2^B(m_{2, 4}) = P_2$ and $\rho_2^B(m_{2, 1}) = I$, so the difference in the right hand side is $p_2 + O(t)$, giving
\begin{equation}\label{eq:second}
	(\Sigma_1(m_{1,1}))_{21} - (\Sigma_1(m_{1,2}))_{21} = -\theta(t)(p_2 + O(t))
\end{equation}
These are two independent computations of the same quantity, arising from the group relator and from the gluing identification respectively. Equating the two expressions \eqref{eq:first}-\eqref{eq:second} gives $\theta(t) = \frac{8s^2}{p_1p_2}t^2 + O(t^3)$. Since $\theta(t) > 0$ for every sufficiently small $t$, $p_1p_2 > 0$. The quotient formula follows from the first-order expressions for $q_1$, $q_2$. 
\end{proof}

\begin{theorem}\label{thm:staggeredtraces}
Suppose $\lambda=\ell_1\ell_2\ell_3\ell_4$ is in positive parabolic normal form along $B$. Then, as $t \to 0$,
\begin{equation*}
	\tr\rho_t(\ell_1\ell_2) = \epsilon_1\epsilon_2\frac{-a_1a_2}{\theta(t)}(1 + O(t)) \ \ \ \ \tr\rho_t(\ell_1\ell_2\ell_3) = \epsilon_1\epsilon_2\epsilon_3\frac{-a_1a_2}{\theta(t)}(1 + O(t))
\end{equation*}
\begin{equation*}
	\tr\rho_t(\lambda) = \epsilon_1\epsilon_2\epsilon_3\epsilon_4\frac{-a_1a_2}{2\theta(t)}(1 + O(t))
\end{equation*}
In particular, $\rho_t(\ell_1\dots\ell_k)$ is hyperbolic for $k = 2, 3, 4$ and all sufficiently small $t$, and $I_\lambda$ has a pole. 
\end{theorem}

\begin{proof}
Set
\begin{equation*}
	\bar A_1 = \epsilon_1A_1 \ \ \ \ \ \bar B_2 = \epsilon_2B_2 \ \ \ \ \ \bar A_3 = \epsilon_3A_3 \ \ \ \ \ \bar B_4 = \epsilon_4B_4
\end{equation*}
where $A_1$, $B_2$, $A_3$, $B_4$ are as in Definition \ref{def:fourlongitude}. Let $\bar A_1 = \left(\begin{smallmatrix}a&b\\c&d\end{smallmatrix}\right)$ and $\bar B_2 = \left(\begin{smallmatrix}e&f\\g&h\end{smallmatrix}\right)$. By Lemma \ref{lma:analyticgauge} and the definition of positive parabolic normal form, as $t \to 0$, we have $c \to a_1$, $g \to a_2$, $a, d, e, h \to 1$, and $b, f \to 0$. For arbitrary matrices $A = (a_{jk})$ and $B = (b_{jk})$, one has the exact identity 
\begin{equation}\label{eq:twofactor}
	\tr(AS_\theta BS_\theta^{-1}) = A_{11}B_{22} + A_{22}B_{11} - \theta A_{12}B_{12} - \theta^{-1}A_{21}B_{21}
\end{equation}
Applying this to $\bar A_1$, $\bar B_2$ and using $cg \to a_1a_2$, $ah + de - \theta bf = O(1)$, together with $\theta \sim t^2$ (Lemma \ref{lma:residuehalf}) gives $\mathrm{tr}(\bar A_1S_\theta\bar B_2S_\theta^{-1}) = \frac{-a_1a_2}{\theta(t)}(1 + O(t))$. Restoring the central signs to the first asymptotic of \eqref{eq:firstasymptotic}, with $S_\theta$ playing the role of $\widehat R$, proves the first formula. 

\medskip

Let $N := \bar A_1S_\theta\bar B_2S_\theta^{-1}$. Direct multiplication gives
\begin{equation*}
	N = \begin{pmatrix}ah-b\theta f & -ag/\theta + be \\ ch-d\theta f & -cg/\theta + de\end{pmatrix}
\end{equation*}
Consequently, $N_{11}, N_{21} = O(1)$, $N_{12} = O(\theta^{-1})$, and $N_{22} = -a_1a_2\theta^{-1}(1+O(t))$. Let $\bar A_3 = \left(\begin{smallmatrix}A&B\\q_1&D\end{smallmatrix}\right)$; then $\bar A_3 \to \left(\begin{smallmatrix}1&b_1\\0&1\end{smallmatrix}\right)$, so $A, D = 1  + O(t)$, $B = O(1)$, and $q_1 = O(t)$ by Lemma \ref{lma:residuehalf}. The trace is
\begin{equation*}
	\tr(\bar A_1S_\theta\bar B_2S_\theta^{-1}\bar A_3) = N_{11}A + N_{12}q_1 + N_{21}B + N_{22}D
\end{equation*}
The terms $N_{11}A$ and $N_{21}B$ are $O(1)$, negligible against terms diverging like $\theta^{-1}$. Since $N_{12} = O(\theta^{-1})$, $q_1 = O(t)$, and $\theta \sim t^2$ by Lemma \ref{lma:residuehalf}, $N_{12}q_1 = O(t/\theta) = O(t^{-1})$, a factor of $t$ smaller than $N_{22}D = -a_1a_2\theta^{-1}(1 + O(t))$, and so it is absorbed into $(1 + O(t))$. Thus,
\begin{equation*}
	\tr(\bar A_1S_\theta\bar B_2S_\theta^{-1}\bar A_3) = \frac{-a_1a_2}{\theta}(1 + O(t))
\end{equation*}
This, along with the second asymptotic of \eqref{eq:firstasymptotic}, proves the second formula after restoring $\epsilon_1\epsilon_2\epsilon_3$.

\medskip

For the full product, write $\bar B_4=\left(\begin{smallmatrix}E&F\\q_2&H\end{smallmatrix}\right)$. A direct multiplication gives
\begin{equation}\label{eq:fourfactor}
	\mathrm{tr}(\bar A_1S_\theta\bar B_2S_\theta^{-1}\bar A_3S_\theta\bar B_4S_\theta^{-1}) = \frac{c g q_1q_2}{\theta^2} -\frac{Aq_2ch+q_1q_2de+q_1Hag+DEcg}{\theta}+O(1)                                      
\end{equation}
Here, all entries other than $q_1,q_2$ have finite limits. Lemma \ref{lma:residuehalf} gives
$\frac{q_1q_2}{\theta}=\frac12+O(t)$. Hence the first term in \eqref{eq:fourfactor} is $a_1a_2(2\theta)^{-1}(1+O(t))$, the last displayed numerator contributes $-a_1a_2\theta^{-1}(1+O(t))$, and the remaining negative-power terms are $o(\theta^{-1})$. Therefore $\operatorname{tr}(\bar A_1S_\theta\bar B_2S_\theta^{-1}\bar A_3S_\theta\bar B_4S_\theta^{-1})=\frac{-a_1a_2}{2\theta}(1+O(t))$. Restoring the four central signs and \eqref{eq:thirdasymptotic} proves the longitude formula. Since $a_1a_2<0$ and $\theta(t)>0$, all three displayed traces have unbounded absolute value, so for $t$ sufficiently small, all three traces approach $\sigma_k\infty$, where $\sigma_k = \epsilon_1\dots\epsilon_k$. In particular, $I_\lambda$ has a pole at the ideal point.
\end{proof}

\subsection{Translation numbers of longitudes} 

Throughout the remainder of this subsection, $\widetilde\rho_u$ denotes the meridian-normalized lift of Lemma \ref{lma:normalizedlift}.

\begin{lemma}\label{lma:heightlimit}
For any $\gamma \in \pi_1(T_1) \cup \pi_1(T_2)$, the limit $\lim_{u \to 0^+}\mathrm{trans}(\widetilde\rho_u(\gamma))$ exists. 
\end{lemma}

\begin{proof}
Consider $\gamma \in \pi_1(T_i)$. The path $\Sigma_i(u)(\gamma) =  C_i(u)\rho_u(\gamma)C_i(u)^{-1}$ extends continuously to $u = 0$ with limit $\rho_i^B(\gamma)$. Since $(0, \delta)$ is contractible, the image of $C_i(u)$ in $PSL_2(\mathbb{R})$ has a continuous lift $\widetilde{C}_i(u)$ to $\widetilde{PSL}_2(\mathbb{R})$. Thus $u \mapsto \widetilde{C}_i(u)\,\widetilde\rho_u(\gamma)\widetilde{C}_i(u)^{-1}$ is then a continuous lift of the image in $PSL_2(\mathbb{R})$. For all sufficiently small $u$, this path lies in an evenly covered neighborhood of its limit in $PSL_2(\mathbb{R})$, so this lift differs from a convergent local lift by a fixed central deck transformation. Thus it converges. Translation number is continuous and conjugacy invariant, proving the claim.
\end{proof}

In a four-longitude decomposition $\lambda = \ell_1\ell_2\ell_3\ell_4$ in parabolic form, let
\begin{equation*}
	h_j = \lim_{u \to 0^+}\mathrm{trans}(\widetilde\rho_u(\ell_j)) \in \mathbb Z
\end{equation*}
These translation numbers are integers because their limiting images are noncentral parabolics.

\begin{lemma}\label{lma:trans}
Let $\lambda=\ell_1\ell_2\ell_3\ell_4$ be a four-longitude decomposition in positive parabolic normal form along $B$. Then for every sufficiently small $t \neq 0$, 
\begin{equation*}
	\mathrm{trans}(\widetilde\rho_t(\ell_1\ell_2)) = h_1 + h_2 \ \ \ \ \ \mathrm{trans}(\widetilde\rho_t(\ell_1\ell_2\ell_3)) = h_1 + h_2 + h_3
\end{equation*}
\begin{equation*}
	\mathrm{trans}(\widetilde\rho_t(\lambda)) = h_1+h_2+h_3+h_4
\end{equation*}
\end{lemma}

\begin{proof}
Since the limiting image of $\ell_j$ is a noncentral parabolic of trace $2\epsilon_j$, Lemma \ref{lma:parity} gives $\epsilon_j = (-1)^{h_j}$. For $k = 2$, for $t$ sufficiently small, continuity of translation number gives
\begin{equation*}
	|\mathrm{trans}(\widetilde\rho_t(\ell_j)) - h_j| < \frac14 \ \ \ \ \ j = 1, 2
\end{equation*}
The defect-one quasimorphism property of translation numbers then implies, for $t$ sufficiently small,
\begin{equation}\label{eq:diff2}
	|\mathrm{trans}(\widetilde\rho_t(\ell_1\ell_2)) - (h_1 + h_2)| < \frac32 < 2
\end{equation} 
By Theorem \ref{thm:staggeredtraces}, for small $t$, $\rho_t(\ell_1\ell_2)$ is hyperbolic, so its translation number is an integer; in addition, by positive parabolic normal form, $-a_1a_2 > 0$, for $t$ sufficiently small, $\rho_t(\ell_1\dots\ell_k)$ is hyperbolic with the same sign as $\epsilon_1\dots\epsilon_k = (-1)^{h_1+\dots+h_k}$, for $k = 2, 3, 4$. By Lemma \ref{lma:parity}, $\mathrm{trans}(\widetilde\rho_t(\ell_1\ell_2))$ and $h_1 + h_2$ have the same parity, so by \eqref{eq:diff2} they are equal. Applying the same argument successively to $(\ell_1\ell_2)\ell_3$ and $(\ell_1\ell_2\ell_3)\ell_4$ gives the result.
\end{proof}

\begin{theorem}\label{thm:h00}
Let $B$ be a real half-branch as in Proposition \ref{prop:branchmodel}, with associated representations $\rho_u$ and meridian-normalized lifts $\widetilde\rho_u$. Suppose the homological longitude admits a four-longitude decomposition $\lambda=\ell_1\ell_2\ell_3\ell_4$ in positive parabolic normal form along $B$. Then for $u$ sufficiently small, the expanding eigenline $v_u\in\mathbb RP^1$ of $\rho_u(\mu)$ is well-defined and depends continuously on $u$. The pairs $(\widetilde\rho_u,v_u)$ form a continuous family of augmented peripherally hyperbolic representations, and $\operatorname{EV}(\widetilde\rho_u,v_u) = \bigl(x(u),y(u);0,h_1+h_2+h_3+h_4\bigr)$, where $x(u)=\log|a_\mu(u)|$ and $y(u)=\log|a_\lambda(u)|$, and $a_\mu(u),a_\lambda(u)$ are the eigenvalues of $\rho_u(\mu),\rho_u(\lambda)$ on $v_u$. Consequently, the image of this family is an augmented holonomy arc contained in $H_{0,h_1+h_2+h_3+h_4}(M)$.
\end{theorem}

\begin{proof}
Lemma \ref{lma:trans} gives $\mathrm{trans}(\widetilde\rho_u(\lambda)) = h_1+h_2+h_3+h_4$. By meridional normalization, $\operatorname{trans}\bigl(\widetilde\rho_u(\mu)\bigr)=0$. By Corollary \ref{cor:realbranches}, $\operatorname{tr}\rho_u(\mu)>2$ for all sufficiently small $u>0$. Hence $\rho_u(\mu)$ has two distinct positive eigenvalues $a_\mu(u)>1$ and $a_\mu(u)^{-1}<1$, and its expanding eigenline $v_u$ is simple and depends continuously on $u$. The peripheral elements $\mu$ and $\lambda$ commute, so $\rho_u(\lambda)$ preserves each of the two eigenspaces of $\rho_u(\mu)$, and in particular preserves $v_u$. Thus, the whole peripheral representation fixes $v_u$. Moreover, the centralizer of a hyperbolic element in $PSL_2(\mathbb R)$ is the split real torus preserving its two fixed points, so the restriction of $\widetilde\rho_u$ to $\pi_1(\partial M)$ is hyperbolic in the sense used in the definition of the holonomy extension locus. Consequently, $(\widetilde\rho_u,v_u)$ is an augmented peripherally hyperbolic representation. Let $a_\lambda(u)$ denote the eigenvalue of $\rho_u(\lambda)$ on $v_u$. By the definition of the eigenvalue-translation map, we have
\begin{equation*}
	\operatorname{EV}(\widetilde\rho_u,v_u)(\mu) = \left(\log|a_\mu(u)|,\operatorname{trans}\bigl(\widetilde\rho_u(\mu)\bigr)\right) = \bigl(x(u),0\bigr)
\end{equation*}
\begin{equation*}
	\operatorname{EV}(\widetilde\rho_u,v_u)(\lambda) = \left(\log|a_\lambda(u)|, \operatorname{trans} \bigl(\widetilde\rho_u(\lambda)\bigr)\right) = \bigl(y(u),h_1+h_2+h_3+h_4\bigr)
\end{equation*}
Therefore every point of the resulting augmented holonomy arc lies in $H_{0,h_1+h_2+h_3+h_4}(M)$.
\end{proof}

\begin{definition}\label{def:heightcancellation}
Let $B$ be a real half-branch and $\lambda = \ell_1\ell_2\ell_3\ell_4$ a four-longitude decomposition. We say that $(\lambda, B)$ satisfies \emph{height cancellation} if $h_1 + h_2 + h_3 + h_4 = 0$.
\end{definition}

\begin{corollary}\label{cor:h00}
Suppose a four-longitude decomposition is in positive parabolic normal form on $B^+$ and $B^-$ and satisfies height cancellation. Then both associated augmented holonomy arcs lie in $H_{0,0}(M)$.
\end{corollary}

\begin{proof}
Apply Theorem \ref{thm:h00} separately to $B^+$ and $B^-$. Height cancellation is exactly $h_1+h_2+h_3+h_4=0$. Theorem \ref{thm:h00} therefore places both associated arcs in $H_{0,0}(M)$.
\end{proof}

\subsection{Proof of Theorem \ref{mainthm:b}}

We next show that the two augmented arcs from Theorem \ref{thm:h00} and Corollary \ref{cor:h00} are unbounded with vertical asymptotes on opposite sides of the ideal point. In $\pi_1(M)$, $\mu_{1,1}$ is conjugate to the knot meridian $\mu$, so their trace functions agree. Along the signed normalization, $I_\mu(t)=2+Qt^2+O(t^3)$, $Q>0$ (Corollary \ref{cor:realbranches}). If positive parabolic normal form holds on both half-branches, Theorem \ref{thm:staggeredtraces}, applied separately to the two sides, shows that $I_\lambda$ has a pole at $x$ and $|I_\lambda(t)|\to \infty$ as $t\to0$ from either side.

\begin{lemma}\label{lma:eigenvalues}
There exist real analytic functions $m_+(t), m_-(t)$ near $t = 0$ such that, for $t \neq 0$, they are the two eigenvalues of $\rho_t(\mu)$, and 
\begin{equation*}
	m_+(t) = 1 + \sqrt{Q}t + O(t^2) \ \ \ \ \ m_-(t) = 1 - \sqrt{Q}t + O(t^2)
\end{equation*}
Moreover, define 
\begin{equation*}
	l_+(t) = \frac{I_{\mu\lambda}(t) - m_-(t)I_\lambda(t)}{m_+(t) - m_-(t)} \ \ \ \ \ l_-(t) = \frac{m_+(t)I_\lambda(t) - I_{\mu\lambda}(t)}{m_+(t) - m_-(t)}
\end{equation*}
Then $l_+$ and $l_-$ are meromorphic functions at $t = 0$. For $t \neq 0$, $l_{\pm}(t)$ are the eigenvalues of $\rho_t(\lambda)$ on the $m_{\pm}(t)$-eigenlines of $\rho_t(\mu)$. In particular, we have $l_+(t) + l_-(t) = I_\lambda(t)$ and $l_+(t)l_-(t) = 1$.
\end{lemma}

\begin{proof}
From the quadratic expansion of $I_\mu(t)$, we have $I_\mu(t)^2 - 4 = 4Qt^2 + O(t^3) = t^2d(t)$ where $d(t)$ is real analytic, and $d(0) = 4Q > 0$. After shrinking the interval, choose the positive real analytic square root $h(t)$ of $d(t)$, so $h(0) = 2\sqrt{Q}$. Set
\begin{equation*}
	m_\pm(t) = \frac{I_\mu(t)\pm th(t)}{2} \Longrightarrow m_\pm(t)^2 - I_\mu(t)m_\pm(t) + 1 = 0
\end{equation*}
Hence, for $t \neq 0$, these are the two eigenvalues of $\rho_t(\mu)$. Their quadratic expansions are 
\begin{equation*}
	m_+(t) = 1 + \sqrt{Q}t + O(t^2) \ \ \ \ \ m_-(t) = 1 - \sqrt{Q}t + O(t^2)
\end{equation*}
as desired. For $t\neq0$, the meridian has distinct eigenvalues. Since $\mu$ and $\lambda$ commute in the peripheral subgroup, $\rho_t(\lambda)$ preserves the two meridian eigenlines. Let its eigenvalues on these lines be $\eta_+(t)$ and $\eta_-(t)$. Then
\begin{equation*}
	I_\lambda(t) = \eta_+(t) + \eta_-(t) \ \ \ \ \ I_{\mu\lambda}(t) = m_+(t)\eta_+(t)+ m_-(t)\eta_-(t)
\end{equation*}
Solving this linear system gives 
\begin{equation*}
	\eta_+(t) = \frac{I_{\mu\lambda}(t)-m_-(t)I_\lambda(t)}{m_+(t)-m_-(t)} \ \ \ \ \ \eta_-(t)
	= \frac{m_+(t)I_\lambda(t)-I_{\mu\lambda}(t)}{m_+(t)-m_-(t)}
\end{equation*}
Thus $\eta_\pm=l_\pm$. Since the trace functions are rational on $\widetilde{\mathcal C}$ and hence meromorphic at the ideal point, and $m_+(t) - m_-(t) = th(t)$, the functions $l_\pm$ are meromorphic at $t=0$. Finally, we have $l_+(t)l_-(t) = \det\rho_t(\lambda) = 1$ for $t \neq 0$, and hence as an identity of meromorphic germs. 
\end{proof}

\begin{theorem}\label{thm:asymptote}
Let $\nu: (-\epsilon, \epsilon) \to \widetilde{\mathcal{C}}(\mathbb{R})$ be the signed real normalization at the ideal point $x$, and let $\chi_t$ be the corresponding characters for $t \neq 0$. Suppose positive parabolic normal form holds on both half-branches. For each $t \neq 0$, choose the common peripheral eigenline $v_t$ on which the meridian eigenvalue has absolute value greater than one, and let $m(t)$ and $l(t)$ be the eigenvalues of $\rho_t(\mu)$ and $\rho_t(\lambda)$ on $v_t$. By the definition of the eigenvalue translation map, the corresponding eigenvalue coordinates are $x(t)=\log|m(t)|$ and $y(t)=\log|l(t)|$. Then $x(t) = \sqrt Q|t| + O(t^2)$; in particular, $x(t) > 0$ for small $t \neq 0$, and there is a sign $\sigma \in \{\pm 1\}$ such that $y(t) \to \sigma\infty$ as $t \to 0^+$, and $y(t) \to -\sigma\infty$ as $t \to 0^-$. Consequently, $-\frac{y(t)}{x(t)}$ tends to opposite infinities on the two half-branches.
\end{theorem}

\begin{proof}
Let $m_\pm(t)$ and $l_\pm(t)$ be the meromorphic peripheral eigenvalue branches constructed in Lemma \ref{lma:eigenvalues}. For $t > 0$ sufficiently small, $m_+(t) > 1 > m_-(t) > 0$, while for $t < 0$ sufficiently small, $m_-(t) > 1 > m_+(t) > 0$, both positive since $I_\mu(t) > 2$. Thus, we have 
\begin{equation*}
	x(t) = \begin{cases}\log m_+(t) & t > 0 \\ \log m_-(t) & t < 0\end{cases} 
\end{equation*}
The expansions of $m_\pm$ from Lemma \ref{lma:eigenvalues} give $x(t) = \sqrt{Q}|t| + O(t^2)$; in particular, $x(t) > 0$ for small $t \neq 0$. 

\medskip

Since $l_+l_- = 1$, the orders at $t = 0$ satisfy $\mathrm{ord}_0(l_+) + \mathrm{ord}_0(l_-) = 0$. On the other hand, $I_\lambda = l_+ + l_-$ has a pole. Therefore, exactly one of $l_+$ and $l_-$ has a pole at $t = 0$, and the other has a zero of the same order. Suppose first that $l_+$ is the pole eigenvalue. Then 
\begin{equation*}
	|l_+(t)| \to +\infty \ \ \ \ \ |l_-(t)| = |l_+(t)|^{-1} \to 0
\end{equation*}
from both sides. Since the meridian-expanding line is the $m_+(t)$ eigenline for $t > 0$ and the $m_-(t)$ eigenline for $t < 0$, and since the longitude coordinate is defined using the absolute value of its eigenvalue, one has
\begin{equation*}
	y(t) = 
	\begin{cases} \log |l_+(t)|&t>0 \\ \log |l_-(t)|&t<0\end{cases}
\end{equation*}
Consequently, $\lim_{t \to 0^+}y(t) = +\infty$ and $\lim_{t \to 0^-}y(t) = -\infty$. If instead $l_-(t)$ were the pole branch, the same argument would give $\lim_{t \to 0^+}y(t) = -\infty$ and $\lim_{t \to 0^-}y(t) = +\infty$. This proves the asserted statement for some $\sigma\in\{\pm1\}$. Thus, $-\frac{y(t)}{x(t)}$ tends to opposite infinities on the two half-branches.
\end{proof}

\begin{corollary}\label{cor:sheets}
Suppose positive parabolic normal form holds on both half-branches. Then the augmented holonomy arc associated to $B$ is unbounded with vertical asymptote and is contained in $H_{0,h_1+h_2+h_3+h_4}(M)$; the two arcs have vertical asymptotes on opposite sides.
\end{corollary}

\begin{proof}
Combine Theorems \ref{thm:h00} and \ref{thm:asymptote}, applied to each half-branch separately. 
\end{proof}

\begin{corollary}\label{cor:slopes}
Assume that positive parabolic normal form and height cancellation hold on both half-branches. Then there exist $R_-<R_+$ such that, for every rational slope $r<R_-$ or $r>R_+$, the group $\pi_1(M(r))$ is left-orderable.
\end{corollary}

\begin{proof}
For $r>0$, parametrize the half-branch on which $y\to-\infty$ by $u\in(0,\delta)$, with $u\to0^+$ at the ideal point, and put $F_r(u)=y(u)+rx(u)$. Choose $u_0>0$ sufficiently small so that the conclusions of Theorem \ref{thm:h00} and Corollary \ref{cor:h00} hold on $(0, u_0)$; the associated augmented arc over $(0, u_0)$ lies in $H_{0,0}(M)$. Since $x(u_0)>0$, one has $F_r(u_0)>0$ for all sufficiently large $r$, whereas $F_r(u)\longrightarrow-\infty$ as $u\to0$. The intermediate value theorem therefore gives $u_r\in(0,u_0)$ such that $F_r(u_r)=0$. Equivalently, $(x(u_r),y(u_r))$, being the $EV$-image of $(\widetilde\rho_{u_r}, v_{u_r})$ with $u_r \in (0, u_0)$, lies in $L_r\cap H_{0,0}(M)$ (here $L_r$ is the filling line of slope $r$, in the notation of Section \ref{sec:h00}). This point is nonzero because $x(u_r)>0$, and it is nonideal because it is the eigenvalue-translation data of the actual augmented representation $(\widetilde\rho_{u_r},v_{u_r})$, rather than a point added only in taking the closure of the holonomy extension locus.

\medskip
	
The same argument on the other half-branch gives such an intersection for every sufficiently large negative $r$. Since $M$ is the exterior of a nontrivial knot in $S^3$, it is irreducible. By Theorem 1.2 of Gordon-Luecke \cite{gordonleucke}, at most three fillings of $M$ are reducible. Enlarging $[R_-,R_+]$ to contain these finitely many slopes, Lemma 3.8 of \cite{gao} applies to every rational $r \notin [R_-, R_+]$: the intersection point produced above is nonzero and nonideal, and $M(r)$ is irreducible, so $\pi_1(M(r))$ is left-orderable. 
\end{proof}

Corollaries \ref{cor:h00} and \ref{cor:slopes} prove Theorem \ref{mainthm:b}. 

\subsection{A criterion from chosen limiting representatives}

The definition of positive parabolic normal form uses the branch-adapted limiting representatives. In the applications, the limiting representations are independently chosen. We provide a practical criterion, based on independently chosen limiting representations which share a trace with the branch-adapted ones, to certify that positive parabolic normal form and height cancellation hold. 

\medskip

Fix irreducible real representatives $\rho_i: \pi_1(T_i) \to SL_2(\mathbb R)$ of the limiting characters. Suppose that, after adapting to the fixed lines of the Conway parabolic and the identity-strand piece, the four-longitude decomposition $\lambda = \ell_1\ell_2\ell_3\ell_4$ satisfies
\begin{equation*}
	\rho_i(m_{i, 3}) = \begin{pmatrix}1 & c_i\\0&1\end{pmatrix} \ \ \ \ \ \rho_i(\ell_i)=\epsilon_i\begin{pmatrix}1 & 0 \\ \alpha_i & 1\end{pmatrix} \ \ \ \ \ \rho_i(\ell_{i+2}) = \epsilon_{i+2}\begin{pmatrix}1&\beta_i\\0&1\end{pmatrix} \ \ \ \ \ i = 1, 2
\end{equation*}
where $c_i, \alpha_i, \beta_i \in \mathbb R^\times$, and $\epsilon_j \in \{\pm 1\}$. 

\medskip

Let $\widetilde\rho_i: \pi_1(T_i) \to \widetilde{PSL}_2(\mathbb R)$ be a meridian-normalized lift of $\rho_i$, with $\mathrm{trans}(\widetilde\rho_i(m_{i,1})) = \mathrm{trans}(\widetilde\rho_i(m_{i,3})) = 0$. The two strand meridian classes generate $H_1(T_i; \mathbb{Z})$, so this normalization determines the lift uniquely. Let 
\begin{equation*}
	\eta_i = \mathrm{trans}(\widetilde\rho_i(\ell_i)) \ \ \ \ \ \eta_{i+2} = \mathrm{trans}(\widetilde\rho_i(\ell_{i+2})) \ \ \ \ \ i = 1, 2
\end{equation*}

\begin{theorem}\label{thm:modelcriterion}
Under the above hypotheses, the four-longitude decomposition is in positive parabolic normal form along one, and hence both, real half-branches if and only if $c_1c_2\alpha_1\alpha_2 < 0$. Moreover, for each real half-branch $B$, there is a sign $\sigma_B \in \{\pm 1\}$ such that its limiting height vector is 
\begin{equation*}
	(h_1, h_2, h_3, h_4) = \sigma_B(\eta_1, \mathrm{sgn}(c_1c_2)\eta_2, \eta_3, \mathrm{sgn}(c_1c_2)\eta_4)
\end{equation*}
Consequently, height-cancellation holds along one, and hence both, real half-branches if and only if $\eta_1+\eta_3 + \mathrm{sgn}(c_1c_2)(\eta_2 + \eta_4) = 0$. 
\end{theorem}

\begin{proof}
Fix a half-branch $B$. Let $\rho_i^B$ be its branch-adapted limiting representatives, and let  
\begin{equation*}
	\rho_i^B(m_{i, 3}) = \begin{pmatrix}1&p_i\\0&1\end{pmatrix} \ \ \ \ \ \rho_i^B(\ell_i) = \epsilon_i\begin{pmatrix}1&0\\a_i&1\end{pmatrix}
\end{equation*}
Choose a real conjugator $D_i \in GL_2(\mathbb R)$ such that $\rho_i^B = D_i\rho_iD_i^{-1}$, and put $d_i = \mathrm{sgn}\det D_i$. An orientation-preserving change of basis would rescale the upper right entries by a positive square, while an orientation-reversing change negates the sign. Hence
\begin{equation}\label{eq:signs}
	\mathrm{sgn}(p_i) = d_i\mathrm{sgn}(c_i) \ \ \ \ \ \mathrm{sgn}(a_i) = d_i\mathrm{sgn}(\alpha_i)
\end{equation}
In particular, $\mathrm{sgn}(p_i a_i) = \mathrm{sgn}(c_i\alpha_i)$. Lemma \ref{lma:residuehalf} gives $p_1p_2 > 0$. Therefore, 
\begin{equation*}
	\mathrm{sgn}(a_1a_2) = \mathrm{sgn}(p_1p_2a_1a_2) = \mathrm{sgn}((p_1a_1)(p_2a_2)) = \mathrm{sgn}(c_1c_2\alpha_1\alpha_2)
\end{equation*}
Positive parabolic normal form is exactly the condition $a_1a_2 < 0$, proving the first assertion. The right-hand side is independent of $B$, so the same criterion holds on both half-branches. For the height statement, let $\widetilde\rho_i^B$ be the limiting lift obtained from restriction to $\pi_1(T_i)$. The comparison determined by $D_i$ preserves translation number when $d_i = 1$. When $d_i = -1$, it uses the lifted outer involution of $\widetilde{PSL}_2(\mathbb R)$ and negates translation number. Thus, conjugating $\widetilde\rho_i$ by $D_i$ multiplies the translation number by $d_i$. Both the lift of $D_i\rho_iD_i^{-1}$ and $\widetilde\rho_i^B$ have translation number zero on the two strand-meridian generators, which generate $H_1(T_i; \mathbb Z)$ so they must agree. Thus, we have $h_i = d_i\eta_i$ and $h_{i+2} = d_i\eta_{i+2}$. Equation \eqref{eq:signs} and $p_1p_2 > 0$ give $d_1d_2 = \mathrm{sgn}(c_1c_2)$. Taking $\sigma_B = d_1$, we have $d_2 = \sigma_B\mathrm{sgn}(c_1c_2)$, and hence
\begin{equation*}
	(h_1, h_2, h_3, h_4) = \sigma_B(\eta_1, \mathrm{sgn}(c_1c_2)\eta_2, \eta_3, \mathrm{sgn}(c_1c_2)\eta_4)
\end{equation*}
Summing the four entries proves the height-cancellation criterion. 
\end{proof}

%% file: twobridge.tex
The constructions of Sections \ref{sec:idealpoints} and \ref{sec:trans} utilize various data from tangle exteriors. In this section we examine two infinite families of tangle exteriors: the \emph{torus-trivial} and \emph{twist-trivial} tangles. We compute the associated data needed to apply the results of the previous two sections. 

\medskip

For each type of tangle, we verify the four local hypotheses of Theorem \ref{thm:main}. We also choose representations which will form real admissible limiting gluing data and two strand traversals. We compute their central signs, parabolic parameters, and translation numbers. These will provide model tangles and representations which will be matched together in Section \ref{sec:families} to form the knots admitting four-longitude decompositions in positive parabolic normal form. 

\medskip

Subsections \ref{subsec:torus-trivial} and \ref{subsec:twist-trivial} follow the same sequence: presentation, quotient, real models, certification conditions, peripheral data, and heights. Subsection \ref{subsec:mirrors} extends all of these to the mirrors and flips, and the output of the section is assembled in Table \ref{tab:decorations}. We conclude the section by showing that the Conway sphere in any knot formed by gluing together two torus- or twist-trivial tangles has an essential Conway sphere, and we prove that all such knots satisfy the hypotheses of Theorem \ref{mainthm:a}, resulting in a precise count of ideal points and real branches originating from the construction of Section \ref{sec:idealpoints}.

\medskip

Condition (1) requires no per-tangle computation: every tangle group in this paper is free of rank 2, and for $F_2 = \langle x, y \rangle$, the trace map $\chi \mapsto (\chi(x), \chi(y), \chi(xy))$ is an isomorphism $X(F_2) \cong \mathbb{C}^3$. This is a classical theorem of Fricke and Vogt; see Theorem A of \cite{goldman}. For the two quotient groups below, condition (2) follows from the standard local descriptions of their nonabelian character varieties.

\subsection{Torus-trivial tangles}\label{subsec:torus-trivial}

\begin{figure}[h]
	\centering
	\includegraphics[scale=.6]{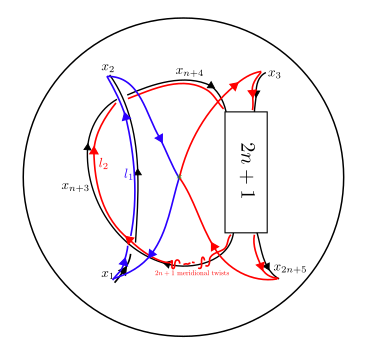}
	\caption{The torus-trivial tangle $T_n^{tor}$; the box stands for $2n+1$ right-handed crossings. The trivial strand has endpoints labeled $x_1, x_2$. The knotted strand, with endpoints $x_3, x_{2n+5}$, closes up to the $(2, 2n+1)$ torus knot. The labeled arcs are the Wirtinger generators appearing in the clasp relations \eqref{eq:type1initial}, together with the four Conway meridians $m_1 = x_1, m_2 = x_2, m_3 = x_3, m_4 = x_{2n+5}$. The two colored curves are as follows: $l_1$ (blue) traverses the trivial strand from $x_1$ to $x_2$, and $l_2$ (red) is a parallel of the knotted strand, corrected by $2n+1$ left-handed meridional twists.}\label{fig:torustrivial}
\end{figure}

Our first family is the \emph{torus-trivial tangles} $T_n^{tor}$, for $n \geq 1$. A generic torus-trivial tangle is pictured in Figure \ref{fig:torustrivial}.

Using the Wirtinger presentation, we compute the fundamental group. The diagram has $2n+5$ arcs; denote the corresponding generators $x_1, \dots, x_{2n+5}$. The relators are then
\begin{equation}\label{eq:type1initial}
	x_2 = x_{n+3}^{-1}x_1x_{n+3} \ \ \ \ \ x_{n+4} = x_2^{-1}x_{n+3}x_2
\end{equation}
For $1 \leq k \leq n$:
\begin{equation}\label{eq:type1ladder1}
	x_{k+3} = x_{n+k+4}x_{k+2}x_{n+k+4}^{-1}
\end{equation}
For $0 \leq k \leq n$:
\begin{equation}\label{eq:type1ladder2}
	x_{n+k+5} = x_{k+3}x_{n+k+4}x_{k+3}^{-1}
\end{equation}
Denote this group $G_n^{tor}$.

\begin{prop}\label{prop:torusgroup}
$G_n^{tor}$ is free on $A = x_{n+3}, H = x_2(x_{n+3}x_{2n+4})^{-n}$. Killing the trivial-strand meridian gives the quotient
\begin{equation*}
	Q_n^{tor} \cong \langle A, B \mid (BA)^nB = A(BA)^n \rangle \cong \langle u, v \mid u^2 = v^{2n+1} \rangle
\end{equation*}
which is the group of the $(2, 2n+1)$-torus knot, denoted $T_{2,2n+1}$. 
\end{prop}

\begin{proof}
Set $A_k = x_{k+3}$ and $B_k = x_{n+k+4}$. Then relations \eqref{eq:type1ladder1} and \eqref{eq:type1ladder2} imply that for any $1 \leq k \leq n$, $A_k = B_kA_{k-1}B_k^{-1}$ and $B_{k+1}=A_kB_kA_k^{-1}$, and so $A_kB_k=B_kA_{k-1}$; applying this one index lower, $B_kA_{k-1} = A_{k-1}B_{k-1}$, so $A_kB_k = A_{k-1}B_{k-1}$. Thus, $A_kB_k$ is actually independent of $k$; call it $D = A_kB_k$. Let $A = A_n = x_{n+3}$. Since $B_k = A_k^{-1}D$ and $A_{k-1} = D^{-1}A_kD$, downward induction gives $A_k = D^{-(n-k)}AD^{n-k}$, and $B_k=D^{-(n-k)}A^{-1}D^{n-k+1}$. In particular, $B_0 = D^{-n}A^{-1}D^{n+1}$. Define $H = x_2D^{-n}$, so that $x_2 = HD^n$. The second relation of \eqref{eq:type1initial} gives $B_0 = x_2^{-1}Ax_2 = D^{-n}H^{-1}AHD^n$. Comparing with the formula for $B_0$ and cancelling the outer power of $D$ yields $A^{-1}D = H^{-1}AH$, or equivalently $D = AH^{-1}AH$. The first relation of \eqref{eq:type1initial} is $x_2 = A^{-1}x_1A$, and hence $x_1=Ax_2A^{-1}=AHD^nA^{-1}$. Thus, every Wirtinger generator is a word in $A, H, D$, and $D = AH^{-1}AH$ eliminates $D$. Thus $G_n^{tor} \cong \langle A, H \rangle$ is free of rank two.

\medskip

Killing the trivial strand is equivalent to setting $x_1 = x_2 = 1$. Hence $H=D^{-n}$. The equation $D=AH^{-1}AH$ becomes $D=AD^nAD^{-n}$, or $D^{n+1} = AD^nA$. Put $B = HAH^{-1}=D^{-n}AD^n$. Then $BA=D^{-n}AD^nA=D$, and by definition, $D^nB=AD^n$. Since $D=BA$, we have $(BA)^nB=A(BA)^n$, implying that $Q_n^{tor}\cong\langle A,B\mid(BA)^nB=A(BA)^n\rangle$. For the second presentation, put $u=D^nB = AD^n$, $v=D$. Then $u^2=(D^nB)(AD^n)=D^{2n+1}=v^{2n+1}$. Conversely, from $u^2=v^{2n+1}$ define $B=v^{-n}u$, $A=u^{-1}v^{n+1}$. Then $BA=v$ and $(BA)^nB=v^nv^{-n}u=u$, whereas $A(BA)^n=u^{-1}v^{2n+1}=u$. Hence the two presentations are equivalent. 
\end{proof}

Let $q_n: G_n^{tor} \to Q_n^{tor}$ be the quotient map. The Conway meridians are identified as:
\begin{equation*}
	m_1 = x_1 \ \ \ \ \ m_2 = x_2 \ \ \ \ \ m_3 = x_3 \ \ \ \ \ m_4 = x_{2n+5}
\end{equation*}
In the free basis of Proposition \ref{prop:torusgroup}, they are $m_1 = AHD^nA^{-1}$, $m_2 = HD^n$, $m_3 = D^{-n}AD^n$, and $m_4 = DA^{-1}$. Indeed, $m_4=B_{n+1}$, while the $k=n$ case of \eqref{eq:type1ladder2} gives $B_{n+1}=A_nB_nA_n^{-1}$. Since $D=A_nB_n$, one has $B_n=A^{-1}D$, and hence $m_4=B_{n+1}=DA^{-1}$. After killing the trivial strand, $H=D^{-n}$, $B=HAH^{-1}=D^{-n}AD^n$, and $D=BA$. Consequently, $q_n(m_3)=q_n(m_4)=B$. 

\medskip

We now introduce the real representations we will compare to the branch-adapted limiting representatives. These representations send the trivial-strand meridian to the identity, so they descend to a boundary-parabolic real representation of $\pi_1(T_{2,2n+1})$. Boundary-parabolic representations of two-bridge knot groups were first studied systematically by Riley \cite{riley}. Define a recursive polynomial $C_0 = 0, C_1 = 1, C_{k+1} = (r+2)C_k-C_{k-1}$. Then let $P_n(r) = (r+1)C_n(r) - C_{n-1}(r)$. Let
\begin{equation*}
	A_r = \begin{pmatrix}1 & 0 \\ r & 1\end{pmatrix} \ \ \ \ \ B_r = \begin{pmatrix}1 & 1 \\ 0 & 1\end{pmatrix} \ \ \ \ \ D_r = B_rA_r = \begin{pmatrix}1 + r & 1 \\ r & 1\end{pmatrix}
\end{equation*}
By the Cayley-Hamilton theorem, we have $D_r^n = C_n(r)D_r - C_{n-1}(r)I$. Thus
\begin{equation*}
	D_r^n = \begin{pmatrix}P_n(r)&C_n(r) \\ rC_n(r) & C_n(r) - C_{n-1}(r)\end{pmatrix}
\end{equation*}
If $r$ is a root of $P_n$, then
\begin{equation}\label{eq:rootidentity}
	D_r^nB_r = A_rD_r^n \ \ \ \ \ D_r^{n+1} = D_r^nB_rA_r = A_rD_r^nA_r
\end{equation}
The first identity shows that $A \mapsto A_r, B \mapsto B_r$ defines a representation $\bar\rho_r: Q_n^{tor} \to SL_2(\mathbb R)$ of the quotient presentation in Proposition \ref{prop:torusgroup}. Let $\rho_r = \bar\rho_r \circ q_n$. Since $q_n(H) = D^{-n}$, one has $\rho_r(A) = A_r$, $\rho_r(H) = D_r^{-n}$, $\rho_r(B) = B_r$, $\rho_r(D) = D_r$. The quotient kills the two trivial-strand meridians, so $\rho_r(m_1) = \rho_r(m_2) = I$. For the surviving strand, $\rho_r(m_3) = D_r^{-n}A_rD_r^n=B_r$ and $\rho_r(m_4) = D_rA_r^{-1}=B_r$. Thus
\begin{equation}\label{eq:torus-meridian-images}
	\rho_r(m_1) = \rho_r(m_2) = I \ \ \ \ \ \rho_r(m_3) = \rho_r(m_4) = P_r = B_r = \begin{pmatrix}1 & 1 \\ 0 & 1\end{pmatrix}
\end{equation}
The matrix $P_r$ is a noncentral trace-$2$ parabolic.

\begin{lemma}\label{lma:roots}
The polynomial $P_n(r)$ has exactly $n$ simple roots $r_{n,j}$, all in $(-4, 0)$.
\end{lemma}

\begin{proof}
Let $x = \frac{r+2}{2}$. Then $C_k(r) = U_{k-1}(x)$, where $U_j$ is the Chebyshev polynomial of second kind, with $U_{-1} = 0$. Hence $P_n(r) = (2x-1)U_{n-1}(x) - U_{n-2}(x)$. Using the identity: $2xU_{n-1}(x) = U_n(x) + U_{n-2}(x)$, we get $P_n(r) = U_n(x) - U_{n-1}(x)$. Put $x = \cos(\theta)$. Then $U_k(\cos\theta) = \frac{\sin((k+1)\theta)}{\sin\theta}$. Thus, $P_n(r) = 0$ if and only if $\sin((n+1)\theta) - \sin(n\theta) = 0$. We also have $\sin((n+1)\theta) - \sin(n\theta) = 2\cos(\frac{(2n+1)\theta}{2})\sin(\frac{\theta}{2})$. The factor $\sin(\theta/2)$ only gives $\theta = 0$, i.e. $x = 1$, where $P_n(0) = 1$, so it's not a root. Thus, the roots come from $\cos(\frac{(2n+1)\theta}{2}) = 0$. Therefore, $\theta_{n,j} = \frac{(2j+1)\pi}{2n+1}$, $j = 0, \dots, n-1$. Thus, the roots are 
\begin{equation*}
	r_{n,j} = 2\cos(\theta_{n,j}) - 2 \ \ \ \ \ j = 0, \dots, n - 1
\end{equation*}
giving the desired $n$ real roots. Since these $n$ angles are all distinct and lie in $(0, \pi)$, the corresponding $r_{n,j}$'s are distinct. Since $P_n$ has degree $n$, these are all the roots; in particular, they are simple.
\end{proof}

Thus, when a $2n+1$ torus-trivial tangle forms the $T_1$ half of a Conway sphere decomposition $T_1 \cup_C T_2$, there are exactly $n$ candidates for an admissible real limiting datum on $T_1$. 

\medskip

For later use, put $U_r = D_r^nB_r = A_rD_r^n$. At $r = r_{n,j}$, the matrix $D_r$ has eigenvalues $e^{\pm i\theta_{n,j}}$. Consequently, 
\begin{equation}\label{eq:minusI}
	U_r^2 = (D_r^nB_r)(A_rD_r^n) = D_r^{2n+1} = -I
\end{equation}
In particular, $\mathrm{tr}(U_r) = 0$. 

\begin{lemma}\label{lma:torustrivial2}
Let $r = r_{n,j}$ be a root of $P_n$, and let $\bar\rho_r: Q_n^{tor} \to SL_2(\mathbb{R})$ be the induced quotient representation. Then hypothesis (2) of Theorem \ref{thm:main} holds for $(T_n^{tor}, \rho_r)$.
\end{lemma}

\begin{proof}
The nonabelian components of the $SL_2(\mathbb C)$-character variety of a $(2, q)$ torus knot are smooth curves parametrized by the meridian trace; see Example 2 and Remark 13 of \cite{torus2q}. Since $\bar\rho_r$ is irreducible and $\bar m = q_n(m_3)$ is a meridian, $I_{\bar m}$ is a local analytic coordinate at $\chi_{\bar\rho_r}$. This proves hypothesis (2). 
\end{proof}

\begin{lemma}\label{lma:torustrivial34}
Let $T_n^{tor}$ be a torus-trivial tangle, and let $\rho_r: G_n^{tor} \to SL_2(\mathbb{R})$ be a standard representation associated to a root $r_{n,j}$ of $P_n(r)$. Then hypotheses (3) and (4) of Theorem \ref{thm:main} hold. 
\end{lemma}

\begin{proof}
Let $a = \rho_r(A) = A_r$, and $p = \rho_r(m_3) = \rho_r(m_4) = B_r$. Since $m_2 = A^{-1}m_1A$, the traversal appearing in Theorem \ref{thm:main} is $\tau_n^{tor}=A^{-1}$. Consequently, $P^{\tau_n^{tor}}=\rho_r(\tau_n^{tor})^{-1}P\rho_r(\tau_n^{tor})=apa^{-1}$. 
For every representation of $G_n^{tor} = \langle A, H \rangle$, one has 
\begin{equation}\label{eq:torustangent}
	\tr((m_2^{-1}A^{-1}m_2-A^{-1})m_3^{-1}A) = 0
\end{equation}
Indeed, using $m_2=HD^n$, $m_3 = D^{-n}AD^n$, $D = AH^{-1}AH$, we obtain $m_2^{-1}A^{-1}m_2=D^{-(n+1)}AD^n$. Hence the two terms in \eqref{eq:torustangent} have traces $\tr(D^{-1}A)=\tr(H^{-1}A^{-1}H)=\tr(A^{-1})$ and $\tr(A^{-1}D^{-n}A^{-1}D^nA)=\tr(D^{-n}A^{-1}D^n)=\tr(A^{-1})$, respectively. Let $u$ be a tangent cocycle at $\rho_r$, and put $X=u(m_1), Y=u(m_2)$. Since $\rho_r(m_2)=I$, differentiating \eqref{eq:torustangent} gives $0=\tr((a^{-1}Y-Ya^{-1})p^{-1}a)$. The relation $m_2=A^{-1}m_1A$ gives $Y=\Ad_{a^{-1}}X$. Substitution and cyclicity of trace yield $0=\tr(X(ap^{-1}a^{-1}-p^{-1})) = -\tr(X(apa^{-1}-p))$, where the last equality uses $p^{-1}=2I-p$. Thus $\tr(u(m_1)(P^{\tau_n^{tor}}-P))=0$, which is hypothesis (3). 

\medskip

For condition (4), the fixed line of $p=B_r$ is $\langle e_1 \rangle$, while the fixed line of $P^{\tau_n^{tor}}=A_rB_rA_r^{-1}$ is $A_r\langle e_1\rangle = \langle \binom{1}{r} \rangle$. These lines are distinct because $r \in (-4, 0)$ by Lemma \ref{lma:roots}. 
\end{proof}

We finish by recording the diagrammatic pieces $l_1, l_2$ used in four-longitude decompositions. For indexing purposes, we will denote them $l_{1,n}^{tor}$ and $l_{2,n}^{tor}$. Following the arrows in Figure \ref{fig:torustrivial}, we have 
\begin{equation*}
	l_{1,n}^{tor} = A = x_{n+3} \ \ \ \ \ l_{2,n}^{tor} = \left(\prod_{k=1}^nx_{n+k+4}^{-1}\right)x_{n+3}^{2n+1}x_2\left(\prod_{k=0}^nx_{k+3}^{-1}\right)
\end{equation*}
It is clear from Figure \ref{fig:torustrivial} that $q_n(l_2)$ represents a homological longitude of the $(2,2n+1)$ torus knot. We record the relevant representation data here. 

\begin{prop}\label{prop:torus-heights}
Let $\theta_{n,j} = \frac{(2j+1)\pi}{2n+1}$, and $r=r_{n,j} = 2\cos\theta_{n,j}-2$ be a real root of $P_n$. Then we have
\begin{equation*}
	\rho_r(l_{1,n}^{tor}) = A_r = \begin{pmatrix}1&0 \\r&1\end{pmatrix} \ \ \ \ \ \rho_r(l_{2,n}^{tor}) = -B_r^{4n+2} = -\begin{pmatrix}1&4n+2\\0&1\end{pmatrix}
\end{equation*}
In particular, the central signs are $(+1, -1)$, and the transverse parameters are $(r, 4n+2)$. Let $\widetilde\rho_{n,j}: Q_n^{tor} \to \widetilde{PSL}_2(\mathbb R)$ be the meridian-normalized lift of $\bar\rho_r$, and put $\widetilde\rho_{n,j} = \widetilde{\bar\rho}_{n,j} \circ q_n$. Then
\begin{equation*}
	\mathrm{trans}(\widetilde\rho_{n,j}(l_{1,n}^{tor})) = 0 \ \ \ \ \ \mathrm{trans}(\widetilde\rho_{n,j}(l_{2,n}^{tor})) = 2j+1
\end{equation*}
\end{prop}

\begin{proof}
Put $q = 2n+1, N = 4n+2$. Recall $A_k = x_{k+3}$, $B_k = x_{n+k+4}$. Since $\rho_r(x_2) = I$, the formulas $\rho_r(A_k)=D_r^{-(n-k)}A_rD_r^{n-k}$ and $\rho_r(B_k) = D_r^{-(n-k)}A_r^{-1}D_r^{n-k+1}$ give the telescoping identities 
\begin{equation*}
	\prod_{k=1}^n\rho_r(B_k)^{-1} = D_r^{-n}A_r^n \ \ \ \ \ \prod_{k=0}^n\rho_r(A_k)^{-1} = A_r^nD_r^{-n}A_r^{-1}
\end{equation*}
Consequently, we have 
\begin{equation*}
	\rho_r(l_{2,n}^{tor}) = D_r^{-n}A_r^{4n+1}D_r^{-n}A_r^{-1} = B_r^ND_r^{-q}
\end{equation*}
Indeed, $B_r=D_r^{-n}A_rD_r^n$ and $D_r^{-(n+1)}=A_r^{-1}D_r^{-n}A_r^{-1}$ by \eqref{eq:rootidentity}. Since $D_r^q=-I$ by \eqref{eq:minusI}, it follows that 
\begin{equation*}
	\rho_r(l_{2,n}^{tor})=-B_r^N = -\begin{pmatrix}1&4n+2\\0&1\end{pmatrix}
\end{equation*}
The equality $\rho_r(l_{1,n}^{tor}) = A_r$ follows directly from the definitions.

\medskip

For the translation number statement, we work in the quotient group $Q_n^{tor}$. Put $q = 2n+1$ and $N = 2q$. Under the isomorphism $Q_n^{tor} \cong \langle u, v \mid u^2 = v^q \rangle$, with $u =D^nB$, and $v = D$, the meridian is $\mu = B = v^{-n}u$, and $f = v^q = D^q$ is central. The telescoping calculation above gives $q_n(l_{2,n}^{tor}) = B^ND^{-q}$, a homological longitude of the torus knot. Let $\widetilde{\bar\rho}_{n,j}$ be the meridian-normalized lift of $\bar\rho_r: Q_n^{tor} \to SL_2(\mathbb R)$, and let $\widetilde\rho_{n,j} = \widetilde{\bar\rho}_{n,j} \circ q_n$. For $(2, q)$ torus knots, Proposition 17.7 of \cite{dunfieldrasmussen} gives the possible absolute longitudinal heights $1, 3, \dots, q-2$. Let $h = 2j+1$. Since $D_r$ has eigenvalues $e^{\pm ih\pi/q}$ and $\mathfrak d(D_r) = 1 - r > 0$, Lemma \ref{lma:realformsign} gives $\mathrm{trans}(\widetilde D_r) \equiv -\frac{h}{q}\pmod{\mathbb Z}$. Since $\widetilde D_r^q$ is central and $\mathrm{trans}(\widetilde B_r) = 0$, $\mathrm{trans}(\widetilde{\bar\rho}_{n,j}(B^ND^{-q})) \equiv h \pmod q$. The height classification implies this integer is $\pm h'$ for some odd $h'$ with $1 \leq h' \leq q - 2$. The congruence forces it to equal $h$: the alternative $h' \equiv -h \pmod q$ would give $h + h' = q$, impossible because $q$ is odd. Thus, $\mathrm{trans}(\widetilde\rho_{n,j}(l_{2,n}^{tor})) = h = 2j+1$. Finally, $l_{1,n}^{tor} = A$ is conjugate in $Q_n^{tor}$ to the meridian $B$, and hence has translation number zero. 
\end{proof}

\subsection{Twist-trivial tangles}\label{subsec:twist-trivial}

Our second family is the \emph{twist-trivial tangles} $T_m^{tw}$, for $m \geq 1$. For $m = 1$, the associated twist knot is the trefoil; for $m \geq 2$, it is hyperbolic.

\begin{figure}[h]
	\centering
	\includegraphics[scale=.6]{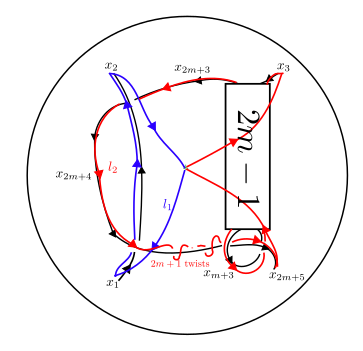}
	\caption{The twist-trivial tangle $T_m^{tw}$; the box stands for $2m-1$ right-handed crossings. The trivial strand has endpoints labeled $x_1$, $x_2$. The knotted strand, with endpoints $x_3$, $x_{2m+5}$, closes up to the $(2m+1)$-crossing twist knot. The labeled arcs are the Wirtinger generators appearing in the clasp and self-clasp relations \eqref{eq:twistinitial}, together with the four Conway meridians $m_1 = x_1$, $m_2 = x_2$, $m_3 = x_3$, $m_4 = x_{2m+5}$. The two colored curves are as follows: $l_1$ (blue) traverses the trivial strand from $x_1$ to $x_2$, and $l_2$ (red) is a parallel of the knotted strand, corrected by the $2m+1$ left-handed meridional twists.}\label{fig:twisttrivial}
\end{figure}

Using the Wirtinger presentation, we compute the fundamental group. The diagram has $2m+5$ arcs; denote the corresponding generators $x_1, \dots, x_{2m+5}$. The relators are then
\begin{equation}\label{eq:twistinitial}
	x_2 = x_{2m+4}x_1x_{2m+4}^{-1} \ \ \ \ \ x_{2m+4} = x_2x_{2m+3}x_2^{-1} \ \ \ \ \ x_{m+4} = x_{2m+5}x_{m+3}x_{2m+5}^{-1} \ \ \ \ \ x_{2m+5} = x_{m+3}x_{2m+4}x_{m+3}^{-1}
\end{equation}
\begin{equation}\label{eq:twistladder1}
	x_{k+3} = x_{2m-k+4}x_{k+2}x_{2m-k+4}^{-1} \ \ \ \ \ 1 \leq k \leq m
\end{equation}
\begin{equation}\label{eq:twistladder2}
	x_{m+k+3} = x_{m-k+4}x_{m+k+2}x_{m-k+4}^{-1} \ \ \ \ \ 2 \leq k \leq m
\end{equation}
Denote these groups $G_m^{tw}$. 

\begin{prop}\label{prop:twistgroup}
The group $G_m^{tw}$ is free on $Y = x_{2m+4}$, $H = x_2(x_{m+4}^{-1}x_{m+3})^m$. Killing the trivial-strand meridian gives the quotient
\begin{equation*}
	Q_m^{tw} = \langle a, b \mid w^ma = bw^m \rangle \ \ \ \ \ w = ba^{-1}b^{-1}a
\end{equation*}
In the conventions of \cite{hosteshanahan}, this is the fundamental group of the twist knot complement $K_{2m}$.
\end{prop}

\begin{proof}
Set $A_j = x_{m-j+3}$, $0 \leq j \leq m$, and $B_j = x_{m+j+4}$, $0 \leq j \leq m+1$. Then $A_0 = x_{m+3}$, $A_m=x_3$, $B_0 = x_{m+4}$, $B_m = x_{2m+4}$, and $B_{m+1} = x_{2m+5}$. Relations $\eqref{eq:twistladder1}$ and $\eqref{eq:twistladder2}$ become $A_j = B_jA_{j+1}B_j^{-1}$, $0 \leq j \leq m - 1$ and $B_j = A_jB_{j-1}A_j^{-1}$, $1 \leq j \leq m - 1$. Let $C_j = B_j^{-1}A_j$; for $0 \leq j \leq m - 2$ we then have $C_{j+1} = B_{j+1}^{-1}A_{j+1} = A_{j+1}B_j^{-1} = B_j^{-1}A_j = C_j$. So all $C_j$ are equal; write $C = B_0^{-1}A_0$. If $A = A_0$, $B = B_0$, and $T = A_m$, then $C = B^{-1}A$, and induction gives $A_j = C^jAC^{-j}$ and $B_j = C^jBC^{-j}$, whenever the relevant indices occur. In particular, $T = C^mAC^{-m}$ and $B_{m-1} = C^{-1}T$. Now set $X = x_2$, $Y = B_m = x_{2m+4}$, and $K = B_{m+1} = x_{2m+5}$. The remaining relations become $Y = XC^{-1}TX^{-1}$, $K = AYA^{-1}$, and $B = KAK^{-1}$. The last two identities and $C=B^{-1}A$ give $C = AYA^{-1}Y^{-1}$. Now define $H = XC^m$. Using $T=C^mAC^{-m}$, we obtain $Y = XC^{-1}TX^{-1} = HC^{-1}AH^{-1}$. Since $C^{-1}A = YAY^{-1}$, this becomes $Y = HYAY^{-1}H^{-1}$. Equivalently, $A = Y^{-1}H^{-1}YHY$. Thus $A$, and consequently $C, T, X, B, K$, and every Wirtinger generator, is a word in $H, Y$. These substitutions are reversible, and no relation remains between $H$ and $Y$. Hence $G_m^{tw} \cong \langle H, Y \rangle$. In the original generators, $H = x_2(x_{m+4}^{-1}x_{m+3})^m$, $Y = x_{2m+4}$. 

\medskip

The quotient $Q_m^{tw}$ is obtained by killing the trivial strand of Figure \ref{fig:twisttrivial}, i.e. imposing $x_1 = x_2 = 1$. The relation $B_m = x_2B_{m-1}x_2^{-1}$ then gives $B_m = B_{m-1} = C^{-1}T$. The relations in \eqref{eq:twistinitial} give $B_{m+1} = AB_mA^{-1}$, $B = B_{m+1}AB_{m+1}^{-1}$. Equivalently, $BAB_m = AB_mA$. Since $B = AC^{-1}$ and $B_m = C^{-1}T$, this reduces to $BT = TA$. Hence $B = TAT^{-1}$, and therefore $C = B^{-1}A = TA^{-1}T^{-1}A$. Together with $T = C^mAC^{-m}$, we obtain $C^mA = TC^m$. Renaming $a = A$, $b = T$, $w = ba^{-1}b^{-1}a$ gives $Q_m^{tw} \cong \langle a, b \mid w^ma = bw^m \rangle$. 

\medskip

To compare with the convention of \cite{hosteshanahan}, set $\alpha=b$, $\beta = a$, and $\omega = \beta^{-1}\alpha\beta\alpha^{-1} = w^{-1}$. Their relation $\alpha\omega^n = \omega^n\beta$ with $n=-m$ becomes $bw^m = w^ma$, which is exactly the displayed relation. Their knot is $K_{-2n} = K_{2m}$. 
\end{proof}

Let $q_m: G_m^{tw} \to Q_m^{tw}$ be the quotient map. The Conway meridians are identified as:
\begin{equation*}
	m_1 = x_1 \ \ \ \ \ m_2 = x_2 \ \ \ \ \ m_3 = x_3 \ \ \ \ \ m_4 = x_{2m+5}
\end{equation*}
In the free basis of Proposition \ref{prop:twistgroup}, they are $m_1 = Y^{-1}HC^{-m}Y$, $m_2 = HC^{-m}$, $m_3 = C^mAC^{-m}$, and $m_4 = AYA^{-1}$; here, $A = Y^{-1}H^{-1}YHY$ and $C = AYA^{-1}Y^{-1}$. After killing the trivial strand, $q_m(A) = a$, $q_m(C) = w$, $q_m(H) = w^m$, and $q_m(Y) = w^{-1}b$, so $q_m(m_3) = q_m(m_4) = b$.  

\medskip

We now introduce the real representations we will compare to the branch-adapted limiting representatives. These representations send the trivial-strand meridian to the identity, so they descend to a boundary-parabolic real representation of the associated twist-knot group. This once again uses the systematic study of two-bridge knots in \cite{riley}, as well as the construction of \cite{hosteshanahan}. For $t \in \mathbb C^*$ and $y \in \mathbb C$, define
\begin{equation*}
	A_{t,y} = \begin{pmatrix}t & 1\\ 0 & t^{-1}\end{pmatrix} \ \ \ \ \
	T_{t,y} = \begin{pmatrix}t & 0\\ y & t^{-1}\end{pmatrix}
\end{equation*}
and put $W_{t,y} = T_{t,y}A_{t,y}^{-1}T_{t,y}^{-1}A_{t,y}$. Put $s = t + t^{-1}$ and $q = \tr(W_{t,y}) = y^2 + (s^2-4)y + 2$. Define $S_0(q) = 0$, $S_1(q) = 1$, and $S_{k+1}(q) = qS_k(q) - S_{k-1}(q)$. Induction gives
\begin{equation}\label{eq:psiidentity}
	W_{t,y}^kA_{t,y} - T_{t,y}W_{t,y}^k = \Psi_k(s, y)(A_{t,y} - T_{t,y}) \ \ \ \ \ \Psi_k(s,y) = (s^2 + y - 3)S_k(q) - S_{k-1}(q)
\end{equation}
Define $\Lambda_{2m-1}(z) = \Psi_m(2, z)$. Equivalently, the odd-indexed polynomials are determined by 
\begin{equation*}
	\Lambda_{-1}(z) = 1 \ \ \ \ \ \Lambda_1(z) = z + 1 \ \ \ \ \
	\Lambda_{j+2}(z) = (z^2+2)\Lambda_j(z) - \Lambda_{j-2}(z) \text{ for odd } j \geq 1
\end{equation*}
For a real root $z$ of $\Lambda_{2m-1}$, put
\begin{equation*}
	A_z = A_{1,z} = \begin{pmatrix}1 & 1 \\ 0 & 1\end{pmatrix} \ \ \ \ \ T_z = T_{1,z} = \begin{pmatrix}1 & 0 \\ z & 1\end{pmatrix} \ \ \ \ \ W_z = W_{1,z} = \begin{pmatrix}1+z&z\\z^2 & z^2 - z + 1\end{pmatrix}
\end{equation*} 
By \eqref{eq:psiidentity}, the assignment $a \mapsto A_z$, $b \mapsto T_z$ defines a representation $\bar\rho_z: Q_m^{tw} \to SL_2(\mathbb R)$. Let $\rho_z = \bar\rho_z \circ q_m$. Then
\begin{equation*}
	\rho_z(A) = A_z \ \ \ \ \ \rho_z(C) = W_z \ \ \ \ \ \rho_z(H) = W_z^m \ \ \ \ \ \rho_z(Y) = W_z^{-1}T_z
\end{equation*}
and consequently 
\begin{equation*}
	\rho_z(m_1)  = \rho_z(m_2) = I \ \ \ \ \ \rho_z(m_3) = \rho_z(m_4) = P_z = T_z = \begin{pmatrix}1 & 0 \\ z & 1\end{pmatrix}
\end{equation*}

\begin{lemma}\label{lma:twistroots}
For every $m \geq 1$, the degree-$(2m-1)$ polynomial $\Lambda_{2m-1}$ is irreducible over $\mathbb Q$ and has exactly one real root $z_m$. One has $z_1 = -1$, while $z_m \in (-1, 0)$ for $m \geq 2$. Moreover, for $m \ge 2$, the associated twist knot is hyperbolic, and there is a root $z_{geo}$ for which $\bar\rho_{z_{geo}}$ is conjugate to a lift of its discrete faithful representation. Its trace field is $\mathbb Q(z_{geo})$; its unique real place sends $z_{geo}$ to $z_m$, and the remaining $2m-2$ roots of $\Lambda_{2m-1}$ are the other Galois conjugates of $z_{geo}$. 
\end{lemma}

\begin{proof}
Let $\Phi_n$ be the parabolic polynomial from \cite{hosteshanahan}. Comparing the presentation of Proposition \ref{prop:twistgroup} with their presentation, and exchanging the two meridians, gives $\Lambda_{2m-1}(z) = \Phi_{-m}(-z)$. The same identity then follows directly from the above recursion, and the cases $m = 1, 2$ agree. For $m \geq 2$, Theorem 1 of \cite{hosteshanahan} gives irreducibility, the geometric root, and the trace-field statement; for $m = 1$, one has $\Lambda_1(z) = z + 1$, $\Phi_{-1}(z) = 1 - z$. For $m \geq 2$, Proposition 1 of \cite{hosteshanahan}, together with the substitution $z = x - x^{-1}$, shows that $\Lambda_{2m-1}$ has exactly one real root; for $m = 1$ this is immediate. Finally, $\Lambda_{2m-1}(0) = 1$ for every $m$, while $\Lambda_{2m-1}(-1) < 0$ for $m \geq 2$. The latter follows inductively from $\Lambda_1(-1) = 0$, $\Lambda_3(-1) = -1$, and $\Lambda_{j+2}(-1) = 3\Lambda_j(-1) - \Lambda_{j-2}(-1)$. Thus $z_m \in (-1, 0)$ for $m \geq 2$. Since $\Lambda_{2m-1}$ is irreducible, its roots are the images of $z_{geo}$ under the embeddings of the trace field; uniqueness of the real root gives the final assertions. Hyperbolicity follows because the knot is prime alternating, and for $m \geq 2$, is not a $(2, q)$-torus knot \cite{menasco}. 
\end{proof}

Thus, a twist-trivial tangle contributes exactly one candidate for an admissible real limiting datum, namely $\rho_{z_m}$. 

\begin{lemma}\label{lma:twist2}
Let $m \geq 1$, let $z = z_m$, and let $\bar\rho_z: Q_m^{tw} \to SL_2(\mathbb R)$ be the quotient representation. Then hypothesis (2) of Theorem \ref{thm:main} holds for $(T_m^{tw}, \rho_z)$.
\end{lemma}

\begin{proof}
Since $z < 0$, the noncentral parabolics $A_z = \bar\rho_z(a)$ and $T_z = \bar\rho_z(b)$ have distinct fixed lines, so $\bar\rho_z$ is irreducible. By the form of the representation established above, the nonabelian character germ at $\chi_{\bar\rho_z}$ is the plane curve germ $\Psi_m(s, y) = 0$, where $s = \tr(a) = \tr(b)$; see Proposition 8 of \cite{twist}. At the boundary-parabolic point $(s, y) = (2, z)$, one has $\partial_y\Psi_m(2,z) = \Lambda_{2m-1}'(z) \neq 0$, because $z$ is a simple root of $\Lambda_{2m-1}$. The implicit function theorem therefore expresses $y$ as a holomorphic function of $s$. Hence the character germ is smooth and $s$ is a local analytic coordinate. Since the surviving meridian is $\bar m = q_m(m_3) = b$, one has $I_{\bar m} = I_b = s$. This proves hypothesis (2).
\end{proof}

\begin{lemma}\label{lma:twist34}
Let $m \geq 1$, let $\rho_z: G_m^{tw} \to SL_2(\mathbb{R})$ be the standard representation associated to the real root $z = z_m$. Then hypotheses (3) and (4) of Theorem \ref{thm:main} hold. 
\end{lemma}

\begin{proof}
Since $m_2 = Ym_1Y^{-1}$, the traversal appearing in Theorem \ref{thm:main} is $\tau_m^{tw} = Y$. Put 
\begin{equation*}
	R = \rho_z(Y) = W_z^{-1}T_z = \begin{pmatrix}1-z&-z\\z&1+z\end{pmatrix} \ \ \ \ \ P= \rho_z(m_3) = T_z
\end{equation*}
For every representation of $G_m^{tw} = \langle H, Y \rangle$, one has
\begin{equation}\label{eq:twisttangent}
	\tr((m_2^{-1}Ym_2-Y)m_3^{-1}Y^{-1}) = 0
\end{equation}
Indeed, $m_2^{-1}Ym_2 = C^{-1}T$, so the first term in \eqref{eq:twisttangent} has trace $\tr(C^{-1}Y^{-1})$. The identities $A = (Y^{-1}H^{-1})Y(Y^{-1}H^{-1})^{-1}$ and $C^{-1}Y^{-1} = (H^{-1}YH)Y^{-1}(H^{-1}YH)^{-1}$ show that $\tr(C^{-1}Y^{-1}) = \tr(Y) = \tr(A) = \tr(T)$. The second term has trace $\tr(YT^{-1}Y^{-1}) = \tr(T)$, proving \eqref{eq:twisttangent}. Let $u$ be a tangent cocycle at $\rho_z$, and put $X_i = u(m_i)$. Since $\rho_z(m_1) = \rho_z(m_2) = I$, differentiating \eqref{eq:twisttangent} gives $0 = \tr((RX_2-X_2R)P^{-1}R^{-1})$. The relation $m_2 = Ym_1Y^{-1}$ gives $X_2 = \Ad_RX_1$. Substitution and cyclicity of trace yield $0 = \tr(X_1(R^{-1}P^{-1}R - P^{-1})) = -\tr(X_1(R^{-1}PR - P))$, where the last equality uses $P^{-1} = 2I-P$. Since $P^{\tau_m^{tw}} = R^{-1}PR$, this is hypothesis (3). 

\medskip

For hypothesis (4), the fixed line of $P = T_z$ is $\langle e_2 \rangle$, while the fixed line of $P^{\tau_m^{tw}}$ is $R^{-1}\langle e_2 \rangle$. Since $R^{-1}e_2 = \binom{z}{1-z}$ and $z \neq 0$, these lines are distinct. 
\end{proof}

We finish by recording the diagrammatic pieces $l_1, l_2$ used in four-longitude decompositions. For indexing purposes, we will denote them $l_{1,m}^{tw}$ and $l_{2,m}^{tw}$. Following the arrows in Figure \ref{fig:twisttrivial}, we have 
\begin{equation*}
	l_1 = Y^{-1} = x_{2m+4}^{-1} \ \ \ \ \ l_2 = \left(\prod_{k=0}^{m-1}x_{2m-k+3}^{-1}\right)x_{2m+5}^{-1}\left(\prod_{k=1}^{m-1}x_{m-k+3}^{-1}\right)x_2^{-1}x_{2m+4}^{2m+1}x_{m+3}^{-1}
\end{equation*}
It is clear from Figure \ref{fig:twisttrivial} that $q_m(l_2)$ represents a homological longitude of the twist knot with $2m+1$ crossings. We record the relevant representation data here.

\begin{prop}\label{prop:twist-height}
Let $S = \left(\begin{smallmatrix}0&-1\\1&1\end{smallmatrix}\right)$. We have
\begin{equation*}
	S^{-1}\rho_z(l_{1,m}^{tw})S = \begin{pmatrix}1 & 0 \\ -z & 1\end{pmatrix} \ \ \ \ \ S^{-1}\rho_z(l_{2,m}^{tw})S = -\begin{pmatrix}1&4-2z\\0&1\end{pmatrix} \ \ \ \ \ S^{-1}P_zS = S^{-1}T_zS = \begin{pmatrix}1 & -z \\ 0 & 1\end{pmatrix}
\end{equation*}
In particular, the central signs are $(+1, -1)$, and the transverse parameters are $(-z, 4-2z)$. Let $\widetilde\rho_z$ be a meridian-normalized lift to $\widetilde{PSL}_2(\mathbb R)$. Then
\begin{equation*}
	\mathrm{trans}(\widetilde\rho_z(l_{1,m}^{tw})) = 0 \ \ \ \ \ \mathrm{trans}(\widetilde\rho_z(l_{2,m}^{tw})) = 1
\end{equation*}
\end{prop}

\begin{proof}
Suppress the subscript $z$, and write $A = A_z$, $T = T_z$, $W = W_z$, and $B = AW^{-1} = TAT^{-1}$. Recall that $W^mA = TW^m$. For $A_j = x_{m-j+3}$, $B_j = x_{m+j+4}$, we have $\rho_z(A_j) = W^jAW^{-j}$ and $\rho_z(B_j) = W^jBW^{-j}$. Moreover, $\rho_z(x_2) = I$, $\rho_z(x_{2m+4}) = W^{-1}T$, and $\rho_z(x_{2m+5})=T$. Since $l_{1,m}^{tw}=x_{2m+4}^{-1}$, 
\begin{equation*}
	\rho_z(l_{1,m}^{tw}) = (W^{-1}T)^{-1}=T^{-1}W=A^{-1}T^{-1}A=\begin{pmatrix}1+z&z\\-z&1-z\end{pmatrix}
\end{equation*}
For the second traversal, put $E_m=B_0B_1\dots B_{m-1}$, and $F_m=A_{m-1}A_{m-2}\dots A_1$. We then have $l_{2,m}^{tw}=E_m^{-1}x_{2m+5}^{-1}F_m^{-1}x_2^{-1}x_{2m+4}^{2m+1}A^{-1}$. Substitution of the preceding Wirtinger formulas gives 
\begin{align*}
	 \rho_z(l_{2,m}^{tw}) = &(W^{m-1}B^{-1}W^{-(m-1)})\dots(WB^{-1}W^{-1})B^{-1}T^{-1}\dots(WA^{-1}W^{-1})\\&\dots(W^{m-1}A^{-1}W^{-(m-1)})(W^{-1}T)^{2m+1}A^{-1}
\end{align*}
Successively cancelling the adjacent powers of $W$, and using $B = AW^{-1}$, $W^{-1}T=A^{-1}TA$, and $W^mA=TW^m$ yields $\rho_z(l_{2,m}^{tw})=W^m(AT^{-1}A^{-1}T)^m$. Put $R=W^m$, and $\overleftarrow W = AT^{-1}A^{-1}T$. Write $R = \left(\begin{smallmatrix}a&b\\c&d\end{smallmatrix}\right)$. The identity $RA=TR$ gives $a=0$, $c=zb$. Since $R$ is a power of $W$, it commutes with $W$; comparing the $(1, 2)$ entries of $RW$ and $WR$ gives $d=(z-2)b$. Finally, $\det R = 1$ gives $-zb^2=1$. Thus, $R=\left(\begin{smallmatrix}0&b\\zb&(z-2)b\end{smallmatrix}\right)$. Let $J=\left(\begin{smallmatrix}0&1\\1&0\end{smallmatrix}\right)$. The map $M \mapsto JM^{\mathsf T}J$ fixes both $A$ and $T$, and therefore sends $W$ to $\overleftarrow W$. Consequently, $\overleftarrow W^m=JR^{\mathsf T}J=\left(\begin{smallmatrix}(z-2)b&b\\zb&0\end{smallmatrix}\right)$. The previous expression then becomes
\begin{equation*}
	\rho_z(l_{2,m}^{tw}) = RJR^{\mathsf T}J = \begin{pmatrix}zb^2&0\\2z(z-2)b^2&zb^2\end{pmatrix} = -\begin{pmatrix}1&0\\2(z-2)&1\end{pmatrix}
\end{equation*}
Conjugation by $S=\left(\begin{smallmatrix}0&-1\\1&1\end{smallmatrix}\right)$ puts $\rho_z(l_{i,m}^{tw})$ in the desired form. 

\medskip

For the translation number statement, work in the quotient twist knot group. Let $\widetilde{\bar\rho}_z: Q_m^{tw} \to \widetilde{PSL}_2(\mathbb R)$ be the meridian-normalized lift of $\bar\rho_z$, and put $\widetilde\rho_z = \widetilde{\bar\rho}_z \circ q_m$. Note that $q_m(l_{2,m}^{tw})$ is a homological longitude of the twist knot. The associated odd twist knot has Seifert genus one. Hence the Milnor-Wood bound featured as Proposition 6.5 of \cite{cullerdunfield} gives $|\mathrm{trans}(\widetilde{\bar\rho}_z(q_m(l_{2,m}^{tw})))| \leq 1$. But this is a noncentral parabolic of trace $-2$, so Lemma \ref{lma:parity} implies that its translation number is odd. Thus the translation number must be $\pm 1$. Both $W = TA^{-1}T^{-1}A$ and $\overleftarrow W=AT^{-1}A^{-1}T$ are commutators, and their matrices have trace $z^2+2 > 2$. Lemma \ref{lma:commutator} bounds the translation numbers of their word lifts by one, while Lemma \ref{lma:parity} makes those translation numbers even. Hence both are zero. By homogeneity, the word lifts of $R=\bar\rho_z(W^m)$, $\overleftarrow R = \bar\rho_z(\overleftarrow W^m)$ are precisely their translation-zero lifts.

\medskip

Retain the notation $R = \bar\rho_z(W^m)$ and $\overleftarrow R = \bar\rho_z(\overleftarrow W^m)$ from above. In the angular convention of Section \ref{sec:h00}, the signs of their entries show that their translation-zero lifts satisfy $\widetilde R(0) = \frac12$, and $\widetilde{\overleftarrow R}(\frac12) = 1$. Since $q_m(l_{2,m}^{tw}) = R\overleftarrow R$, deck equivariance gives $\widetilde\rho_z(l_{2,m}^{tw})(\frac12) = \widetilde R(1) = \frac32$. Its projective image is parabolic and fixes the line represented by $1/2$, and hence $\mathrm{trans}(\widetilde\rho_z(l_{2,m}^{tw})) = 1$. Finally, $\rho_z(l_{1,m}^{tw}) = A_z^{-1}T_z^{-1}A_z$ is conjugate to the inverse of the meridian $T_z$, so $\mathrm{trans}(\widetilde\rho_z(l_{1,m}^{tw})) = 0$. 
\end{proof}

\subsection{Mirrors and flips}\label{subsec:mirrors}

The families of knots in Section \ref{sec:families} each combine two of the four tangle variants constructed below; this subsection transfers all of the hypotheses of Theorem \ref{thm:main} as well as the peripheral data. We here describe \emph{mirrors} and \emph{flips} of torus- and twist-trivial tangles. The mirror reverses every crossing of the diagram; at the level of the Wirtinger presentation, it inverts each meridian generator (Lemma \ref{lma:mirror-transfer}). The flip reverses the two crossings between the trivial strand and the knotted strand. In this subsection, we describe how mirrors and flips transform the certification data and translation numbers of the preceding two subsections.

\begin{figure}[h]
	\centering
	\begin{subfigure}[t]{0.48\textwidth}
		\centering
		\includegraphics[scale=.55]{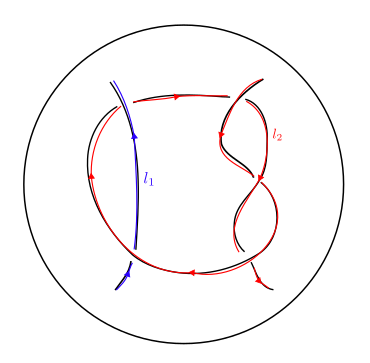}
		\caption{The trefoil-trivial tangle $T_1^{tor}$.}\label{fig:trefoiltrivialdefault}
	\end{subfigure}\hfill
	\begin{subfigure}[t]{0.48\textwidth}
		\centering
		\includegraphics[scale=.55]{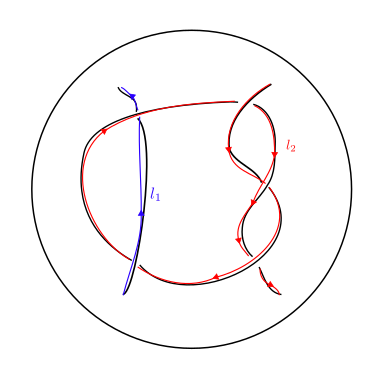}
		\caption{The flipped trefoil-trivial tangle $(T_1^{tor})^{fl}$.}\label{fig:trefoiltrivialflipped}
	\end{subfigure}
	\medskip
	\begin{subfigure}[t]{0.48\textwidth}
		\centering
		\includegraphics[scale=.55]{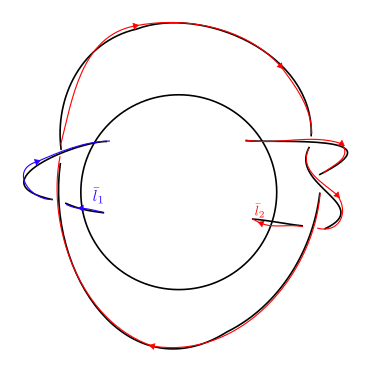}
		\caption{The mirrored trefoil-trivial tangle $\overline{T_1^{tor}}$.}\label{fig:trefoiltrivialmirror}
	\end{subfigure}\hfill
	\begin{subfigure}[t]{0.48\textwidth}
		\centering
		\includegraphics[scale=.55]{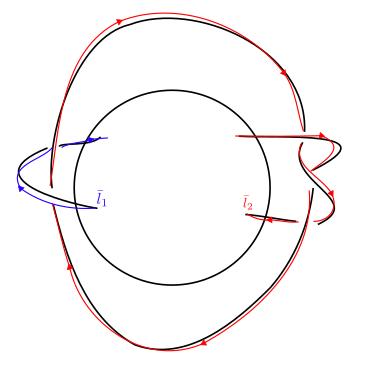}
		\caption{The mirror-flipped trefoil-trivial tangle $\overline{(T_1^{tor})^{fl}}$.}\label{fig:trefoiltrivialmirrorflipped}
	\end{subfigure}
	\caption{The trefoil-trivial tangle, its flip, its mirror, and its mirror-flip. The panels show bare diagrams; labels, orientations, and traversals for the standard model are fixed in Figure \ref{fig:torustrivial}. Panels (c)–(d) are drawn in the outside presentation, as the mirrored tangles appear in the family gluings.}
	\label{fig:mirrorflips}
\end{figure}

\begin{definition}\label{def:decorations}
Let $T$ be a torus- or twist-trivial tangle, presented by the standard diagram of Figure \ref{fig:torustrivial} or Figure \ref{fig:twisttrivial}. The \emph{mirror} $\overline{T}$ is obtained by reversing every crossing, and the \emph{flip} $T^{fl}$ is the tangle presented by the diagram obtained by reversing the two clasp crossings between the trivial and knotted strands. The two operations commute, and their composite $\overline{T}^{fl} = \overline{T^{fl}}$ reverses exactly the box crossings. The four \emph{tangle variants} $T$, $T^{fl}$, $\overline{T}$, $\overline{T}^{fl}$ for the trefoil-trivial tangle are shown in Figure \ref{fig:mirrorflips}. 
\end{definition}

Mirrored tangles occupy the outside slot of the family gluings. The following lemma transfers conditions (1)-(4) of Theorem \ref{thm:main} and the representation data.

\begin{lemma}\label{lma:mirror-transfer}
Let $T$ have a Wirtinger presentation on $x_1, \dots, x_N$, and let $\overline{T}$ be its mirror with Wirtinger generators $\bar x_1, \dots, \bar x_N$. Then $\iota: \pi_1(\overline{T}) \to \pi_1(T)$, $\iota(\bar x_j) = x_j^{-1}$, is an isomorphism. Let $\Omega: SL_2(\mathbb R) \to SL_2(\mathbb R)$ negate all off-diagonal entries, so $\Omega = \Ad_H$ where $H = \left(\begin{smallmatrix}1&0\\0&-1\end{smallmatrix}\right)$. Put $\overline\rho = \Omega\circ\rho\circ\iota$. Then:
\begin{enumerate}
	\item $(\overline T, \overline\rho)$ satisfies hypotheses (1)-(4) of Theorem \ref{thm:main} if and only if $(T, \rho)$ does. 
	\item The nontrivial Conway parabolic is $\overline P = \Omega(P^{-1})$. 
	\item The traversals $\bar l_i = \iota^{-1}(l_i)$ are the ones in Figure \ref{fig:mirrorflips}. If
	\begin{equation*}
		\rho(l_1) = \epsilon_1\begin{pmatrix}1 & 0 \\ \alpha & 1\end{pmatrix} \ \ \ \ \ \rho(l_2) = \epsilon_2\begin{pmatrix}1 & \beta \\ 0 & 1\end{pmatrix}
	\end{equation*}
	then we have 
	\begin{equation*}
		\overline\rho(\bar l_1) = \epsilon_1\begin{pmatrix}1 & 0 \\ -\alpha & 1\end{pmatrix} \ \ \ \ \ \overline\rho(\bar l_2) = \epsilon_2\begin{pmatrix}1 & -\beta \\ 0 & 1\end{pmatrix}
	\end{equation*}
	The central signs are unchanged, and the transverse parameters negate.
	\item If $\widetilde\rho$ is a lift of $\rho$ to $\widetilde{PSL}_2(\mathbb R)$, then $\widetilde{\overline\rho} := \widetilde\Omega\circ\widetilde\rho\circ\iota$ is a lift of $\overline\rho$, and $\operatorname{trans}\bigl(\widetilde{\overline\rho}(\bar g)\bigr) = -\operatorname{trans}\bigl(\widetilde\rho(g)\bigr)$ whenever $\iota(\bar g) = g$.
\end{enumerate}
\end{lemma}

\begin{proof}
Mirroring replaces a right-handed Wirtinger relation $x'=yxy^{-1}$ by $\bar x'=\bar y^{-1}\bar x\bar y$. Under $\iota$ the latter becomes the inverse of the former, so $\iota$ is an isomorphism. Notice that meridians invert, i.e. $\overline\rho(\bar m_j) = \Omega(\rho(m_j)^{-1})$, whereas a transported word $\bar g = \iota^{-1}(g)$ satisfies $\overline\rho(\bar g) = \Omega(\rho(g))$ with no additional inversion. Precomposition by $\iota$ and postcomposition by $\Omega$ identify the relevant representations, character germs and quotients obtained by killing $m_1$ and $\bar m_1$. This transfers hypotheses (1) and (2). Let $\tau = l_1^{-1}$ and $\bar\tau = \bar l_1^{-1}$. The corresponding tangent cocycles satisfy $\bar u(\bar m_1) = -\Ad_Hu(m_1)$, $\overline P^{\bar\tau} - \overline P = -\Ad_H(P^\tau - P)$, because $Q^{-1} = 2I - Q$ for every trace-2 parabolic $Q$. The two signs cancel in the trace pairing, proving hypothesis (3). Inversion does not change a parabolic fixed line, so carries the two lines in hypothesis (4) injectively. Applying $\iota$ to the mirrored endpoint relations gives the inverses of the original relations, which proves the assertion about $\bar l_1, \bar l_2$. Conclusions (2) and (3) follow because $\Omega$ flips off-diagonal signs. The translation-number formula follows from Lemma \ref{lma:outerreverses}. 
\end{proof}

We next treat flips. The operation is depicted in the below figure. 

\begin{figure}[h]
	\centering
	\includegraphics[width=\textwidth]{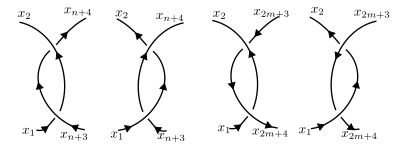}
	\caption{From left to right: the clasp of the torus-trivial, flipped torus-trivial, twist-trivial, and flipped twist-trivial tangles. Only the crossings of the clasp change.}\label{fig:clasp}
\end{figure}

Let $T_n^{tor}$ be a torus-trivial tangle. Under the flip relation, the Wirtinger presentation only changes the clasp relations from \eqref{eq:type1initial} to 
\begin{equation}\label{eq:type1initialflip}
	x_2 = x_{n+4}x_1x_{n+4}^{-1} \ \ \ \ \ x_{n+4} = x_1x_{n+3}x_1^{-1}
\end{equation}
Let $T_m^{tw}$ be a twist-trivial tangle. Under the flip relation, the Wirtinger presentation changes the first two relations from \eqref{eq:twistinitial} to 
\begin{equation}\label{eq:twistinitialflip}
	x_2 = x_{2m+3}^{-1}x_1x_{2m+3} \ \ \ \ \ x_{2m+4} = x_1^{-1}x_{2m+3}x_1
\end{equation}
Nothing else about the presentations changes. For the torus case, retain the notation $A_k = x_{k+3}$, $B_k = x_{n+k+4}$, $D = A_kB_k$, $A = A_n = x_{n+3}$, $B_0 = x_{n+4}$. For the twist case, retain the notation $A = x_{m+3}$, $C = x_{m+4}^{-1}x_{m+3}$, $T = x_3$, and put $Z = x_{2m+3} = C^{-1}T$, $Y = x_{2m+4}$. 

\begin{lemma}\label{lma:flipgroups}
Let $T$ denote $T_n^{tor}$ or $T_m^{tw}$, and let $T^{fl}$ be its flip. 
\begin{enumerate}
	\item The flipped torus-trivial group is free of rank two on $A = x_{n+3}, H^{fl} = x_1^{-1}D^{-n}$, with $D = A(H^{fl})^{-1}AH^{fl}$. The flipped twist-trivial group is free of rank two on $Y = x_{2m+4}$, $H^{fl} = x_1^{-1}C^m$, with $A = Y^{-1}(H^{fl})^{-1}YH^{fl}Y$ and $C = AYA^{-1}Y^{-1}$.
	\item Killing the trivial-strand meridian gives the same quotients $Q_n^{tor}$ and $Q_m^{tw}$ for standard and flipped tangles, with the same presentations as Propositions \ref{prop:torusgroup} and \ref{prop:twistgroup}. 
	\item Let $r$ be a root of $P_n$. The assignments $\rho_r^{fl}(A) = A_r$, $\rho_r^{fl}(H^{fl}) = D_r^{-n}$ define a representation of the flipped torus-trivial group. It satisfies 
	\begin{equation*}
		\rho_r^{fl}(m_1) = \rho_r^{fl}(m_2) = I \ \ \ \ \ \rho_r^{fl}(m_3) = \rho_r^{fl}(m_4) = B_r
	\end{equation*} Let $z$ be a real root of $\Lambda_{2m-1}$. The assignments $\rho_z^{fl}(Y) = W_z^{-1}T_z$, $\rho_z^{fl}(H^{fl}) = W_z^m$ define a representation of the flipped twist-trivial group. It satisfies 
	\begin{equation*}
		\rho_z^{fl}(m_1) = \rho_z^{fl}(m_2) = I \ \ \ \ \ \rho_z^{fl}(m_3) = \rho_z^{fl}(m_4) = T_z
	\end{equation*}
	In both cases, every knot-strand Wirtinger generator has the same image as $\rho_r$ and $\rho_z$; indeed, both representations descend to the same $\bar\rho_r$ and $\bar\rho_z$ on the quotients. 
	\item Let $l_{1,n}^{tor, fl} = x_{n+4}^{-1}$, and $l_{1,m}^{tw, fl} = x_{2m+3}$ be the identity-strand pieces in the flipped diagrams. Then 
	\begin{equation*}
		\rho_r^{fl}(l_{1,n}^{tor,fl}) = \rho_r(l_{1,n}^{tor})^{-1} \ \ \ \ \ \rho_z^{fl}(l_{1,m}^{tw, fl}) = \rho_z(l_{1,m}^{tw})^{-1}
	\end{equation*}
	If $l_2^{fl}$ denotes the knot-strand piece in the flipped diagram, its quotient image agrees with that of $l_2$. Consequently, $\rho^{fl}(l_2^{fl}) = \rho(l_2)$. 
	\item The pair $(T^{fl}, \rho^{fl})$ satisfies hypotheses (1)-(4) of Theorem \ref{thm:main} if and only if $(T, \rho)$ does. 
	\item Let $\widetilde{\bar\rho}: Q \to \widetilde{PSL}_2(\mathbb R)$ be the meridian-normalized lift of the quotient, and define the standard and flipped tangle lifts by composition with their quotient maps. Then 
	\begin{equation*}
		\mathrm{trans}(\widetilde\rho^{fl}(l_1^{fl})) = 0 \ \ \ \ \ \mathrm{trans}(\widetilde\rho^{fl}(l_2^{fl})) = \mathrm{trans}(\widetilde\rho(l_2))
	\end{equation*} 
\end{enumerate}
\end{lemma}

\begin{proof}
We first prove the free-group statements. In the torus case, the unchanged relations \eqref{eq:type1ladder1}-\eqref{eq:type1ladder2} give $B_0 = D^{-n}A^{-1}D^{n+1}$. By \eqref{eq:type1initialflip}, $B_0= x_1Ax_1^{-1}$. Since $H^{fl}=x_1^{-1}D^{-n}$ and $x_1 = D^{-n}(H^{fl})^{-1}$, substitution gives $D^{-n}A^{-1}D^{n+1} = D^{-n}(H^{fl})^{-1}AH^{fl}D^n$. Cancelling the outer powers of $D$ yields $A^{-1}D = (H^{fl})^{-1}AH^{fl}$, or equivalently $D = A(H^{fl})^{-1}AH^{fl}$. The remaining clasp relation $x_2 = B_0x_1B_0^{-1}$ eliminates $x_2$. Thus, every Wirtinger generator is a word in $A, H^{fl}$, the substitutions are reversible, and no relation remains. In the twist case, the unchanged relations \eqref{eq:twistladder1}-\eqref{eq:twistladder2} give $T = C^mAC^{-m}$, $Z = C^{-1}T$, while the unchanged latter relations of \eqref{eq:twistinitial} give $C = AYA^{-1}Y^{-1}$. Put $H^{fl} = x_1^{-1}C^m$, $x_1^{-1}=H^{fl}C^{-m}$. The second relation of \eqref{eq:twistinitialflip} gives 
\begin{equation*}
	Y = x_1^{-1}Zx_1 = H^{fl}C^{-m}C^{-1}(C^mAC^{-m})C^m(H^{fl})^{-1} = H^{fl}C^{-1}A(H^{fl})^{-1}
\end{equation*}
Since $C^{-1}A = YAY^{-1}$, this is equivalent to $A = Y^{-1}(H^{fl})^{-1}YH^{fl}Y$. The first relation of \eqref{eq:twistinitialflip} eliminates $x_2$. Thus the flipped twist group is free on $Y, H^{fl}$. 

\medskip

Killing $m_1$ also kills $m_2$, since they are conjugate. In the torus case this set $H^{fl} = D^{-n}$, and in the twist case it sets $H^{fl} = C^m$. The standard and flipped quotient presentations are therefore identical, with the surviving-strand generators identified literally. This proves (2). 

\medskip

For a torus-trivial tangle, the terminal relation gives $\rho_r^{fl}(D) = A_rD_r^nA_rD_r^{-n} = D_r$ by the root identity \eqref{eq:rootidentity}. Hence $\rho_r^{fl}(x_1^{-1}) = \rho_r^{fl}(H^{fl}D^n) = D_r^{-n}D_r^n = I$. It follows from the two clasp relations that both trivial-strand meridians map to $I$. The unchanged relators then give the standard image of every knot-strand generator. In particular, $\rho_r^{fl}(m_3) = \rho_r^{fl}(m_4) = B_r$. For a twist-trivial tangle, the free-group formulas give $\rho_z^{fl}(A) = A_z$, $\rho_z^{fl}(C) = W_z$, and $\rho_z^{fl}(T) = T_z$. Moreover, $\rho_z^{fl}(x_1^{-1}) = \rho_z^{fl}(H^{fl}C^{-m}) = W_z^mW_z^{-m} = I$. Again, both trivial-strand meridians map to $I$, while all knot-strand generators map to their usual images. 

\medskip

In the torus case, $\rho_r^{fl}(x_{n+4}) = A_r$, so $\rho_r^{fl}(l_{1,n}^{tor,fl}) = A_r^{-1} = \rho_r(l_{1,n}^{tor})^{-1}$. In the twist case, the second flipped clasp relation and $\rho_z^{fl}(x_1) = I$ give $\rho_z^{fl}(x_{2m+3}) = \rho_z^{fl}(x_{2m+4}) = W_z^{-1}T_z$. This is the inverse of $\rho_z(l_{1,m}^{tw})$. The knot-strand traversal is the same in the quotient, so the corresponding matrices for $l_2$ agree.

\medskip

We now transfer hypotheses (1)-(4) of Theorem \ref{thm:main}. Put $w = HD^n$ in the torus case and $w = HC^{-m}$ in the twist case. In the standard tangle, one has $w = m_2$, whereas in the flipped tangle one has $w = m_1^{-1}$. Let $\tau_+$ and $\tau_-$ denote the conjugating words in $m_2 = \tau_\pm m_1\tau_\pm^{-1}$ for the standard and flipped tangles. Explicitly, $(\tau_+, \tau_-) = (A^{-1}, x_{n+4})$ in the torus case, and $(\tau_+, \tau_-) = (Y, x_{2m+3}^{-1})$ in the twist case. We have $\rho(\tau_-) = \rho(\tau_+)^{-1}$. Hypotheses (1) and (2) are identical under the common free group and quotient identifications. For hypothesis (3), put $R = \rho(\tau_+)$, and $V = u(w)$. Since $\rho(w) = I$, the cocycle rule gives $u(m_1) = \mathrm{Ad}_{R^{-1}}V$ in the standard tangle, and $u(m_1) = -V$ in the flipped tangle. Moreover, $P^{\tau_+} = R^{-1}PR$, while $P^{\tau_-} = RPR^{-1}$. Cyclicity of trace gives $\tr(u(m_1)(P^{\tau_\pm}-P)) = -\tr(V(RPR^{-1}-P))$. Thus hypothesis (3) is equivalent for standard and flipped tangles. For hypothesis (4), the fixed lines of $R^{-1}PR$ and $RPR^{-1}$ are distinct from $\mathrm{Fix}(P)$ simultaneously, since $R\mathrm{Fix}(P) = \mathrm{Fix}(P) \iff R^{-1}\mathrm{Fix}(P) = \mathrm{Fix}(P)$. This proves conclusion (5); the converse follows because flipping is an involution. 

\medskip

Finally, the quotient equality for the knot-side pieces implies the equality of their translation numbers under the quotient-induced meridian-normalized lifts. The identity-strand pieces are conjugate to a meridian or its inverse in the quotient, so their translation numbers are zero.
\end{proof}

We record the relevant data for the criteria here. The general slogan is that mirrors negate translation numbers of $l_2$ and both parabolic parameters, while flips negate the parabolic parameter of $l_1$. 

\medskip

Write the standard data as
\begin{equation*}
	\rho(m_3) = \begin{pmatrix}1&c\\0&1\end{pmatrix} \ \ \ \ \ \rho(l_1) = \begin{pmatrix}1&0 \\ \alpha&1\end{pmatrix} \ \ \ \ \ \rho(l_2) = -\begin{pmatrix}1&\beta \\ 0 & 1\end{pmatrix}
\end{equation*}
Let $(\mathrm{trans}(\widetilde\rho(l_1)), \mathrm{trans}(\widetilde\rho(l_2))) = (0, h_K)$. For the torus models, $(c, \alpha, \beta, h_K) = (1, r, 4n+2, 2j+1)$ at the root $r = r_{n,j}$. For the twist models, $(c, \alpha, \beta, h_K) = (-z, -z, 4-2z, 1)$. 

\begin{table}[H]
	\centering
	\begin{tabular}{lcccc}
		\toprule
		variant
		& pieces
		& $(c;\alpha,\beta)$
		& central signs
		& model heights\\
		\midrule
		$T$
		& $(l_1,l_2)$
		& $(c;\alpha,\beta)$
		& $(+1,-1)$
		& $(0,h_K)$\\
		$T^{\mathrm{fl}}$
		& $(l_1^{\mathrm{fl}},l_2^{\mathrm{fl}})$
		& $(c;-\alpha,\beta)$
		& $(+1,-1)$
		& $(0,h_K)$\\
		$\overline T$
		& $(\bar l_1,\bar l_2)$
		& $(c;-\alpha,-\beta)$
		& $(+1,-1)$
		& $(0,-h_K)$\\
		$\overline T^{\mathrm{fl}}$
		& $(\bar l_1^{\mathrm{fl}},
		\bar l_2^{\mathrm{fl}})$
		& $(c;\alpha,-\beta)$
		& $(+1,-1)$
		& $(0,-h_K)$\\
		\bottomrule
	\end{tabular}
	\caption{Peripheral model data for the four tangle variants. A flip negates only the identity-strand parameter. Mirroring negates both piece parameters and translation numbers, while the Conway peripheral parameter $c$ and the central signs remain unchanged. The sign of $\beta$ is recorded for completeness; the positivity criterion of Theorem \ref{thm:modelcriterion} uses only $c$ and $\alpha$.}
	\label{tab:decorations}
\end{table}

\subsection{Producing real ideal points}

Here, we show that the Conway boundary of every torus- or twist-trivial tangle variant is incompressible. It follows that any gluing of two such tangles forming a knot has an essential Conway sphere. This is the final fact needed to apply Theorem \ref{thm:main}, and to be in the situation of Section \ref{sec:trans}. 

\begin{lemma}\label{lma:essential}
Let $T$ be a torus- or twist-trivial tangle variant, and let $C' \subset \partial T$ be its Conway boundary. Then $C'$ is incompressible in $T$. Consequently, if two such variants are glued along their Conway boundaries to form a knot exterior $M = T_1 \cup_{C'}T_2$, then $C'$ is essential in $M$. 
\end{lemma}

\begin{proof}
Suppose that $D$ is a compressing disk for $C'$ in $T$, so $\partial D$ is a simple closed curve in $C'$ which does not bound a disk in $C'$. Such a simple closed curve in a four-holed sphere is either parallel to a boundary component or separates the four boundary components into two pairs. 

\medskip

Suppose $\partial D$ is parallel to a boundary component. Then $[\partial D] \in \pi_1(T)$ is conjugate to an endpoint meridian or its inverse, which the disk $D$ renders trivial in $\pi_1(T)$. On the other hand, we have $H_1(T; \mathbb{Z}) \cong \mathbb{Z}^2$, freely generated by the two strand meridian classes; in particular every strand meridian is nonzero, a contradiction. 

\medskip

Suppose $\partial D$ separates the punctures into pairs. The properly embedded disk $D$ separates $B^3$ into two balls. If each side contained one endpoint of each strand, then each strand would meet $D$, which is impossible since $D$ is disjoint from the strands. Hence $\partial D$ separates the endpoint pairs of the two strands, and $D$ splits the tangle. Write $(B_i, \alpha_i)$ for the two resulting one-string tangles with exteriors $E_i$. Then $T = E_1 \cup_DE_2$ and by van Kampen's theorem, $\pi_1(T) \cong \pi_1(E_1) * \pi_1(E_2)$. The two factors embed in the free group $\pi_1(T) \cong F_2$, and are therefore free. A meridian of $\alpha_i$ normally generates $\pi_1(E_i)$, since filling in the strand recovers the ball $B_i$. Thus, quotienting $\pi_1(T)$ by the normal closure of the trivial-strand meridian kills one free factor and leaves the other. However, we know that this quotient has abelianization $\mathbb Z$ and is nonabelian, since it admits the displayed irreducible representation. A free group with the abelianization $\mathbb Z$ is infinite cyclic and hence abelian, a contradiction. Thus $C'$ is incompressible in $T$.

\medskip

Now, let $M = T_1 \cup_{C'}T_2$ be a knot exterior obtained by gluing two such variants. Any compressing disk for $C'$ in $M$ has interior contained in one component of $M \setminus C'$, and would hence give a compressing disk in one of the tangles. Thus $C'$ is incompressible in $M$. The knot exterior $M$ is irreducible, and $\partial C'$ is contained in the torus $\partial M$. By Lemma 1.10 of \cite{hatcher}, a connected incompressible surface with boundary on torus boundary components is either essential or a boundary-parallel annulus. Since $C'$ is a four-punctured sphere, it is not an annulus. Hence $C'$ is essential in $M$. 
\end{proof}

We now prove the general detection result which counts the number of ideal points to which Theorem \ref{mainthm:a} applies.

\begin{corollary}\label{cor:generaldetection}
Let $M = T_1 \cup_CT_2$ be a knot complement in $S^3$ where $T_i$ are torus- or twist-trivial tangles. For every pair of real roots of $P_n$ or $\Lambda_{2m-1}$, there is a distinct real ideal point detecting the Conway sphere, with the corresponding limiting characters and two-sided real branch. If $T_i$ are both torus-trivial tangles, there are $nm$ such ideal points. The number of admissible data in the mixed twist-torus case is the parameter of the torus-trivial side. If $T_i$ are both twist-trivial tangles, there is one such ideal point.
\end{corollary}

\begin{proof}
Hypotheses (1)-(4) of Theorem \ref{thm:main} are properties of the individual tangle variants, verified throughout the previous three subsections. By Lemma \ref{lma:essential}, the Conway sphere is essential, so $M$ satisfies the standing hypotheses. Each pair of roots yields an admissible real limiting datum: the boundary characters agree since every boundary word has trace 2 at such a datum, and the limiting boundary representations are nonconjugate by Lemma \ref{lma:unrealized}. Theorem \ref{thm:main} and Lemma \ref{lma:detection} now produce the ideal point and its detection of $C$. Indeed, in the normalized models, $r = \tr(A_rB_r) - 2$, $z = \tr(A_zT_z) - 2$, so distinct roots determine distinct quotient characters, and hence distinct limiting characters and ideal points. Conversely, every admissible representation kills the trivial-strand meridian and hence factors through the corresponding torus- or twist-knot quotient. After conjugating the two surviving meridians to the standard form, the quotient relation forces the remaining parameter to be a root of $P_n$, respectively $\Lambda_{2m-1}$. Lemmas \ref{lma:roots} and \ref{lma:twistroots} list all real roots. Thus the preceding data exhaust the admissible real limiting data of Definition \ref{def:admissible}.
\end{proof}

%% file: families.tex
We now define the main families of knots considered in this paper and prove Theorem \ref{mainthm:c}.

\begin{figure}[h]
	\centering
	\begin{subfigure}{0.32\textwidth}
		\includegraphics[width=\linewidth]{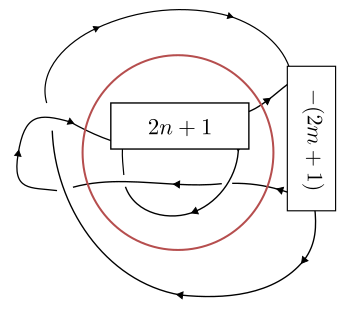}
		\caption{torus-torus: $T^{tor}_n \cup_C \overline{T^{tor}_m}$}\label{fig:familieswhiskers-a}
	\end{subfigure}\hfill
	\begin{subfigure}{0.32\textwidth}
		\includegraphics[width=\linewidth]{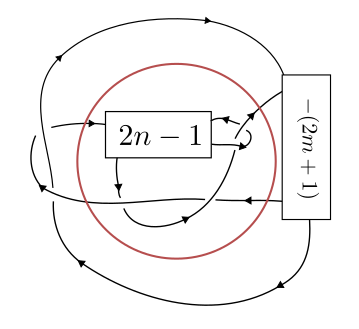}
		\caption{twist-torus: $T^{tw}_n \cup_C \overline{\bigl(T^{tor}_m\bigr)^{\mathrm{fl}}}$}\label{fig:familieswhiskers-b}
	\end{subfigure}\hfill
	\begin{subfigure}{0.32\textwidth}
		\includegraphics[width=\linewidth]{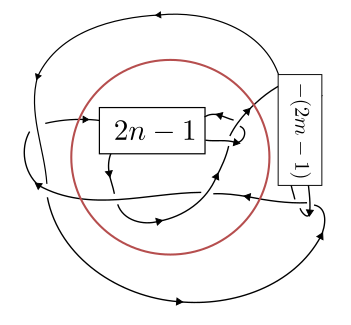}
		\caption{twist-twist: $T^{tw}_n \cup_C \overline{T^{tw}_m}$}\label{fig:familieswhiskers-c}
	\end{subfigure}
	\caption{The base torus-torus, twist-torus, and twist-twist knots, with parameters $n, m \geq 1$. The flipped counterparts of these families flip the clasps at the bottom of the interior of the red ball, and to the left in the exterior of the red ball. Knots are oriented as drawn.}\label{fig:families}
\end{figure}

\begin{definition}\label{def:families}
Let $n, m \geq 1$. The three family types, each with a \emph{base} and \emph{flipped-pair} class, are:
\begin{align*}
	\text{torus-torus:} \quad
	& T^{tor}_n \cup_C \overline{T^{tor}_m}
	&& \bigl(T^{tor}_n\bigr)^{\mathrm{fl}} \cup_C
	\overline{\bigl(T^{tor}_m\bigr)^{\mathrm{fl}}} \\
	\text{twist-torus:} \quad
	& T^{tw}_n \cup_C \overline{\bigl(T^{tor}_m\bigr)^{\mathrm{fl}}}
	&& \bigl(T^{tw}_n\bigr)^{\mathrm{fl}} \cup_C
	\overline{T^{tor}_m} \\
	\text{twist-twist:} \quad
	& T^{tw}_n \cup_C \overline{T^{tw}_m}
	&& \bigl(T^{tw}_n\bigr)^{\mathrm{fl}} \cup_C
	\overline{\bigl(T^{tw}_m\bigr)^{\mathrm{fl}}}
\end{align*}
In each case the first tangle variant occupies the inside ball, the mirrored variant the outside ball, and the boundary Conway spheres are identified as in Figure \ref{fig:families}, matching the strand endpoints so that the union of the four strand arcs is a single circle. In each type, the flipped-pair gluing is obtained from the base gluing by flipping both slots. The resulting knots are the \emph{$(n, m)$ knots} of their type and class, collectively the \emph{family knots}. 
\end{definition}

For any family knot, a four-longitude decomposition is obtained as follows. Let $\ell_1$ begin at the puncture identified as $\mu_{1, 1} \sim \mu_{2, 4}$ in Section \ref{sec:idealpoints}. Then traverse toward $\mu_{1, 2} \sim \mu_{2, 1}$; this forms $\ell_1$. Comparing to Figures \ref{fig:torustrivial} and \ref{fig:twisttrivial} will show that this traversal is precisely $l_1$. Define $\ell_2$, $\ell_3$, $\ell_4$ in the following traversal order. In general, $\ell_1$, $\ell_2$ will correspond to $l_1$ of the tangle and $\ell_3$, $\ell_4$ will correspond to $l_2$. Note that the meridional adjustments on $\ell_3$ and $\ell_4$ make the total writhe 0, thus making the final product $\lambda = \ell_1\ell_2\ell_3\ell_4$ a homological longitude. Call this particular selection the \emph{canonical four-longitude decomposition}. This is demonstrated in the below diagram.

\begin{figure}[h]
	\centering
	\includegraphics[scale=.8]{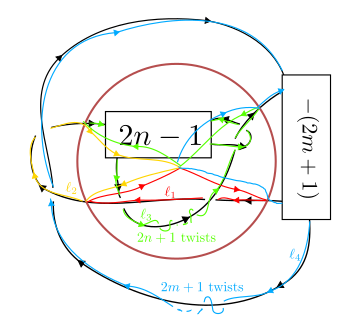}
	\caption{The canonical four-longitude decomposition for a base twist-torus family knot.}
\end{figure}

The order of the canonical four-longitude decomposition is thus
\begin{equation*}
	(\ell_1, \ell_2, \ell_3, \ell_4) = (l_{1, in}, l_{1, out}, l_{2, in}, l_{2, out})
\end{equation*}

\begin{prop}\label{prop:matching-template}
Let $K$ be a family knot with exterior $M$, and let $\lambda = \ell_1\ell_2\ell_3\ell_4$ be its canonical four-longitude decomposition. For every admissible real limiting datum described in Section \ref{sec:twobridge}, this decomposition is in positive parabolic normal form along both associated real half-branches. Moreover, for either half-branch $B$, there is a sign $\sigma_B \in \{\pm 1\}$ such that, for all sufficiently small $u > 0$:
\begin{enumerate}
	\item For a torus-torus datum $(r_{n,j}, r_{m,k})$ $\mathrm{trans}(\widetilde\rho_u(\lambda)) = 2\sigma_B(j-k)$. In particular, height cancellation holds if and only if $j=k$.
	\item For a twist-torus datum whose torus root is $r_{m,k}$, $\mathrm{trans}(\widetilde\rho_u(\lambda)) = -2\sigma_B k$. In particular, height cancellation holds if and only if $k = 0$. 
	\item For the unique twist-twist datum, $\mathrm{trans}(\widetilde\rho_u(\lambda)) = 0$. 
\end{enumerate}
\end{prop}

\begin{proof}
The first two longitudinal pieces are the identity-strand pieces of the inside and outside variants, while the last two are the knot-strand pieces, and all longitudinal pieces are oriented as they are displayed in Figures \ref{fig:torustrivial} and \ref{fig:twisttrivial}. Thus, the relevant data are precisely those computed in Section \ref{sec:twobridge}. The identity-strand pieces are conjugate in the quotient models to meridians or inverse meridians, so their quotient-induced normalized lifts satisfy $\eta_1 = \eta_2 = 0$. The arrows in the family diagrams agree with the local traversal orientations used in Section \ref{sec:twobridge}, including the mirror words, so no additional inversion of a piece is required. The parameters $\beta_i$ are nonzero and the central signs of $\ell_3$ and $\ell_4$ are $-1$, while the central signs of $\ell_1$ and $\ell_2$ are $+1$. We are thus in the situation of Theorem \ref{thm:modelcriterion}. The following table, assembled from the data in Section \ref{sec:twobridge}, contains all of the relevant information.
\[
\begin{array}{c|cc|cc|cc}
	\text{family class}
	&
	c_1&\alpha_1
	&
	c_2&\alpha_2
	&
	\eta_3&\eta_4
	\\ \hline
	\text{torus--torus, base}
	&
	1&r_{n,j}
	&
	1&-r_{m,k}
	&
	2j+1&-(2k+1)
	\\
	\text{torus--torus, flipped}
	&
	1&-r_{n,j}
	&
	1&r_{m,k}
	&
	2j+1&-(2k+1)
	\\
	\text{twist--torus, base}
	&
	-z_n&-z_n
	&
	1&r_{m,k}
	&
	1&-(2k+1)
	\\
	\text{twist--torus, flipped}
	&
	-z_n&z_n
	&
	1&-r_{m,k}
	&
	1&-(2k+1)
	\\
	\text{twist--twist, base}
	&
	-z_n&-z_n
	&
	-z_m&z_m
	&
	1&-1
	\\
	\text{twist--twist, flipped}
	&
	-z_n&z_n
	&
	-z_m&-z_m
	&
	1&-1 .
\end{array}
\]
The root calculations of Section \ref{sec:twobridge} give $r_{n,j} < 0$, $r_{m,k} < 0$, $z_n < 0$, $z_m < 0$. Thus all $c_i$ in the table are positive, while $\alpha_1$ and $\alpha_2$ have opposite signs in every row. Consequently, $c_1c_2\alpha_1\alpha_2 < 0$ in all family classes, and Theorem \ref{thm:modelcriterion} gives positive parabolic normal form along both half-branches. Moreover, $c_1c_2 > 0$ in every row. The height criterion of Theorem \ref{thm:modelcriterion} therefore gives, for some $\sigma_B\in \{\pm 1\}$, $(h_1, h_2, h_3, h_4) = \sigma_B(0, 0, \eta_3, \eta_4)$. Lemma \ref{lma:trans} now gives $\mathrm{trans}(\widetilde\rho_u(\lambda)) = \sigma_B(\eta_3 + \eta_4)$. For a torus-torus datum this is $\sigma_B((2j+1) - (2k+1)) = 2\sigma_B(j-k)$. For a twist-torus datum, it is $\sigma_B(1-(2k+1)) = -2\sigma_Bk$. For the twist-twist datum, it is zero. 
\end{proof}

\begin{proof}[Proof of Theorem \ref{mainthm:c}]
(1) is Corollary \ref{cor:generaldetection}. For (2), Corollary \ref{cor:realbranches} gives $\mathrm{tr}\rho_u(\mu) > 2$ on both half-branches associated to every admissible datum. Proposition \ref{prop:matching-template} places the canonical four-longitude decomposition in positive parabolic normal form. Theorem \ref{thm:staggeredtraces} therefore gives $|\tr\rho_u(\lambda)| \to \infty$. Corollary \ref{cor:sheets} shows that the meridian-normalized lifts define two unbounded augmented holonomy arcs with vertical asymptotes on opposite sides. For (3), in the torus-torus case, the admissible data are indexed by $0 \leq j < n$, $0 \leq k < m$. By Proposition \ref{prop:matching-template}, height cancellation holds exactly when $j=k$, giving $\min(n,m)$ branches. In the twist-torus case it holds exactly for $k=0$, giving one branch, and in the twist-twist case it holds for the unique datum. For these and only these constructed branches, the longitudinal translation number is zero. Corollary \ref{cor:h00} places both associated augmented arcs in $H_{0,0}(M)$, and Corollary \ref{cor:slopes} implies that there exist $R_- < R_+$ such that $\pi_1(M(r))$ is left-orderable for every rational $r \notin [R_-, R_+]$.
\end{proof}

%% file: future.tex
The overall method of this paper has two distinct inputs. The first is that an admissible real limiting datum satisfying the four hypotheses of Theorem \ref{thm:main} produces a real ideal point and a two-sided real branch. The second is that an exact four-longitude decomposition, together with the sign and height criterion of Theorem \ref{thm:modelcriterion}, determines the longitudinal translation number and hence whether the associated augmented arcs lie in $H_{0, 0}$. Future extensions of the methods in this paper may enlarge either input, or may vary the way in which the two local tangle models are glued. 

\subsection{Other tangles}

The torus-trivial and twist-trivial tangles considered in this paper provide particularly accessible examples because their groups are free of rank two, killing one strand produces a two-bridge knot group, and the relevant real boundary-parabolic characters can be described explicitly. None of these features is intrinsic to any of the arguments in this paper. It would be useful to characterize tangle exteriors which admit admissible real limiting representations which satisfy hypotheses (1)-(4) of Theorem \ref{thm:main}. It would also be interesting to prove that any of these conditions hold all of the time or generically; hypothesis (3) seems to be the most difficult to assess, while the other three have natural geometric explanations. It would also be useful to determine which of these tangle exteriors admit exact longitudinal pieces satisfying the criterion of Theorem \ref{thm:modelcriterion}.

\medskip

A broader extension would allow limiting boundary data outside Definition \ref{def:admissible}. For example, one may ask whether a comparable real node construction exists when the killed-strand image is nontrivial, when one limiting tangle restriction is reducible, or when the natural limiting representation does not factor through a two-bridge knot quotient. Such examples would require a generalization of the nonrealized gluing datum used here.

\subsection{The Conway knot} 

The Conway knot $K11n34$ appears to provide a particularly concrete example beyond the symmetric setting of this paper. In its computed holonomy extension locus, there are the expected two arcs tending to slope $\pm\infty$. However, only the arc for which the slope $y/x$ tends to $-\infty$ lies in $H_{0,0}$, while the other has nonzero longitudinal translation number. The computed locus of the Kinoshita-Terasaka knot $K11n42$ coincides with that of $K11n34$. The two knots are Conway mutants, and mutation is known to preserve substantial character-variety data, but we know of no theorem forcing equality of their holonomy extension loci; see \cite{tillmann}. 

\begin{figure}[ht]
	\centering
	\includegraphics[scale=.4]{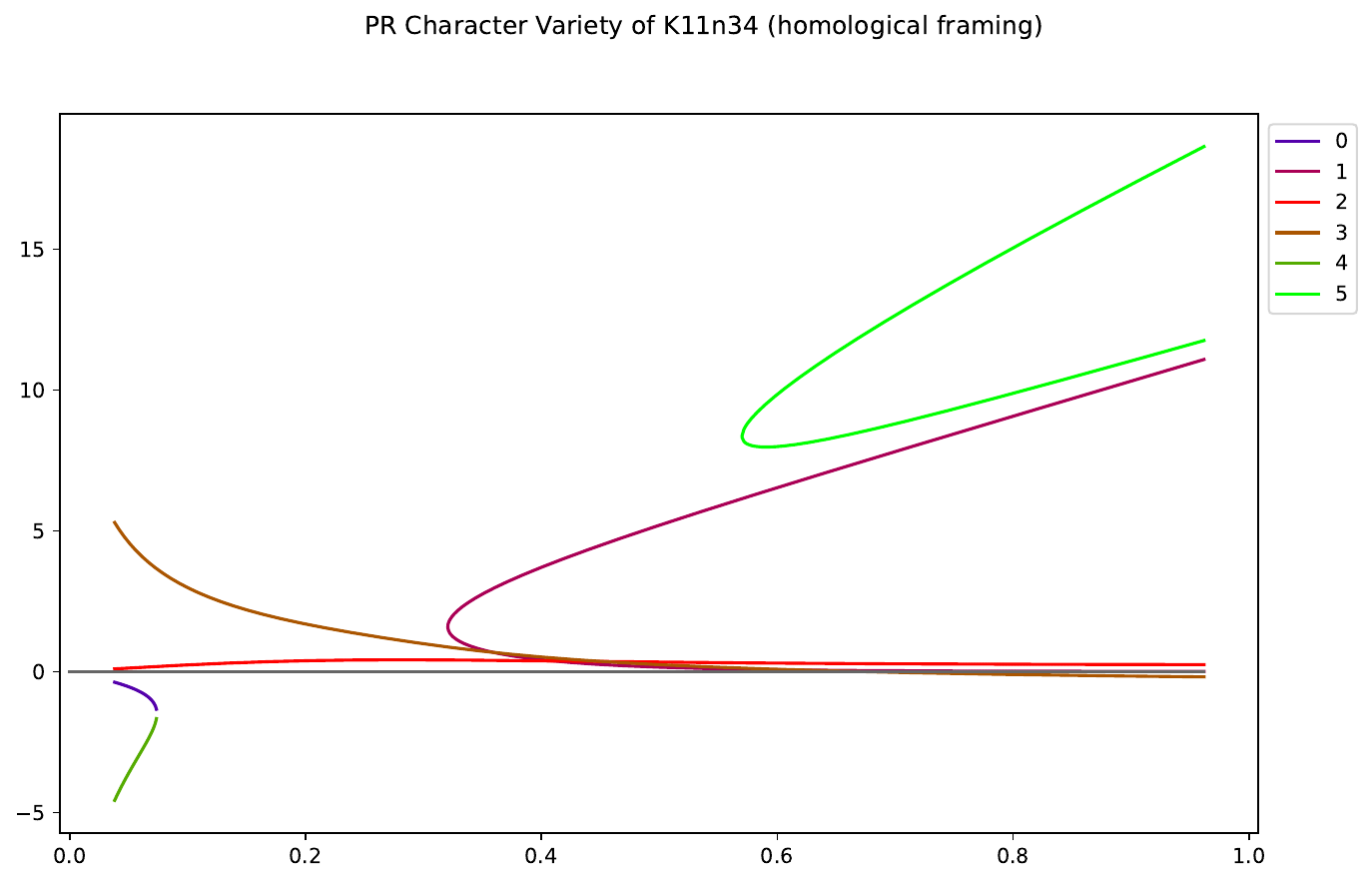}
	\caption{The locus for the Conway knot K11n34. Only the arc with slope approaching $-\infty$ is computed to have translation number $0$. The computed locus for the Kinoshita--Terasaka knot K11n42 is identical.}
	\label{fig:conway}
\end{figure}

\begin{conjecture}\label{conj:conway}
Let $M$ be the exterior of the Conway knot $K11n34$ or of the Kinoshita-Terasaka knot $K11n42$. There exists $a \in \mathbb{Q}$ such that $\pi_1(M(r))$ is left-orderable for every $r \in \mathbb{Q}$ with $r > a$.
\end{conjecture}

A proof would require a one-sided version of our peripheral criteria: one real half-branch would have longitudinal translation number zero, while the other would be allowed to lie on a nonzero sheet of the holonomy extension locus. To take things further, we ask if there is a general family of knots or tangles beginning with the Conway or Kinoshita-Terasaka knots and tangles for which we can show orderability for every sufficiently large-slope Dehn filling, similarly to how the family knots in this paper grew out of torus-trivial and twist-trivial tangles.

\subsection{Several and higher-punctured spheres}

The construction in this paper detects and exploits one chosen essential Conway sphere. For an exterior containing several disjoint Conway spheres, one could seek compatible limiting representations on all complementary tangle regions and a simultaneous real gluing deformation. The local character variety would then involve several interacting nodal equations, while the translation-number part of the construction would require summing multiple longitudinal pieces across the entire decomposition graph.

\medskip

More generally, one may consider essential $2n$-punctured spheres for $n > 2$. In \cite{lmz} the authors speculate that an $L$-space knot is $n$-string prime for every $n$. Combined with the $L$-space conjecture, this predicts that a knot exterior containing an essential $2n$-punctured sphere should have left-orderable fundamental group for every rational filling. The case $n = 1$ is the theorem of Boyer-Gordon-Hu \cite{boyerhu} for composite knots. Conjecture \ref{conj:main} is the case $n = 2$. 

\medskip

A first step for $n > 2$ is to construct real ideal points detecting such surfaces and to understand their limiting characters. The Fricke calculations in Section \ref{sec:idealpoints} would have to be replaced by higher-dimensional boundary character varieties; in particular, the boundary character variety would be the character variety of $\pi_1(S^2 \setminus 2n \text{ points}\}) \cong F_{2n-1}$. If real branches can be constructed, the remaining problem is to compute the translation number of the homological longitude along those branches. This preserves the two-step structure of this paper and offers a possible route to progress on the general $n$-string version of Conjecture \ref{conj:main}. 

%% file: computations.tex
\section{Computations in the proof of Theorem \ref{thm:main}}

\begin{lemma}\label{lma:discriminantterms}
The constant, linear, and quadratic homogeneous parts of $\Delta_1$, $\Delta_2$ vanish. The cubic homogeneous parts are
\begin{equation*}
	\Delta_1^{(3)} = 4X_0\bigl(P_0^2 - 2P_0R_0 - 2P_0Y_0 + R_0^2 - 2R_0Y_0 + X_0Y_0 + Y_0^2\bigr)
\end{equation*}
\begin{equation*}
	\Delta_2^{(3)} = 4Y_0\bigl(P_0^2 - 2P_0R_0 - 2P_0X_0 + R_0^2 - 2R_0X_0 + X_0^2 + X_0Y_0\bigr)
\end{equation*}
\end{lemma}

\begin{proof}
Regard $\Delta_1$ as a polynomial in $x$ alone:
\begin{equation*}
	\Delta_1 = (y^2-4)x^2 + \bigl(4p^2+4r^2-4pry\bigr)x + \bigl(-8p^2+8pry-8r^2-4y^2+16\bigr) =: Ax^2+Bx+C.
\end{equation*}
Setting $x=2$ in the underlying quadratic $z^2+(xy-2pr)z+G_1$ gives $z^2+2(y-pr)z+(y-pr)^2=(z+y-pr)^2$, a perfect square; hence $\Delta_1|_{x=2}=0$, i.e.\ $4A+2B+C=0$, and $(x-2)$ divides $\Delta_1$. Synthetic division gives quotient $\Psi_1 := Ax+(2A+B)$ with zero remainder:
\begin{equation*}
	\Delta_1 = (x-2)\Psi_1, \qquad \Psi_1 = 4p^2-4pry+4r^2+xy^2-4x+2y^2-8.
\end{equation*}
Substituting $p=2+P_0$, $r=2+R_0$, $x=2+X_0$, $y=2+Y_0$ term by term,
\begin{equation*}
	4p^2 = 16+16P_0+4P_0^2, \qquad 4r^2 = 16+16R_0+4R_0^2,
\end{equation*}
\begin{equation*}
	-4pry = -32-16P_0-16R_0-16Y_0-8P_0R_0-8P_0Y_0-8R_0Y_0-4P_0R_0Y_0,
\end{equation*}
\begin{equation*}
	xy^2 = 8+4X_0+8Y_0+2Y_0^2+4X_0Y_0+X_0Y_0^2, \qquad -4x=-8-4X_0, \qquad 2y^2 = 8+8Y_0+2Y_0^2, \qquad -8=-8.
\end{equation*}
Summing, the constant terms give $16+16-32+8-8+8-8=0$; the linear terms give $(16P_0-16P_0)+(16R_0-16R_0)+(-16Y_0+8Y_0+8Y_0)+(4X_0-4X_0)=0$; the quadratic terms give
\begin{equation*}
	\Psi_1^{(2)} = 4P_0^2-8P_0R_0-8P_0Y_0+4R_0^2-8R_0Y_0+4X_0Y_0+4Y_0^2;
\end{equation*}
the cubic terms give $\Psi_1^{(3)} = -4P_0R_0Y_0+X_0Y_0^2$. Since $\Delta_1 = X_0\Psi_1$ exactly and $\Psi_1^{(0)} = \Psi_1^{(1)}=0$, the polynomial $\Delta_1$ vanishes to order three and $\Delta_1^{(3)} = X_0\Psi_1^{(2)}$, as desired. For $\Delta_2$: since $G_2(p,r,x,y) = G_1(p,r,y,x)$ and $2pr-xy$ is symmetric in $x,y$, we have $\Delta_2(p,r,x,y) = \Delta_1(p,r,y,x)$. The statements for $\Delta_2$, including $\Delta_2^{(3)}$, follow from those for $\Delta_1$ by exchanging $X_0 \leftrightarrow Y_0$.
\end{proof}

\begin{lemma}\label{lma:crossjet}
Suppose $\tr(Y_i) = \tr(Y_i^\tau) = Q_i$, $\tr(X_i^2) = \tr((X_i^\tau)^2) = 2Q_i$, and $\tr(X_iX_i^\tau) = -2Q_i$. Then $\tr\rho_{i,s}(m_{i,1}^{-1}m_{i,2}) = 2 + 4Q_is^2 + O(s^3)$. 
\end{lemma}

\begin{proof}
Write 
\begin{equation*}
	M = I + sX_i + s^2Y_i + O(s^3) \ \ \ \ \ K = I + sX_i^\tau + s^2Y_i^\tau + O(s^3)
\end{equation*}
Then $M^{-1} = I - sX_i + s^2(X_i^2 - Y_i) + O(s^3)$, and 
\begin{equation*}
	M^{-1}K = I + s(X_i^\tau - X_i) + s^2(Y_i^\tau - X_iX_i^\tau + X_i^2 - Y_i) + O(s^3)
\end{equation*}
The linear term of the trace vanishes since $X_i$ and $X_i^\tau$ are traceless. The quadratic term is 
\begin{equation*}
	\tr(Y_i^\tau) - \tr(X_iX_i^\tau) + \tr(X_i^2) - \tr(Y_i) = Q_i - (-2Q_i) + 2Q_i - Q_i = 4Q_i
\end{equation*}
\end{proof}

\begin{remark}
Since $\tr(A^{-1}) = \tr(A)$ for any $A \in SL_2(\mathbb R)$, we have, for any $A, B \in SL_2(\mathbb R)$,
\begin{equation*}
	\tr(AB^{-1}) = \tr((AB^{-1})^{-1}) = \tr(BA^{-1}) = \tr(A^{-1}B)
\end{equation*}
by cyclicity of trace. Consequently, Lemma \ref{lma:crossjet} applies equally to $y_2 = \tr(\mu_{2,1}\mu_{2,2}) = \tr(m_{2,1}m_{2,2}^{-1}) = \tr(m_{2,1}^{-1}m_{2,2})$ with $i = 2$, giving $y_2 = 2 + 4Q_2s^2 + O(s^3)$.
\end{remark}

\section{Proof of Proposition \ref{prop:branchmodel}}\label{sec:algebra}

\begin{proof}
After shrinking, $B$ is a connected component of $(\widetilde{\mathcal C}(\mathbb R) \cap U) \setminus \{x\}$ for a small semialgebraic neighborhood $U$ of $x$; hence $B$ is semialgebraic by Proposition 2.4.5 of \cite{bochnak-coste-roy}. Let $R(M) = \mathrm{Hom}(\pi_1(M), SL_2(\mathbb C))$, let $\tau: R(M) \to X(M)$ be the character map, and let $\pi: \widetilde{\mathcal C} \to \overline{\mathcal C}$ be the normalization. In a real projective compactification of $R(M)$, consider the semialgebraic set 
\begin{equation*}
	\mathcal S = \{(z, \rho) \in B \times R(M)(\mathbb R) \mid \tau(\rho) = \pi(z)\} \subset \widetilde{\mathcal C}(\mathbb R) \times \overline{R(M)}(\mathbb R)
\end{equation*}
By Corollary \ref{cor:realbranches}, its projection to $B$ is surjective. Choose $z_n \in B$ with $z_n \to x$ and $(z_n, \rho_n) \in \mathcal S$. After passing to a subsequence, projective compactness gives $\rho_n \to y$, so $(x, y) \in \overline{\mathcal S}$. Proposition 8.1.13 of \cite{bochnak-coste-roy} gives a Nash arc $(z(s), \rho(s))$ approaching that point, with $z(s) \in B$ and $\chi_{\rho(s)} = \pi(z(s))$ for $s > 0$. Let $D$ be the Zariski closure of $\rho((0, \delta))$. Its vanishing ideal is prime, and its coordinate functions are algebraic over $\mathbb R(s)$. Since its character is nonconstant, $D$ is an irreducible real algebraic curve. Let $\widetilde D$ be its smooth projective model. Since the characters of the arc lie in $\mathcal C$ and are nonconstant, the image of $D$ is Zariski dense in $\mathcal C$. The construction of Section 1.3 of \cite{cullershalen} gives a tautological representation $P: \pi_1(M) \to SL_2(F)$, where $F = \mathbb R(D)$, and the character map extends to a nonconstant finite morphism $\widetilde\tau: \widetilde D \to \widetilde{\mathcal C}$. After shrinking, the arc $\rho(s)$ avoids the finite locus over which the normalization of the projective closure of $D$ is not an isomorphism, and hence lifts to an arc $\widetilde\rho(s) \in \widetilde D(\mathbb R)$. By properness, this arc converges to a real point $\widetilde y \in \widetilde D(\mathbb R)$. For sufficiently small $s>0$, the character of $\rho(s)$ avoids the finite locus over which $\pi$ is not injective. Hence $\widetilde\tau(\widetilde\rho(s))=z(s)$, and continuity gives $\widetilde\tau(\widetilde y)=x$. A real local parameter $u$ at $\widetilde y$ embeds $F$ into $\mathbb R((u))_{\mathrm{conv}}$, and the resulting representation is $\rho_B$. If $t$ is the signed parameter at $x$, then $\widetilde\tau^*t = u^e\eta(u)$ for some $e \geq 1$ and an analytic unit $\eta$. After a real analytic change of $u$, this becomes $t = \sigma u^e$, where $\sigma$ records the chosen half-branch. 

\medskip

Let $A = \mathcal O_{\widetilde D, \widetilde y}$. Since $\widetilde y$ is a real smooth point of $\widetilde D$, $A$ is a discrete valuation ring with residue field $\mathbb R$. For every $\gamma \in \pi_1(T_i)$, $\tr P(\gamma) = \widetilde\tau^{*}I_\gamma$ is regular at $\widetilde y$, because $I_\gamma$ is finite at $x$; its residue is $\chi_i(\gamma)$. Thus $\tr P(\gamma) \in A$. Since the specialization character $\chi_i$ is irreducible, the generic restriction of $P$ to $\pi_1(T_i)$ is irreducible, and Lemma 1.3.1 of \cite{cullershalen} shows that $P|_{\pi_1(T_i)}$ is absolutely irreducible. Lemma 1.4.3 of \cite{cullershalen} therefore gives $G_i \in GL_2(F)$ such that $G_iP(\pi_1(T_i))G_i^{-1} \subset SL_2(A)$. After shrinking $B$, $\det(G_i)$ has constant sign. Multiplying by $\mathrm{diag}(1, -1)$ if necessary makes this determinant positive. Because $\det G_i$ has positive leading coefficient along $B$, its positive square root lies in a finite real extension. Adjoining the positive square roots of $\det G_1$ and $\det G_2$ gives a finite real Puiseux extension $K/\mathbb R((u))_{\mathrm{conv}}$ with residue field $\mathbb R$. Then $C_i = (\det G_i)^{-1/2}G_i \in SL_2(K)$ has the same conjugation action, and reduction gives $\rho_i^B: \pi_1(T_i) \to SL_2(\mathbb{R})$ with character $\chi_i$. Since $\chi_i$ is irreducible, the reduced representation $\rho_i^B$ is irreducible.
\end{proof}

\section{Computations in the proof of Lemma \ref{lma:residuehalf}}

\begin{lemma}\label{lma:commutatorasymptotics}
Let $p, b \in \mathbb R^\times$, let $Z = \left(\begin{smallmatrix}\alpha&\beta\\\zeta&\delta\end{smallmatrix}\right)$ be an arbitrary real $2\times 2$ matrix, and let $Y = \left(\begin{smallmatrix}\upsilon&\eta\\0&-\upsilon\end{smallmatrix}\right)$ be traceless with $Y_{21} = 0$. Put $P = \left(\begin{smallmatrix}1&p\\0&1\end{smallmatrix}\right)$ and $P_b = \left(\begin{smallmatrix}1&b\\0&1\end{smallmatrix}\right)$, and let
\begin{equation*}
	M(t) = \bigl(I + tY + t^2Y^{(2)} + O(t^3)\bigr)P \ \ \ \ \
	K(t) = \bigl(I + tZ + t^2Z^{(2)} + O(t^3)\bigr)P_b
\end{equation*}
be real analytic, with $Y^{(2)}, Z^{(2)}$ arbitrary. Then
\begin{equation*}
	M^{-1}K^{-1}MK = I + t\begin{pmatrix}p\zeta & * \\ 0 & -p\zeta\end{pmatrix} + O(t^2) \ \ \ \ \
	MK^{-1}M^{-1}K = I + t\begin{pmatrix}-p\zeta & * \\ 0 & p\zeta\end{pmatrix} + O(t^2)
\end{equation*}
Moreover,
\begin{equation*}
	\bigl(M^{-1}K^{-1}MK\bigr)_{21} = -\zeta\bigl(p\zeta + 2\upsilon\bigr)t^2 + O(t^3)
\end{equation*}
and this coefficient is independent of $\eta$, of $\alpha, \beta, \delta$, of $b$, and of $Y^{(2)}, Z^{(2)}$.
\end{lemma}

\begin{proof}
For a matrix $A(t) = A_0 + tA_1 + t^2A_2 + O(t^3)$ with $A_0$ invertible,
\begin{equation*}
	A(t)^{-1} = A_0^{-1} - tA_0^{-1}A_1A_0^{-1} + t^2\bigl(A_0^{-1}A_1A_0^{-1}A_1A_0^{-1} - A_0^{-1}A_2A_0^{-1}\bigr) + O(t^3)
\end{equation*}
Applying this to $M(t) = P + tYP + t^2Y^{(2)}P$ and simplifying with $P^{-1}(YP)P^{-1} = P^{-1}Y$ gives
\begin{equation*}
	M(t)^{-1} = P^{-1} - tP^{-1}Y + t^2\bigl(P^{-1}Y^2 - P^{-1}Y^{(2)}\bigr) + O(t^3)
\end{equation*}
and likewise for $K(t)^{-1}$ with $Y, p, Y^{(2)}$ replaced by $Z, b, Z^{(2)}$. Since $P$ and $P_b$ are both upper unitriangular they commute, so both $M^{-1}K^{-1}MK$ and $MK^{-1}M^{-1}K$ have constant term $I$ and hence begin at order $t$.
	
\medskip
	
A first-order multiplication gives
\begin{equation*}
	[t]\bigl(M^{-1}K^{-1}MK\bigr) = \begin{pmatrix}p\zeta & p(\delta-\alpha)+p^2\zeta+2b(p\zeta+\upsilon)\\ 0&-p\zeta \end{pmatrix}
\end{equation*}
and
\begin{equation*}
	[t]\bigl(MK^{-1}M^{-1}K\bigr) = \begin{pmatrix} -p\zeta & p(\alpha-\delta)+p^2\zeta-2b(p\zeta+\upsilon)\\ 0&p\zeta \end{pmatrix}
\end{equation*}
This proves the two first-order assertions.

\medskip

For the second-order calculation, write
\begin{equation*}
	M^{-1}=A_0+tA_1+t^2A_2+O(t^3) \ \ \ \ \
	K^{-1}=B_0+tB_1+t^2B_2+O(t^3)
\end{equation*}
\begin{equation*}
	M=C_0+tC_1+t^2C_2+O(t^3) \ \ \ \ \
	K=D_0+tD_1+t^2D_2+O(t^3)
\end{equation*}
using the coefficients displayed above.  Put $y_{21}^{(2)}=(Y^{(2)})_{21}$, $z_{21}^{(2)}=(Z^{(2)})_{21}$. The ten contributions to $[t^2](M^{-1}K^{-1}MK)_{21}$, listed in the order
\begin{equation*}
	\begin{gathered}
		A_2B_0C_0D_0 \ \ \ \ \
		A_0B_2C_0D_0 \ \ \ \ \
		A_0B_0C_2D_0 \ \ \ \ \ 
		A_0B_0C_0D_2 \ \ \ \ \ 
		A_1B_1C_0D_0 \\
		A_1B_0C_1D_0 \ \ \ \ \ 
		A_1B_0C_0D_1 \ \ \ \ \
		A_0B_1C_1D_0 \ \ \ \ \
		A_0B_1C_0D_1 \ \ \ \ \
		A_0B_0C_1D_1
	\end{gathered}
\end{equation*}
have $(2,1)$-entries
\begin{equation*}
	\begin{gathered}
		-y_{21}^{(2)} \ \ \ \ \
		(\alpha+\delta)\zeta-z_{21}^{(2)} \ \ \ \ \
		y_{21}^{(2)} \ \ \ \ \
		z_{21}^{(2)} \ \ \ \ \
		-\upsilon\zeta \\
		0 \ \ \ \ \
		\upsilon\zeta \ \ \ \ \
		-\upsilon\zeta \ \ \ \ \
		-(\alpha+\delta)\zeta-p\zeta^2 \ \ \ \ \
		-\upsilon\zeta
	\end{gathered}
\end{equation*}
Summing gives
\begin{equation*}
	[t^2]\bigl(M^{-1}K^{-1}MK\bigr)_{21} = -p\zeta^2-2\upsilon\zeta = -\zeta(p\zeta+2\upsilon)
\end{equation*}
This also makes explicit the cancellation of all dependence on $\eta,\alpha,\beta,\delta,b,Y^{(2)}$, and $Z^{(2)}$.
\end{proof}

%% file: conway.bib
@article {sikora,
	AUTHOR = {Sikora, Adam S.},
	TITLE = {Character varieties},
	JOURNAL = {Trans. Amer. Math. Soc.},
	FJOURNAL = {Transactions of the American Mathematical Society},
	VOLUME = {364},
	YEAR = {2012},
	NUMBER = {10},
	PAGES = {5173--5208},
	ISSN = {0002-9947,1088-6850},
	MRCLASS = {14D20 (14L24 53D30 57M50)},
	MRNUMBER = {2931326},
	MRREVIEWER = {Benjamin\ M. S. Martin},
	DOI = {10.1090/S0002-9947-2012-05448-1},
	URL = {https://doi.org/10.1090/S0002-9947-2012-05448-1},
}

@article {lmz,
	AUTHOR = {Lidman, Tye and Moore, Allison H. and Zibrowius, Claudius},
	TITLE = {{$L$}-space knots have no essential {C}onway spheres},
	JOURNAL = {Geom. Topol.},
	FJOURNAL = {Geometry \& Topology},
	VOLUME = {26},
	YEAR = {2022},
	NUMBER = {5},
	PAGES = {2065--2102},
	ISSN = {1465-3060,1364-0380},
	MRCLASS = {57K18 (57K30 57R58)},
	MRNUMBER = {4520302},
	MRREVIEWER = {Thilo\ Kuessner},
	DOI = {10.2140/gt.2022.26.2065},
	URL = {https://doi.org/10.2140/gt.2022.26.2065},
}

@article {paoluzziporti,
	AUTHOR = {Paoluzzi, Luisa and Porti, Joan},
	TITLE = {Conway spheres as ideal points of the character variety},
	JOURNAL = {Math. Ann.},
	FJOURNAL = {Mathematische Annalen},
	VOLUME = {354},
	YEAR = {2012},
	NUMBER = {2},
	PAGES = {707--726},
	ISSN = {0025-5831,1432-1807},
	MRCLASS = {57M25 (20C99 57M50)},
	MRNUMBER = {2965258},
	MRREVIEWER = {John\ G.\ Ratcliffe},
	DOI = {10.1007/s00208-011-0748-y},
	URL = {https://doi.org/10.1007/s00208-011-0748-y},
}

@article {gao,
	AUTHOR = {Gao, Xinghua},
	TITLE = {Orderability of homology spheres obtained by {D}ehn filling},
	JOURNAL = {Math. Res. Lett.},
	FJOURNAL = {Mathematical Research Letters},
	VOLUME = {29},
	YEAR = {2022},
	NUMBER = {5},
	PAGES = {1387--1427},
	ISSN = {1073-2780,1945-001X},
	MRCLASS = {20F34 (20F65 57K30 57M07)},
	MRNUMBER = {4589362},
	MRREVIEWER = {Michael\ Heusener},
	DOI = {10.4310/mrl.2022.v29.n5.a4},
	URL = {https://doi.org/10.4310/mrl.2022.v29.n5.a4},
}

@article {cullerdunfield,
	AUTHOR = {Culler, Marc and Dunfield, Nathan M.},
	TITLE = {Orderability and {D}ehn filling},
	JOURNAL = {Geom. Topol.},
	FJOURNAL = {Geometry \& Topology},
	VOLUME = {22},
	YEAR = {2018},
	NUMBER = {3},
	PAGES = {1405--1457},
	ISSN = {1465-3060,1364-0380},
	MRCLASS = {57M60 (20F60 57M05 57M25)},
	MRNUMBER = {3780437},
	MRREVIEWER = {Wolfgang\ H.\ Heil},
	DOI = {10.2140/gt.2018.22.1405},
	URL = {https://doi.org/10.2140/gt.2018.22.1405},
}

@article {boyerhu,
	AUTHOR = {Boyer, Steven and Gordon, Cameron McA. and Hu, Ying},
	TITLE = {J{SJ} decompositions of knot exteriors, {D}ehn surgery and the
	{$L$}-space conjecture},
	JOURNAL = {Selecta Math. (N.S.)},
	FJOURNAL = {Selecta Mathematica. New Series},
	VOLUME = {31},
	YEAR = {2025},
	NUMBER = {1},
	PAGES = {Paper No. 3, 31},
	ISSN = {1022-1824,1420-9020},
	MRCLASS = {57M50 (57K10 57M05)},
	MRNUMBER = {4833878},
	MRREVIEWER = {Quach thi C\^am V\^an},
	DOI = {10.1007/s00029-024-00999-3},
	URL = {https://doi.org/10.1007/s00029-024-00999-3},
}

@article {dunfield,
	AUTHOR = {Dunfield, Nathan M.},
	TITLE = {Examples of non-trivial roots of unity at ideal points of
	hyperbolic {$3$}-manifolds},
	JOURNAL = {Topology},
	FJOURNAL = {Topology. An International Journal of Mathematics},
	VOLUME = {38},
	YEAR = {1999},
	NUMBER = {2},
	PAGES = {457--465},
	ISSN = {0040-9383},
	MRCLASS = {57M50 (57N10)},
	MRNUMBER = {1660313},
	MRREVIEWER = {Hugh\ M.\ Hilden},
	DOI = {10.1016/S0040-9383(98)00035-4},
	URL = {https://doi.org/10.1016/S0040-9383(98)00035-4},
}

@article {tillmann,
	AUTHOR = {Tillmann, Stephan},
	TITLE = {Character varieties of mutative 3-manifolds},
	JOURNAL = {Algebr. Geom. Topol.},
	FJOURNAL = {Algebraic \& Geometric Topology},
	VOLUME = {4},
	YEAR = {2004},
	PAGES = {133--149},
	ISSN = {1472-2747,1472-2739},
	MRCLASS = {57M27 (14E05 57M50 57N10)},
	MRNUMBER = {2059186},
	MRREVIEWER = {Bruno\ P.\ Zimmermann},
	DOI = {10.2140/agt.2004.4.133},
	URL = {https://doi.org/10.2140/agt.2004.4.133},
}

@article {WangPuncturedJSJTori,
	AUTHOR = {Wang, Yi},
	TITLE = {Punctured {JSJ} tori and tautological extensions of {A}zumaya
	algebras},
	JOURNAL = {Algebr. Geom. Topol.},
	FJOURNAL = {Algebraic \& Geometric Topology},
	VOLUME = {26},
	YEAR = {2026},
	NUMBER = {6},
	PAGES = {2049--2077},
	ISSN = {1472-2747,1472-2739},
	MRCLASS = {57K31 (16H05)},
	MRNUMBER = {5092454},
	DOI = {10.2140/agt.2026.26.2049},
	URL = {https://doi.org/10.2140/agt.2026.26.2049},
}

@misc{WangLimitingCharacters,
	author        = {Wang, Yi},
	title         = {Limiting characters at ideal points detecting twice-punctured tori},
	year          = {2024},
	eprint        = {2404.06388},
	archivePrefix = {arXiv},
	primaryClass  = {math.GT},
	note          = {arXiv preprint},
	doi           = {10.48550/arXiv.2404.06388}
}

@misc{WangDetectedSeifertSurfaces,
	author        = {Wang, Yi},
	title         = {Detected {Seifert} surfaces and intervals of left-orderable surgeries},
	year          = {2025},
	eprint        = {2509.08127},
	archivePrefix = {arXiv},
	primaryClass  = {math.GT},
	note          = {arXiv preprint},
	doi           = {10.48550/arXiv.2509.08127}
}

@article {cullershalen,
	AUTHOR = {Culler, Marc and Shalen, Peter B.},
	TITLE = {Varieties of group representations and splittings of
	{$3$}-manifolds},
	JOURNAL = {Ann. of Math. (2)},
	FJOURNAL = {Annals of Mathematics. Second Series},
	VOLUME = {117},
	YEAR = {1983},
	NUMBER = {1},
	PAGES = {109--146},
	ISSN = {0003-486X},
	MRCLASS = {57N10},
	MRNUMBER = {683804},
	MRREVIEWER = {G.\ Peter\ Scott},
	DOI = {10.2307/2006973},
	URL = {https://doi.org/10.2307/2006973},
}

@article {hosteshanahan,
	AUTHOR = {Hoste, Jim and Shanahan, Patrick D.},
	TITLE = {Trace fields of twist knots},
	JOURNAL = {J. Knot Theory Ramifications},
	FJOURNAL = {Journal of Knot Theory and its Ramifications},
	VOLUME = {10},
	YEAR = {2001},
	NUMBER = {4},
	PAGES = {625--639},
	ISSN = {0218-2165,1793-6527},
	MRCLASS = {57M27 (57M25 57M50)},
	MRNUMBER = {1831680},
	MRREVIEWER = {Mattia\ Mecchia},
	DOI = {10.1142/S0218216501001049},
	URL = {https://doi.org/10.1142/S0218216501001049},
}

@article {dunfieldrasmussen,
	AUTHOR = {Dunfield, Nathan M. and Rasmussen, Jacob},
	TITLE = {A unified {C}asson-{L}in invariant for the real forms of
	{${\rm SL}(2)$}},
	JOURNAL = {Geom. Topol.},
	FJOURNAL = {Geometry \& Topology},
	VOLUME = {29},
	YEAR = {2025},
	NUMBER = {8},
	PAGES = {4055--4188},
	ISSN = {1465-3060,1364-0380},
	MRCLASS = {57K10 (20F60 53A20 57K31 57M05)},
	MRNUMBER = {4999726},
	DOI = {10.2140/gt.2025.29.4055},
	URL = {https://doi.org/10.2140/gt.2025.29.4055},
}

@incollection {goldman,
	AUTHOR = {Goldman, William M.},
	TITLE = {Trace coordinates on {F}ricke spaces of some simple hyperbolic
	surfaces},
	BOOKTITLE = {Handbook of {T}eichm\"uller theory. {V}ol. {II}},
	SERIES = {IRMA Lect. Math. Theor. Phys.},
	VOLUME = {13},
	PAGES = {611--684},
	PUBLISHER = {Eur. Math. Soc., Z\"urich},
	YEAR = {2009},
	ISBN = {978-3-03719-055-5},
	MRCLASS = {30F60 (30F45 57M50)},
	MRNUMBER = {2497777},
	MRREVIEWER = {Bruno\ P.\ Zimmermann},
	DOI = {10.4171/055-1/16},
	URL = {https://doi.org/10.4171/055-1/16},
}

@article {bgw,
	AUTHOR = {Boyer, Steven and Gordon, Cameron McA. and Watson, Liam},
	TITLE = {On {L}-spaces and left-orderable fundamental groups},
	JOURNAL = {Math. Ann.},
	FJOURNAL = {Mathematische Annalen},
	VOLUME = {356},
	YEAR = {2013},
	NUMBER = {4},
	PAGES = {1213--1245},
	ISSN = {0025-5831,1432-1807},
	MRCLASS = {57M27 (55N35)},
	MRNUMBER = {3072799},
	MRREVIEWER = {Dale\ P. O. Rolfsen},
	DOI = {10.1007/s00208-012-0852-7},
	URL = {https://doi.org/10.1007/s00208-012-0852-7},
}

@misc{CullerDunfieldPE,
	author = {Culler, Marc and Dunfield, Nathan M.},
	title = {{PE}: Peripherally parabolic and elliptic representations of 3-manifold groups},
	howpublished = {\url{https://github.com/3-manifolds/PE}},
	note = {Version 0.2.3; accessed 27 June 2026},
}

@incollection {shalen,
	AUTHOR = {Shalen, Peter B.},
	TITLE = {Representations of 3-manifold groups},
	BOOKTITLE = {Handbook of geometric topology},
	PAGES = {955--1044},
	PUBLISHER = {North-Holland, Amsterdam},
	YEAR = {2002},
	ISBN = {0-444-82432-4},
	MRCLASS = {57M07 (20E08 57M25 57M50)},
	MRNUMBER = {1886685},
	MRREVIEWER = {Kimihiko\ Motegi},
}

@article {ccgls,
	AUTHOR = {Cooper, D. and Culler, M. and Gillet, H. and Long, D. D. and
	Shalen, P. B.},
	TITLE = {Plane curves associated to character varieties of
	{$3$}-manifolds},
	JOURNAL = {Invent. Math.},
	FJOURNAL = {Inventiones Mathematicae},
	VOLUME = {118},
	YEAR = {1994},
	NUMBER = {1},
	PAGES = {47--84},
	ISSN = {0020-9910,1432-1297},
	MRCLASS = {57N10 (57M25)},
	MRNUMBER = {1288467},
	MRREVIEWER = {Serge\ L.\ Tabachnikov},
	DOI = {10.1007/BF01231526},
	URL = {https://doi.org/10.1007/BF01231526},
}

@article {chesebrotillmann,
	AUTHOR = {Chesebro, Eric and Tillmann, Stephan},
	TITLE = {Not all boundary slopes are strongly detected by the character
	variety},
	JOURNAL = {Comm. Anal. Geom.},
	FJOURNAL = {Communications in Analysis and Geometry},
	VOLUME = {15},
	YEAR = {2007},
	NUMBER = {4},
	PAGES = {695--723},
	ISSN = {1019-8385,1944-9992},
	MRCLASS = {57N10 (57M50)},
	MRNUMBER = {2395254},
	MRREVIEWER = {Thilo\ Kuessner},
	DOI = {10.4310/cag.2007.v15.n4.a2},
	URL = {https://doi.org/10.4310/cag.2007.v15.n4.a2},
}

@article {ohtsuki,
	AUTHOR = {Ohtsuki, Tomotada},
	TITLE = {Ideal points and incompressible surfaces in two-bridge knot
	complements},
	JOURNAL = {J. Math. Soc. Japan},
	FJOURNAL = {Journal of the Mathematical Society of Japan},
	VOLUME = {46},
	YEAR = {1994},
	NUMBER = {1},
	PAGES = {51--87},
	ISSN = {0025-5645,1881-1167},
	MRCLASS = {57M25 (20E08 57M07 57N10)},
	MRNUMBER = {1248091},
	MRREVIEWER = {Mark\ Brittenham},
	DOI = {10.2969/jmsj/04610051},
	URL = {https://doi.org/10.2969/jmsj/04610051},
}

@book{bochnak-coste-roy,
	author    = {Bochnak, Jacek and Coste, Michel and Roy, Marie-Fran{\c{c}}oise},
	title     = {Real Algebraic Geometry},
	series    = {Ergebnisse der Mathematik und ihrer Grenzgebiete. 3. Folge},
	volume    = {36},
	publisher = {Springer-Verlag},
	address   = {Berlin},
	year      = {1998},
	doi       = {10.1007/978-3-662-03718-8}
}

@article {ghys,
	AUTHOR = {Ghys, \'Etienne},
	TITLE = {Groups acting on the circle},
	JOURNAL = {Enseign. Math. (2)},
	FJOURNAL = {L'Enseignement Math\'ematique. Revue Internationale. 2e
	S\'erie},
	VOLUME = {47},
	YEAR = {2001},
	NUMBER = {3-4},
	PAGES = {329--407},
	ISSN = {0013-8584},
	MRCLASS = {37C85 (37E10 57S05)},
	MRNUMBER = {1876932},
	MRREVIEWER = {Grant\ Cairns},
}

@article {gordonleucke,
	AUTHOR = {Gordon, C. McA. and Luecke, J.},
	TITLE = {Reducible manifolds and {D}ehn surgery},
	JOURNAL = {Topology},
	FJOURNAL = {Topology. An International Journal of Mathematics},
	VOLUME = {35},
	YEAR = {1996},
	NUMBER = {2},
	PAGES = {385--409},
	ISSN = {0040-9383},
	MRCLASS = {57N10},
	MRNUMBER = {1380506},
	MRREVIEWER = {Darren\ D.\ Long},
	DOI = {10.1016/0040-9383(95)00016-X},
	URL = {https://doi.org/10.1016/0040-9383(95)00016-X},
}

@inproceedings{hatcher,
	title={Notes on Basic 3-Manifold Topology},
	author={Allen E. Hatcher},
	year={2001},
	url={https://api.semanticscholar.org/CorpusID:9792594}
}

@article {weil,
	AUTHOR = {Weil, Andr\'e},
	TITLE = {Remarks on the cohomology of groups},
	JOURNAL = {Ann. of Math. (2)},
	FJOURNAL = {Annals of Mathematics. Second Series},
	VOLUME = {80},
	YEAR = {1964},
	PAGES = {149--157},
	ISSN = {0003-486X},
	MRCLASS = {22.70 (32.65)},
	MRNUMBER = {169956},
	MRREVIEWER = {S.\ Murakami},
	DOI = {10.2307/1970495},
	URL = {https://doi.org/10.2307/1970495},
}

@misc{SnapPy,
	author={Culler, Marc and Dunfield, Nathan M. and Goerner,
	Matthias and Weeks, Jeffrey R.},
	title={Snap{P}y, a computer program for studying the geometry and topology of $3$-manifolds},
	howpublished={Available at \url{http://snappy.computop.org} (08/24/2026)},
}

@article {riley,
	AUTHOR = {Riley, Robert},
	TITLE = {Parabolic representations of knot groups. {I}},
	JOURNAL = {Proc. London Math. Soc. (3)},
	FJOURNAL = {Proceedings of the London Mathematical Society. Third Series},
	VOLUME = {24},
	YEAR = {1972},
	PAGES = {217--242},
	ISSN = {0024-6115,1460-244X},
	MRCLASS = {55A25},
	MRNUMBER = {300267},
	MRREVIEWER = {D.\ W. L. Sumners},
	DOI = {10.1112/plms/s3-24.2.217},
	URL = {https://doi.org/10.1112/plms/s3-24.2.217},
}

@book {hartshorne,
	AUTHOR = {Hartshorne, Robin},
	TITLE = {Algebraic geometry},
	SERIES = {Graduate Texts in Mathematics},
	VOLUME = {No. 52},
	PUBLISHER = {Springer-Verlag, New York-Heidelberg},
	YEAR = {1977},
	PAGES = {xvi+496},
	ISBN = {0-387-90244-9},
	MRCLASS = {14-01},
	MRNUMBER = {463157},
	MRREVIEWER = {Robert\ Speiser},
}

@article {menasco,
	AUTHOR = {Menasco, W.},
	TITLE = {Closed incompressible surfaces in alternating knot and link
	complements},
	JOURNAL = {Topology},
	FJOURNAL = {Topology. An International Journal of Mathematics},
	VOLUME = {23},
	YEAR = {1984},
	NUMBER = {1},
	PAGES = {37--44},
	ISSN = {0040-9383},
	MRCLASS = {57M25},
	MRNUMBER = {721450},
	MRREVIEWER = {Cameron\ McA.\ Gordon},
	DOI = {10.1016/0040-9383(84)90023-5},
	URL = {https://doi.org/10.1016/0040-9383(84)90023-5},
}

@book {morse,
	AUTHOR = {Zolkadek, Henryk},
	TITLE = {The monodromy group},
	SERIES = {Instytut Matematyczny Polskiej Akademii Nauk. Monografie
	Matematyczne (New Series) [Mathematics Institute of the Polish
	Academy of Sciences. Mathematical Monographs (New Series)]},
	VOLUME = {67},
	PUBLISHER = {Birkh\"auser Verlag, Basel},
	YEAR = {2006},
	PAGES = {xii+580},
	ISBN = {978-3-7643-7535-5; 3-7643-7535-3},
	MRCLASS = {32S40 (14D05 32S65 34M35 34M40 37F75 58K10)},
	MRNUMBER = {2216496},
	MRREVIEWER = {Jan\ Stevens},
}

@article {torus2q,
	AUTHOR = {Dubois, Jerome and Kashaev, Rinat},
	TITLE = {On the asymptotic expansion of the colored {J}ones polynomial
	for torus knots},
	JOURNAL = {Math. Ann.},
	FJOURNAL = {Mathematische Annalen},
	VOLUME = {339},
	YEAR = {2007},
	NUMBER = {4},
	PAGES = {757--782},
	ISSN = {0025-5831,1432-1807},
	MRCLASS = {57M27 (58J28 58J37)},
	MRNUMBER = {2341899},
	MRREVIEWER = {Hitoshi\ Murakami},
	DOI = {10.1007/s00208-007-0109-z},
	URL = {https://doi.org/10.1007/s00208-007-0109-z},
}

@article {twist,
	AUTHOR = {Dubois, Jerome and Huynh, Vu and Yamaguchi, Yoshikazu},
	TITLE = {Non-abelian {R}eidemeister torsion for twist knots},
	JOURNAL = {J. Knot Theory Ramifications},
	FJOURNAL = {Journal of Knot Theory and its Ramifications},
	VOLUME = {18},
	YEAR = {2009},
	NUMBER = {3},
	PAGES = {303--341},
	ISSN = {0218-2165,1793-6527},
	MRCLASS = {57M25 (57M27)},
	MRNUMBER = {2514847},
	MRREVIEWER = {Frank\ H.\ Lutz},
	DOI = {10.1142/S0218216509006951},
	URL = {https://doi.org/10.1142/S0218216509006951},
}
